\documentclass[11pt,reqno]{amsart}

\usepackage[sort&compress]{natbib}

\usepackage[OT2, T1]{fontenc}
\usepackage{url}
\usepackage{amsmath}
\usepackage{array}
\usepackage{graphicx}
\usepackage{amsfonts}
\usepackage{mathtools}
\usepackage{amssymb}
\usepackage{amstext}
\usepackage{amsthm}
\usepackage{ascii}
\usepackage{hyperref}
\usepackage{colonequals}
\usepackage{enumitem}
\usepackage[alphabetic,lite]{amsrefs}
\usepackage{cleveref}
\usepackage[all,cmtip]{xy}
\usepackage{fullpage}
\usepackage{appendix}
\usepackage{scalerel}
\usepackage{mathrsfs}
\usepackage{comment}
\usepackage{tikz-cd}
\usepackage{galois}
\usepackage{amsmath}

\numberwithin{equation}{section}

\newtheorem{theorem}{Theorem}
\newtheorem{lemma}[theorem]{Lemma}
\newtheorem{proposition}[theorem]{Proposition}
\newtheorem{corollary}[theorem]{Corollary}
\newtheorem{Question}[theorem]{Question}
\newtheorem*{thm}{Theorem}
\newtheorem{conj}[theorem]{Conjecture}

\theoremstyle{definition}

\newtheoremstyle{boldremark}
    {\dimexpr\topsep/2\relax} % space above
    {\dimexpr\topsep/2\relax} % space below
    {}          % body font
    {}          % indent amount
    {\bfseries} % theorem head font
    {.}         % punctuation after theorem head
    {.5em}      % space after theorem head
    {}          % theorem hed spec. (empty = "normal")

\theoremstyle{boldremark}
\newtheorem{definition}{Definition}
\newtheorem{remark}[theorem]{Remark}

\newtheorem*{claim}{Claim}
\newtheorem{example}[theorem]{Example}

\newcommand{\bC}{\mathbb{C}}
\newcommand{\bD}{\mathbb{D}}

\newcommand{\bF}{\mathbb{F}}
\newcommand{\bG}{\mathbb{G}}
\newcommand{\bH}{\mathbb{H}}

\newcommand{\bK}{\mathbb{K}}

\newcommand{\bN}{\mathbb{N}}

\newcommand{\bQ}{\mathbb{Q}}

\newcommand{\bZ}{\mathbb{Z}}

\newcommand{\bbD}{\mathbf{D}}

\newcommand{\bbL}{\mathbf{L}}

\newcommand{\bbQ}{\mathbf{Q}}

\newcommand{\cH}{\mathcal{H}}

\newcommand{\cJ}{\mathcal{J}}

\newcommand{\cL}{\mathcal{L}}
\newcommand{\cM}{\mathcal{M}}
\newcommand{\cN}{\mathcal{N}}
\newcommand{\cO}{\mathcal{O}}
\newcommand{\cP}{\mathcal{P}}

\newcommand{\cZ}{\mathcal{Z}}

\newcommand{\fm}{\mathfrak{m}}

\newcommand{\et}{{\text{\'et}}}

\newcommand{\cris}{{\mathrm{cris}}}

\DeclareMathOperator{\GSp}{GSp}

\DeclareMathOperator{\sO}{SO}
\DeclareMathOperator{\GSpin}{GSpin}

\DeclareMathOperator{\Spf}{Spf}

\DeclareMathOperator{\End}{End}
\DeclareMathOperator{\Hom}{Hom}

\DeclareMathOperator{\Spec}{Spec}
\DeclareMathOperator{\Ext}{Ext}

\DeclareMathOperator{\Span}{Span}

\DeclareMathOperator{\Fil}{Fil}

\DeclareMathOperator{\Id}{Id}

\DeclareMathOperator{\Tr}{Tr}
\DeclareMathOperator{\Frac}{Frac}

\DeclareMathOperator{\rk}{rk}
\DeclareMathOperator{\tor}{tor}

\DeclareMathOperator{\Frob}{Frob}
\DeclareMathOperator{\gr}{gr}

\DeclareMathOperator{\AO}{NO}
\DeclareMathOperator{\SSS}{SS}
\DeclareMathOperator{\even}{even}
\DeclareMathOperator{\odd}{odd}

\DeclareMathOperator{\kk}{\mathbf{k}}
\DeclareMathOperator{\Pic}{Pic}
\DeclareMathOperator{\Cl}{Cl}
\DeclareMathOperator{\KS}{KS}

\begin{document}

\begin{center}
       \vspace*{1cm}
\Large\textbf{{Splitting of almost ordinary abelian surfaces in families and the $S$-integrality conjectures}}
        \vspace*{0.5cm}
\begin{center}
\normalsize{Ruofan Jiang}  
\end{center}
\vspace*{0.5cm}
\normalsize\end{center}
\textbf{Abstract}. Let $A$ be a non-isotrivial almost ordinary abelian surface with possibly bad reductions over a global function field of odd characteristic $p$.
Suppose $\Delta$ is an infinite set of positive integers, such that $\left(\frac{m}{p}\right)=1$ for $\forall m\in \Delta$. If $A$ does not admit any global real multiplication, we prove the existence of infinitely many places modulo which the reduction of $A$ has endomorphism ring containing $\mathbb{Z}[x]/(x^2-m)$ for some $m\in \Delta$. This implies that there are infinitely many places modulo which $A$ is not simple, generalizing the main result of \cite{MAT} to the non-ordinary case. As an another application, we also  generalize the $S$-integrality theorem for elliptic curves over number fields, as proved in \cite{BR08}, to the setting of abelian surfaces over global function fields.\\[5pt]
\noindent\textit{2010 Mathematics Subject Classification}: 11G10, 14G25.

\address{Dept.\ of Mathematics, University of California--Berkeley}
\email{ruofanjiang@berkeley.edu}
\date{\today}
%\maketitle
\numberwithin{theorem}{section}

\tableofcontents
\section{Introduction}
Let $F$ be a global function field of odd characteristic $p$ and $A/F$ be a non-isotrivial  abelian surface with possibly bad reductions. Assume that $A$ is generically almost ordinary, i.e., the kernel of the multiplication by $p$ map $[p]:A\rightarrow A$ consists of $p$ $\overline{F}$-points. Equivalently, for all but finitely many places of $F$, the reductions of $A$ are almost ordinary, i.e.,  the kernels of $[p]$ on the reductions consist of $p$ $\overline{\bF}_p$-points. In this paper, we prove:
\begin{theorem}[Theorem~\ref{T:interMain}]\label{TT:T2}
Let $\Delta$ be an infinite set of positive integers such that $\left(\frac{m}{p}\right)=1$ for $\forall m\in \Delta$. If $A$ does not admit any global real multiplication, then there are infinitely many places modulo which the endomorphism ring of the reduction of $A$ contains $\mathbb{Z}[x]/(x^2-m)$ for some $m\in \Delta$. 
\end{theorem}
The theorem unifies and generalizes two previous results from different fields. The first is a result pertaining to the simpleness and real multiplications of the reductions of generically ordinary abelian surfaces over global function fields, see \cite{MAT}. The second is the $S$-integrality theorem for non-torsion points over elliptic curves, as proved in \cite{BR08}. 
%\subsection{Previous works, motivations and generalizations} We begin by a brief survey of previous works that are involved. In the meantime, we make various generalizations and point out how our main theorem fits in the story.
\subsection{Real multiplication and splitting of reductions}
The question of the splitting of reduction of abelian varieties over a global fields has a long history. In \cite{MP08}, Murty and Patankar conjectured that an absolutely simple abelian variety over a number field has absolutely simple
reduction for a density one set of primes up to a finite extension if and only if its endomorphism ring is
commutative. The conjecture was solved for abelian varieties of dimension 2 and 6 with geometric endomorphism ring $\mathbb{Z}$ (\cite{Cha97}), and established in full generality by Zywina (\cite{Zyw14}) conditionally on the Mumford--Tate conjecture. Though the set of primes over which a given abelian variety does not have  absolutely simple reduction is conjecturally of density 0, it remains an interesting question whether such primes form an infinite set.   

In \cite{ST20}, this question is answered affirmatively for abelian surfaces with real multiplication. More precisely, it is proved that an abelian surface over a number field $E$ with real multiplication has non-absolutely-simple reductions over infinitely many places of $E$. Later on, a similar question is raised and proved for abelian surfaces over global function field of characteristic $p\geq 5$ in \cite{MAT}. More precisely, the authors proved that a non-isotrivial generically ordinary abelian surface $A$ over a global function field $F$ with good reduction everywhere admits non-simple reductions over infinitely many places of $F$. The condition that $A$ has good reduction everywhere is further removed in \cite{ST22}. The only case that remains unsolved is for non-isotrivial abelian surfaces over global function field which are generically almost ordinary. However, when $A$ is generically almost ordinary, it is not always true that $A$ admits non-simple reductions over infinitely many places of $F$. We have the following counterexample:
\begin{example}\label{ex:counterex}
  Let $m$ be a positive integer such that $\left(\frac{m}{p}\right)=1$ but $m$ is not a perfect square. 
  Consider a Hilbert modular surface associated to the real quadratic field $\mathbb{Q}(\sqrt{m})$.
  The almost ordinary locus $U$ of the  Hilbert modular surface is one dimensional. We claim that, this gives rise to an almost ordinary abelian surface over a global function field whose reduction is simple for all but finitely many places. Indeed, let $P$ be any nonsupersingular $\overline{\bF}_p$-point on $U$ and let $\mathscr{A}_P$ be the reduction. If $\mathscr{A}_P$ splits, then it is isogenous to a product of an ordinary elliptic curve and a supersingular elliptic curve. The algebra $\End^0(\mathscr{A}_P)$ is then a product of a complex quadratic extension of $\bQ$ and a quaternion nonsplit at $p$. 
But then $\bQ(\sqrt{m})\subsetneq\End^0(\mathscr{A}_P)$, a contradiction. Therefore $\mathscr{A}_P$ must be simple. See Corollary~\ref{T:simple} for another argument.  
\end{example}
In this paper, we will investigate the question of exceptional splitting of generically almost ordinary abelian surfaces over global function fields. We will show that Example~\ref{ex:counterex} is the only obstruction towards the existence of infinitely many places of non-simple reductions. The following is a consequence of Theorem \ref{TT:T2}:
%\subsection{Reduction of abelian surfaces over function fields} We begin by stating the main consequence of the paper: 
\begin{theorem}[Theorem~\ref{T:maintheorem}]\label{TT:1}
 Let $A$ be a non-isotrivial generically almost ordinary abelian surface with possibly bad reductions over a global function field of odd characteristic $p$. Let $l$ be a positive integer such that $\left(\frac{l}{p}\right)=1$, and suppose that $A$ does not admit any global real multiplication, then
 \begin{enumerate}
     \item there are infinitely many places modulo which $A$ admits real multiplication by $\mathbb{Q}[x]/(x^2-l)$.
     \item  there exist infinitely places modulo which $A$ is not simple.
 \end{enumerate}
\end{theorem}
There is a geometric interpretation of Theorem \ref{TT:1}. Let $\mathcal{A}_{2,\overline{\mathbb{F}}_p}$ be the Siegel moduli stack of principally polarized abelian surfaces over $\overline{\mathbb{F}}_p$.  The set of special divisors over $\mathcal{A}_{2,\overline{\mathbb{F}}_p}$ is parameterized by positive integers $m$. In fact, for each $m$ there is a  special divisor $Z(m)$ parametrizing abelian surfaces whose endomorphism ring contains
$\mathbb{Z}[x]/(x^2-m)$. Let $m=lk^2, k\in \mathbb{Z}^+$, if a point lies on $Z(m)$ then it corresponds to an abelian surface with real multiplication by $\mathbb{Q}[x]/(x^2-l)$. In particular, when $l=1$, an $\overline{\mathbb{F}}_p$-point of $Z(m)$ corresponds to a non-simple Abelain surface.
Therefore Theorem \ref{TT:1} is a consequence of the following: for any positive integer $l$ with $\left(\frac{l}{p}\right)=1$, and any nonconstant morphism  $U\rightarrow\mathcal{A}_{2,\overline{\mathbb{F}}_p}$ from a smooth connected (affine) curve $U$ such that the image lies generically in the almost ordinary stratum but not on any special divisor, there are infinitely many $\overline{\mathbb{F}}_p$-points of $U$ that lie on a special divisor of form $Z(lk^2), k\in\mathbb{Z}$.

It is worth noting that, as a variant of the problem of exceptional splitting for the reductions, similar Picard rank jump problems for reductions were solved for generically ordinary K3 surfaces over a global function field with everywhere good reduction in \cite{MST22}. The condition of having everywhere good reduction is further removed in  \cite{ST22}. In the future work we will use the techniques developed in this paper to solve these questions for K3 surfaces lying generically in higher Newton strata. 
\subsection{$S$-integrality conjectures}
The $S$-integrality theorem for elliptic curves is a classical result proved in \cite{BR08}. Let $k$ be a number field, $S$ be a finite set of places of $k$ containing all Archimedean places, and $O_{k,S}$ be the ring of $S$-integers of $k$. Let $E/k$ be an elliptic curve and $\mathcal{E}_S$ be a model of $E$ over $O_{k,S}$. Then the main theorem of \textit{loc.cit} says that if $\alpha\in E(\overline{k})$ is non-torsion, then there are only finitely many torsion points of $E(\overline{k})$ whose closures in $\mathcal{E}_S$ are disjoint from the closure of $\alpha$. The property that a point $\xi\in E(\overline{k})$ having its closure in $\mathcal{E}_S$ disjoint from the closure of $\alpha$, is called "$S$-integral to $\alpha$" by the authors, hence the name. As a natural generalization, one can consider a similar question for abelian varieties. The precise statement conjecture was raised
by Su-ion Ih:\begin{conj}[$S$-integrality conjecture for abelian varieties, see Conj 1.2 of \textit{loc.cit}]\label{TT:3}
Pick an effective divisor $D$ of $A$, defined over $\overline{k}$, at least one of its components is not a translation of an abelian subvariety by a torsion point. Let $\overline{D}$ be the closure of $D$ in $\mathcal{A}_S$. The set consisting of all torsion points of $A(\overline{k})$ whose closures in $\mathcal{A}_S$ are disjoint from $\overline{D}$, is not Zariski dense in $A$. 
\end{conj}
It is easy to see that \textit{loc.cit} solves Conjecture~\ref{TT:3} for the case of elliptic curves. To generalize the problem, one can replace $k$ by a global function field $F$, or replace the intersections of the closures of divisors and torsion points by the intersections of the closures of cycles with complementary absolute dimensions. 

However, there is even a more mysterious generalization that one can make by exploiting the analogy between abelian varieties and Shimura varieties. This analogy is fruitful in generalizing unlikely intersection results from abelian varieties to Shimura varieties, e.g.\ , Mordell--Lang conjecture and André--Oort conjecture. Under this analogy, torsion points correspond to CM points,  translations of abelian subvarieties by torsion points correspond to special cycles. After the main result of \cite{BR08}, similar questions for Shimura varieties are raised. Let $\mathcal{S}_S$ be an integral model of a Shimura variety over $O_{k,S}$: %(\cite{KM09})
\begin{Question}[$S$-integrality questions for Shimura varieties]\label{Qes:1}
Let $\Gamma\subseteq \mathcal{S}_{\overline{k}}$ be a subvariety, which is not a special cycle. What can one say about the scarcity of the special cycles $\delta\subseteq \mathcal{S}_{\overline{k}}$ such that (1) the closures $\overline{\Gamma}$ and $\overline{\delta}$ in $\mathcal{S}_S$ are of complementary absolute dimensions and (2) the intersection $\overline{\Gamma}\cap \overline{\delta}$ is empty? In particular, one can ask:
\begin{enumerate}
    \item Suppose that $\Gamma$ is a non-special divisor. Let $T_\Gamma$ be the set of CM points in $\mathcal{S}(\overline{k})$ whose closures in $\mathcal{S}_S$ are disjoint from $\overline{\Gamma}$. Is it true that $T_\Gamma$ is  not Zariski dense in $\mathcal{S}_{\overline{k}}$?
    \item \label{q12}Suppose that $\Gamma$ is a non-CM point. Let $T_\Gamma$ be the set of special divisors whose closures in $\mathcal{S}_S$ are disjoint from $\overline{\Gamma}$. Is $T_\Gamma$ finite ? 
\end{enumerate}
\end{Question}
One can also formulate a function field analogue of Question~\ref{Qes:1}. Take $\mathcal{S}_{\overline{\mathbb{F}}_p}$ as the mod $p$ special fiber of the canonical integral model of a Shimura variety. Let $\tau:U\rightarrow\mathcal{S}_{\overline{\mathbb{F}}_p}$ be a non-constant morphism from a smooth connected curve to $\mathcal{S}_{\overline{\mathbb{F}}_p}$ and $S$ be a finite subset of $U(\overline{\mathbb{F}}_p)$. The curve $U$ is an analogue of the ring of $S$-integers $O_{k,S}$, and $\tau$ is an analogue of the closure of a point $\xi\in \mathcal{S}(\overline{k})$ in the canonical integral model of a Shimura variety over $O_{k,S}$. We can then state an analogue of Question~\ref{Qes:1}(\ref{q12}): 
\begin{Question}[$S$-integrality question for mod $p$ Shimura varieties]\label{Q:1}
With notation as above, suppose that the image of $U$ does not lie on any special subvariety, are there only finitely many special divisors disjoint from the image of $U-S$? 
\end{Question}
In this paper, we investigate Question~\ref{Q:1} when $\mathcal{S}_{\overline{\mathbb{F}}_p}=\mathcal{A}_{2,\overline{\mathbb{F}}_p}$. In this case, a special divisor is 
of form $Z(m)$, which parametrizes abelian surfaces whose endomorphism rings contain $\mathbb{Z}[x]/(x^2-m)$. The answer to Question~\ref{Q:1} is \textbf{yes} when (a) the image of $U$ is generically almost ordinary, and (b) the intersections with the image of $U-S$ are taken among certain subclass of special divisors. As we will see, restricting to a subclass of special divisors as in (b) is necessary. The following result is a  consequence of Theorem \ref{TT:T2}:
\begin{theorem}[Theorem \ref{T:unlikelyint}]
 \label{QQ:1}
Suppose that $\mathcal{S}_{\overline{\mathbb{F}}_p}\simeq \mathcal{A}_{2,\overline{\mathbb{F}}_p}$ and that the image of $U$ lies generically on the almost ordinary locus. Let $S\subseteq U(\overline{\mathbb{F}}_p)$ be a finite set containing all the points with images in the supersingular loci. Then following are true: 
\begin{enumerate}
    \item\label{QQ:11}  All special divisors $Z(m)$ such that $\left(\frac{m}{p}\right)=-1$ are disjoint from $U-S$.
    \item \label{QQ:12} Only finitely many special divisors $Z(m)$ such that $\left(\frac{m}{p}\right)=1$ are disjoint from $U-S$.
\end{enumerate}
\end{theorem}
Therefore, Question~\ref{Q:1} is true, if we make some mild assumptions on $S$, and make the further restriction that the numbers $m$ parametrizing the special divisors $Z(m)$ enjoy the property $\left(\frac{m}{p}\right)=1$. On the other hand, Theorem~\ref{QQ:1}(\ref{QQ:11}) shows that Question~\ref{Q:1} is in general false. However, the phenomenon is also interesting, as it is "unlikely" that infinitely many special divisors all intersect the image of $U$ at a fixed set of points.

\subsection{Outline of the proof} One can obtain Theorem~\ref{TT:T2} by arithmetic 
intersection theory over the Siegel moduli stack $\mathcal{A}_{2,\overline{\mathbb{F}}_p}$. Let $p$ be an odd prime and $\mathcal{A}_{2,\overline{\mathbb{F}}_p}$ be the moduli stack of principally polarized abelian surfaces over the algebraic closure of ${\mathbb{F}}_p$, regarded as a mod $p$ GSpin Shimura variety. A non-isotrivial almost ordinary abelian surface over a global function field can then be viewed as a smooth connected affine curve $U$ which maps nonconstantly to $\mathcal{A}_{2,\overline{\mathbb{F}}_p}$ with generically almost ordinary image. As noted earlier, for each $m$, there is a special divisor $Z(m)$ parametrizing abelian surfaces whose endomorphism rings contain
$\mathbb{Z}[x]/(x^2-m)$. We now restate Theorem~\ref{TT:T2} into the following intersection theoretical result: 
\begin{theorem}\label{TT:2}
 Let $U$ be a smooth connected (affine) curve with a nonconstant map $U\rightarrow \mathcal{A}_{2,\overline{\mathbb{F}}_p}$ sending the generic point to the almost ordinary stratum, whose image does not lie in any special divisor. Let $\Delta$ be an infinite set of positive integers $m$ such that  $\left(\frac{m}{p}\right)=1$, there are infinitely many $\overline{\mathbb{F}}_p$-points of $U$ lying on $\bigcup_{m\in \Delta} Z(m)$.
\end{theorem}

The overall strategy for proving Theorem~\ref{TT:2} is a generalization of that in \cite{MAT} with several additional difficulties. Fix a toroidal compatification $\mathcal{A}^{\text{tor}}_{2,\overline{\mathbb{F}}_p}$. The closure of $Z(m)$ in $\mathcal{A}^{\text{tor}}_{2,\overline{\mathbb{F}}_p}$ will also be called special divisors. By abuse of notation, they are also denoted $Z(m)$. Let $C$ be the smooth compatification of $U$. The morphism of $U\rightarrow \mathcal{A}_{2,\overline{\mathbb{F}}_p}$ can be uniquely extended to a morphism $C\rightarrow\mathcal{A}^{\tor}_{2,\overline{\mathbb{F}}_p}$. We then study the intersection of $C$ with $Z(m)$ in 
$\mathcal{A}^{\text{tor}}_{2,\overline{\mathbb{F}}_p}$. 

On one hand, $C\cdot Z(m)$ can be understood globally via  arithmetic Borcherds' theory for GSpin Shimura varieties, which relates the intersection numbers $C\cdot Z(m)$ to coefficients of certain Eisenstein series. 

On the other hand, the intersection $C\cdot Z(m)$ can be understood locally at each intersection point: The local intersection at a point is essentially the number of deformations of a formal special endomorphism of the $p$-divisible group corresponding to that point. More precisely, let $P$ be a point on $U$, the local intersection at $P$ can be expressed in terms of the number of deformations of special endomorphisms of $\mathscr{A}_P[p^\infty]$ to finite length quotient rings of $U^{/P}$ (the formal germ of $U$ at $P$). It turns out that the computation of the intersection number boils down, via a deformation result of $p$-divisible groups developed by Faltings--Kisin--Moonen (see \S\ref{3.1} for more details), to a lattice point counting problem in the quadratic lattice of special endomorphisms.

Let $\Delta$ be an infinite set of $m$ such that $\left(\frac{m}{p}\right)=1$. We prove that,\begin{enumerate}
    \item \label{171}  when $m\rightarrow \infty$, the local contribution of intersection numbers at supersingular  points is asymptotically much smaller than the global intersection $C\cdot Z(m)$.
    \item\label{172} the local contribution of intersection numbers at a non-supersingular or bad reduction  point is bounded above by an integer independent of $m$.
\end{enumerate}
Given (\ref{171}) and (\ref{172}), one is able to deduce the unboundeness of $\# ({C}\cap {Z}(m))$ as $m\rightarrow \infty$, which implies Theorem~\ref{TT:2}.

Unlike the generically ordinary case tackled in \cite{MAT}, (\ref{171}) is not always true in the generically almost ordinary setting. For example, (\ref{171}) is not true when $C$ lies on a Hilbert modular surface as in Example~\ref{ex:counterex}. In fact, the obstruction for (\ref{171}) being true is the existence of a formal special endomorphism of the pullback $p$-divisible group over $C$ in the formal neighbourhood of a superspecial point. In other words, the existence of a formal special endomorphism implies that the local and global intersections are asymptotically equal (see Remark~\ref{lastrmk}). However, our second result (Theorem~\ref{TT:4}) claims that the existence of a formal special endomorphism would guarantee that $C$ lies on a special divisor $Z(m)$, which is exactly the case excluded in Theorem~\ref{TT:2}. Furthermore, the independence of $m$ claim we made in (\ref{172}) is much stronger than that in \cite{MAT}. In fact, we again need to use the algebraization result (Theorem~\ref{TT:4}) to achieve this, see Proposition~\ref{T:intindep}.

As a final remark, we note that the deformation problem of special endomorphisms in the  
characteristic p setting is much more subtle than its characteristic 0 counterpart. Let $k$ be an algebraically closed field, $U/k$ be a connected smooth affine curve, and $P$ be a closed point of $U$. Consider an abelian surface $\mathscr{A}/U$ such that $\mathscr{A}[p^\infty]$ does not admit any global special endomorphism. If $\mathrm{char} k=0$, then Hodge theory implies that there is a sharp cutoff, i.e., a uniform bound $N$, such that no special endomorphism of $\mathscr{A}_{P}[p^{\infty}]$ can deform to $\Spf k[[t]]/t^{N+1}$, where $t$ is a uniformizer of $U^{/P}$. However, when $\mathrm{char} k>0$, there is no such uniform bound, and it is much harder to bound the number of deformations. In the worst case, one may have a special endomorphism $\omega$ that deforms to $\Spf k[[t]]$. When $\mathscr{A}$ is generically almost ordinary, our algebraization theorem (Theorem~\ref{TT:4}) would guarantee that this does not happen. In general, it is not at all clear whether this happens or not, see Question~\ref{Qes:tt}. 

%Another difficulty in our case lies in bounding the local intersections. The behavior is more subtle than the generic ordinary case. This is the reason why we need a more explicit description of the decay.  

\subsection{Algebraization of formal special endomorphisms}
Let $U$ be a connected smooth (not necessarily projective) curve with a nonconstant morphism $U\rightarrow \mathcal{A}_{2,\overline{\mathbb{F}}_p}$. Let $K(U)$ be the function field of $U$, and $\widehat{K}(U)$ be its completion at a place. Let $\mathscr{A}[p^\infty]$ be the universal $p$-divisible group over $\mathcal{A}_{2}$ and let $\tilde{\omega}$ be a nonzero special endomoprhism of $\mathscr{A}_{\Spec \widehat{K}(U)}[p^\infty]$. We show that:
\begin{theorem}[Theorem~\ref{T: algebraicC}]\label{TT:4}
Notation as above. Suppose the image of $U$ is generically almost ordinary, then the image of $U$ lies on a special divisor.
\end{theorem}
Though we just state the theorem for a curve $U$ mapping to almost ordinary strata of $\mathcal{A}_{2,\overline{\mathbb{F}}_p}$, our method in Theorem~\ref{T: algebraicC} works in the case where $U$ is an arbitrary smooth variety, and the almost ordinary stratum of $\mathcal{A}_{2,\overline{\mathbb{F}}_p}$ is replaced by the almost supersingular stratum (the one with height ${1}/{n}$) of mod $p$ GSpin Shimura varieties, see Remark~\ref{remarkOth}. Inspired by this, we ask a more general algebraization question: 
\begin{Question}\label{Qes:tt}
Let $X$ be a smooth subvariety of a mod p Shimura variety $\mathcal{S}_{\overline{\mathbb{F}}_p}$ lying generically on an Newton strata $\mathcal{S}_{\overline{\mathbb{F}}_p}^\eta$,
$P\in X$ and $\omega$ be an endomorphism of $\mathscr{A}[p^{\infty}]$. If $X^{/P}$ coincides with the support of $\mathcal{D}(\omega)\cap \mathcal{S}_{\overline{\mathbb{F}}_p}^\eta$,  does there exist a shimura subvariety $ \mathcal{T}\subset \mathcal{S}$ such that $X$ is the support of a irreducible component of $\mathcal{T}_{\overline{\mathbb{F}}_p}^\eta$ ?
\end{Question}
The question can be considered as a generalization of Chai's Tate linear conjecture (see \cite[Conjecture 7.2, Remark 7.2.1, Proposition 5.3, Remark 5.3.1]{Cha03}).
When $P$ is ordinary, for example, it will be answered affirmatively if the Tate linear conjecture is true. Theorem~\ref{TT:4} and its fore-mentioned generalization answers the question in the case where $\mathcal{S}$ is the GSpin Shimura variety associated to $\GSpin(2, 2n-1)$, $\eta$ is the height $1/n$ strata, and $\omega$ is a special endomorphism.

\subsection{Organization of the paper} In \S\ref{cha2}, we review the definition of $\mathcal{A}_2$ as a GSpin Shimura variety and review various background theories. In \S\ref{compactification} we review the theory of toroidal compactifications and logarithmic geometry. In \S\ref{sec:AIT}, we review arithmetic intersection theory. In \S\ref{algebraization} we prove Theorem~\ref{TT:4}. In \S\ref{cha3}, we prove the decay results. In \S\ref{reduction} we prove Theorem~\ref{TT:2} and the other main theorems. 
\subsection{Notation}\label{subsec:Notation} We use $p$ to denote a odd prime, use $\mathbf{k}$ to denote $\overline{\mathbb{F}}_p$ and $W$ its ring of Witt vectors. We write $K=W[\frac{1}{p}]$. For a fraction ideal $\mathfrak{a}\subseteq W[[x,y,z]]$, and an element $f\in K[[x,y,z]]$ we use the notation $f=O(\mathfrak{a})$ to indicate $f\in \mathfrak{a}$. In particular, $f\in O(1)$ means $f\in W[[x,y,z]]$.

We will use $\sigma$ to denote the semilinear endmorphism on $K[[x,y,z]]$ which is the Frobnious on $K$ and sends $x,y,z$ to $x^p,y^p,z^p$. For $i\geq 0$ and matrix $M$ with entries in $K[[x,y,z]]$, we use $M^{[i]}$ to denote $\sigma^i(M)$.\\[10pt]
\textbf{Acknowledgements}  The author thanks Ananth Shankar for suggesting this problem and his help, thanks Asvin G, Qiao He, Jiaqi Hou and Salim Tayou for valuable discussions, and thanks Jordan Ellenberg for pointing out some imprecision in the introduction. The author also thanks Keerthi Madapusi Pera for answering a question on toroidal compactifications and Martin Olsson for answering a question on log geometry.  The author is partially supported by the NSF grant DMS-2100436. 

\section{Preliminaries}\label{cha2}
We review several constructions. In \S\ref{2.1} we recall the construction of $\mathcal{A}_2$ as a $\GSpin$ Shimura variety. In \S\ref{compactification} we review the theory of toroidal compatifications. In section \S\ref{3.1} we review the deformation theory of the universal $p$-divisible group over $\mathcal{A}_{2,\kk}$ at a non-necessary ordinary point $P$. In \S\ref{C:intersection} we review the intersection theory over $\mathcal{A}^\text{tor}_{2,\kk}$ and the modularity of special divisors.  

\subsection{$\mathcal{A}_2$ as GSpin Shimura variety}\label{2.1}
We begin by reviewing the definition and basic properties of GSpin Shimura varieties, the main references are \cite{KR00}, \cite{MP15} and \cite{AGHMP17}. Let $p\geq 3$ be a prime. For an integer $b\geq 1$, let $(L,Q)$ be a quadratic $\mathbb{Z}$-lattice of rank $b+2$ with an even bilinear form $(,):L\otimes L\rightarrow \mathbb{Z}$
such that for $x\in L$, $Q(x)=(x,x)/2\in \mathbb{Z}$, and that $(L,Q)$ is self dual at $p$. Let $G:=\GSpin(L\otimes \mathbb{Z}_{(p)},Q)$ be the group of spinor similitude of $L\otimes \mathbb{Z}_{(p)}$. $G$ is a reductive group over $\mathbb{Z}_{(p)}$, and is a subgroup of $\Cl(L\otimes \mathbb{Z}_{(p)})^{\times}$, where $\Cl(\cdot)$ is the Clifford algebra.  The group $G(\mathbb{R})$ acts on the symmetric space 
$$\mathcal{D}_L=\{z\in\mathbb{P}(L_\mathbb{C})|(z,z)=0,(z,\overline{z})<0\}$$
via the natural morphism $G_\mathbb{Q}\rightarrow \sO(L\otimes \mathbb{Q})$. This gives rise to a GSpin Shimura datum $(G_\mathbb{Q},\mathcal{D})$ with refex field $\mathbb{Q}$. Let $\mathbb{K}\subseteq G(\mathbb{A}_f)\cap \Cl(L\otimes \widehat{\mathbb{Z}})^\times$ be a compact open subgroup such that $\mathbb{K}_p=G(\mathbb{Z}_p)$. Then we have a Deligne--Mumford stack $\mathcal{S}h:=\mathcal{S}h(G_\mathbb{Q},\mathcal{D}_L)_\mathbb{K}$ over $\mathbb{Q}$, such that $\mathcal{S}h(\mathbb{C})=G(\mathbb{Q})\backslash \mathcal{D}_L\times G(\mathbb{A}_f)/ \mathbb{K}$. $\mathcal{S}h$ is called the GSpin Shimura variety associated to $G_{\bbQ}$ with level structure $\bK$. By \cite[Thm 2.3.8]{KM09}), $\mathcal{S}h$ admits a smooth canonical integral model $\mathcal{S}_\mathbb{K}$ over $\mathbb{Z}_{(p)}$.

Let $H=\Cl(L)$ with the action of itself on the right. Equip $\Cl(L\otimes \mathbb{Z}_p)$ with the action of ${G}$ on the left. There exists a choice of symplectic form on $H$ that gives rise to a map $G_\mathbb{Q}\rightarrow \GSp(H_\mathbb{Q})$, which induces a morphism of Shimura data, hence a morphism of Shimura varieties and their integral models.  
Pulling back the universal abelian scheme over the Siegel modular variety yields a Kuga--Satake abelian scheme $\mathcal{A}^{\KS}\rightarrow \mathcal{S}$ with left $\Cl(L)$-action and $\mathbb{Z}/2\mathbb{Z}$-grading $\mathcal{A}^{\KS}\simeq \mathcal{A}^{\KS+}\times\mathcal{A}^{\KS-}$. Let $\mathbf{H}_\text{B},\mathbf{H}_{\text{dR}}, \mathbf{H}_{l,\text{ét}}$ be the integral betti, de Rham, $l$-adic étale ($l\neq p$) relative first cohomology of $\mathcal{A}^{\KS}\rightarrow \mathcal{S}$ and $\mathbf{H}_{\text{cris}}$ be the first integral crystalline cohomology of $\mathcal{A}_{\mathbb{F}_p}^{\KS}\rightarrow \mathcal{S}_{\mathbb{F}_p}$. 

The action of $L$ on $H$ induces a $G_\mathbb{Q}$ invariant embedding $L_\mathbb{Q}\hookrightarrow \End_{\Cl(L)}(H_\mathbb{Q})$. Correspondingly one has a sheaf $\mathbf{L}_{\bullet}$ where $ \bullet=\text{B},\text{dR},\{l,\text{ét}\}, {\text{cris}}$, equipped with a natural quadratic form $\mathbf{Q}$ such that $f\comp f=\mathbf{Q}(f)\Id$ for a section $f$ of $\mathbf{L}_\bullet$, and embeddings $\mathbf{L}_\bullet\hookrightarrow \End_{\Cl(L)}(\mathbf{H}_{\bullet})$, which are compatible with various $p$-adic Hodge theoretic comparison maps, see \cite[\S 4.3]{AGHMP17} and \cite[Proposition 3.11, 3.12, 4.7]{MP15}.
\subsubsection{Siegel moduli of abelian surfaces}\label{subsubsec:SA}
Now we briefly review the construction of $\mathcal{A}_2$ as a GSpin Shimura variety. Let $b=3$ and $Q=x_0^2+x_1x_2-x_3x_4$. By the general defintion of GSpin Shimura varieties, the resulting Shimura variety $\mathcal{S}$ admits a Kuga--Satake abelian scheme $\mathcal{A}^{\KS}\rightarrow \mathcal{S}$ of dimension 16. According to \cite[Remark 2.2.2]{MAT}, the abelian subscheme $\mathcal{A}^{\KS+}\simeq \mathscr{A}^{4}$ for an abelian scheme $\mathscr{A}\rightarrow \mathcal{S}$ of relative dimension 2, which makes $\mathcal{S}$ the Siegel modular variety  $\mathcal{A}_{2,(p)}=\mathcal{A}_2\times \mathbb{Z}_{(p)}$ with universal abelian scheme $\mathscr{A}$. Let $\dagger$ be the Rosati involution on $\mathscr{A}$. 
Let $\mathbf{H}'_\bullet, \bullet=\text{B},\text{dR},\{l,\text{ét}\}, {\text{cris}}$ be the {}{sheaf} arising from first relative $\bullet$-cohomology of $\mathscr{A}$. Then there is an embedding $\mathbf{L}_\bullet\hookrightarrow \End(\mathbf{H}'_{\bullet})$, which realizes $\mathbf{L}_\bullet$ as the 
$\dagger$-fixed sub-sheaf with trace 0. {}{Here the trace is the reduced trace on $\End(\mathscr{A})_{\bQ}$ (cf. \cite[Definition 2.2.1]{MAT}).} The Shimura variety $\mathcal{A}_{2,(p)}$ is equipped with a class of special divisors:\\

\begin{definition}\label{ddd1} Let $S$ be a $\mathcal{A}_{2,(p)}$-scheme, and $P\in \mathcal{A}_{2,(p)}(\kk)$: \begin{enumerate}
    \item  A special endomorphism of $\mathscr{A}_S$ is an element $s\in \End(\mathscr{A}_S)$ such that $s^{\dagger}=s$ and $\Tr s=0$,
    \item 
    If $S$ is an $\mathcal{A}_{2,\kk}$-scheme. Then a special endomorphism of $\mathscr{A}_S[p^\infty]$ is an element of $\End(\mathscr{A}_S[p^\infty])$ whose crystalline realization lies in $\mathbf{L}_{\cris,S}$ {}{(cf. \cite[Definition 2.2.4]{MAT})},
    \item A formal special endomoprhism of $\mathscr{A}[p^\infty]$ in the formal neighborhood of $P$ is an endomorphism of $\End(\mathscr{A}_{\mathcal{A}_{2,\overline{\mathbb{F}}_p}^{/P}}[p^\infty])$ whose crystalline realization lies in $\mathbf{L}_{\cris,\mathcal{A}_{2,\overline{\mathbb{F}}_p}^{/P}}$,
    \item {}{Following \cite[Definition 2.2.6]{MAT}}, for an integer $m>0$, the special divisor $\mathcal{Z}(m)$ is a Deligne--Mumford stack which is étale locally an effective Cartier divisor  over 
    $\mathcal{A}_{2,(p)}$. It is defined by the functor  $$\mathcal{Z}(m)(S)=\{s\in \End(\mathscr{A}_S)\,\, \text{special }|s\comp s=m\}$$
for any $\mathcal{A}_{2,(p)}$-scheme $S$.  By \cite[Proposition 4.5.8]{AGHMP17}, $\mathcal{Z}(m)$ is flat over $\mathbb{Z}_{(p)}$, hence $\mathcal{Z}(m)_{\bF_p}$ is a divisor of $\mathcal{A}_{2,\bF_p}$. We will use $Z(m)$ to denote $\mathcal{Z}(m)_{\kk}$. {}{According to the paragraph after \cite[Proposition 6.5.2]{HBK17}, we can ``take the image\footnote{This is the terminology used by the authors in \cite{MAT}}'' of $\mathcal{Z}(m)$ in $\mathcal{A}_{2,(p)}$,  viewing it as a Cartier divisor on $\mathcal{A}_{2,(p)}$ in the usual sense. This ``image'' will still be denoted by $\mathcal{Z}(m)$ and, if not otherwise specified, we will always use $\mathcal{Z}(m)$ to denote the ``image''. If we want to distinguish them, we use ``$^{\mathrm{O}}\cZ(m)$'' to denote the original special divisor before taking the image. The same convention applies to $Z(m)$.}
\item We  use $L''_P:=\End_{\text{special}}(\mathscr{A}_P)$ to denote the $\mathbb{Z}$-lattice of special endomorphisms of $\mathscr{A}_P$. The notation $L''$ is chosen to be compatible with \cite{MAT}. We equip $L''_P$ with a positive definite quadratic pairing $Q'$ by $s\comp s =[Q'(s)]\in \End(\mathscr{A}_P)\otimes \bQ$.\\ 
\end{enumerate}
\end{definition}

\begin{remark}\label{rmk:specialHilbert}
{The name ``special divisor'' also suggests the fact that $\cZ(m)_{\bQ}$ is a special subvariety in the sense of Moonen (\cite{BMI98}). In fact, from \cite[Page 434, Proposition 4.5.8]{AGHMP17}, we know that $^{\mathrm{O}}\cZ(m)_{\bQ}$ is a disjoint union of a bunch of Hilbert modular surfaces $\cH_{\lambda,\bQ}$ defined by GSpin Shimura data $(L_\lambda,Q)$ of signature $(2,2)$, where each $L_\lambda\subseteq L_{\bQ}$ is a certain signature $(2,2)$ $\bZ$-lattice orthogonal to an element $\lambda\in L$ with $Q(\lambda)=m$. When $m$ is coprime to $p$ (which is the case that is relevant to us), then each $\cH_{\lambda,\bQ}$ admits an integral canonical model $\cH_{\lambda}$ over $\bZ_{(p)}$, and $^{\mathrm{O}}\cZ(m)$ is a union of them, see also \cite[\S 2.4]{ST22}.}
\end{remark}
\begin{remark}
The pairing $Q'$ should not be confused with $Q$. Nevertheless, $Q$ and $Q'$ are intimately related, see \cite[Lemma 2.3.2, Remark 3.1.4]{MAT}. \\
\end{remark}

Special endomorphisms are relevant to the study of the real multiplications and non-simpleness of an abelian surface:  \begin{lemma}\label{T:non-simple}
{}{If $P\in \mathcal{A}_{2,(p)}(\kk)$ lies in $Z(m)(\kk)$, then $\End\mathscr{A}_P$ contains $\mathbb{Z}[x]/(x^2-m)$. Moreover,  $P\in Z(k^2)(\kk)$ for some $k\in \bZ^+$ if and only if $\mathscr{A}_P$ is non-simple.} %Furthermore, $\mathscr{A}_P[p^{\infty}]$ is non-simple if and only if $P\in Z(m)(\kk)$ for some $m$ with $\left(\frac{m}{p}\right)=1$.
\end{lemma}
\proof Let $s$ be the special endomoprhism such that $s\comp s=m$, then $$\mathbb{Z}[x]/(x^2-m)\simeq \mathbb{Z}[s]\hookrightarrow \End\mathscr{A}_P.$$
%Conversely, suppose that $\mathbb{Z}[x]/(x^2-m)\hookrightarrow \End\mathscr{A}_P$, then $x$ is an endomorphism of trace 0 such that $x\comp x=[m]$. %{}{If $\mathscr{A}_P$ is simple, then it is ordinary or almost ordinary. By \cite[Remark 5.9(2)]{Oortnotes}, $\End(\mathscr{A}_P)_{\bQ}$ is a field. Now $(x^\dagger)^2=m^\dagger=m$, so $x^\dagger$ is either $x$ or $-x$. It follows that $x^{\dagger}=x$, by the fact that $\Tr(xx^\dagger)>0$. Therefore $x$ is a special endomorphism. If $\mathscr{A}_P$ is not simple, then it is isogenous to a product of two elliptic curves $E_1\times E_2$. Then $([m]_{E_1},-[m]_{E_2})$ is a special endomorphism. In both cases, $P\in Z(m)(\kk)$.}
If $s$ is a special endomorphism squaring to $k^2$, then $\ker(s-[k])$ is a non-trivial subvariety of $\mathscr{A}_P$. Conversely, if $\mathscr{A}_P$ is non-simple, then $\bQ\oplus \bQ\subset \End(\mathscr{A}_P)_{\bQ}$. There is a sufficiently large $k\in \bZ^+$, such that $s=(k,-k)\in  \End(\mathscr{A}_P)$ is a special endomorphism squaring to $k^2$. $\hfill\square$

\subsubsection{The mod $p$ fiber} The mod $p$ fiber $\mathcal{A}_{2,\kk}$ is a smooth threefold, and admits a Newton stratification. Let $\mathcal{A}_{2,\kk}^{\SSS}\subseteq \mathcal{A}_{2,\kk}^{\AO}\subseteq \mathcal{A}_{2,\kk}$ be the supersingular and non-ordinary strata. It is well known that $\mathcal{A}_{2,\kk}^{\SSS}$ is a projective 1-dimensional scheme, while $\mathcal{A}_{2,\kk}^{\AO}$ is a noncompact surface. There are rich  results regarding the intersection of $Z(m)$ with these higher Newton strata. In the following, a $\kk$-point of $\mathcal{A}_{2,\kk}^{\SSS}$ is called superspecial, if $\mathscr{A}_P$ is a superspecial abelian surface. 
\begin{lemma}\label{T:Rank0}
Let $P$ be a $\kk$-point of $\mathcal{A}_{2,\kk}^{\AO}\setminus\mathcal{A}_{2,\kk}^{\SSS}$\footnote{By $\mathcal{A}_{2,\kk}^{\AO}\setminus\mathcal{A}_{2,\kk}^{\SSS}$ we mean the open stratum which is the complement of $\mathcal{A}_{2,\kk}^{\SSS}$ in $\mathcal{A}_{2,\kk}^{\AO}$.}, then $\rk L''_P\leq 1$.
\end{lemma}
\proof %{}{ By Dieudonné theory, $L''_P\otimes \bZ_p= \mathbf{L}_{\text{cris},P}(W)^{\varphi=1}$}, where $\varphi$ is the crystalline Frobenius over $\mathbf{L}_{\text{cris},P}(W)$. 
This is because the slope 0 part of $\mathbf{L}_{\text{cris},P}$ is of rank 1.$\hfill\square$
\begin{corollary}
\label{T:simple}
Let $P$ be a $\kk$-point of $\mathcal{A}_{2,\kk}^{\AO}\setminus\mathcal{A}_{2,\kk}^{\SSS}$. If $m$ is not a perfect square and $P\in Z(m)(\kk)$, then the abelian surface $\mathscr{A}_P$ is simple. 
\end{corollary}
\proof If $\mathscr{A}_P$ splits, then $\bQ\oplus \bQ\subseteq\End(\mathscr{A}_P)_{\bQ}$. But $P\in Z(m)(\kk)$, so $\bQ(\sqrt{m})\subseteq \End(\mathscr{A}_P)_{\bQ}$. Therefore there is a rank two $\mathbb{Z}$-lattice of special endomorphisms of $\mathscr{A}_P$, contradicting Lemma~\ref{T:Rank0}. $\hfill\square$

\begin{lemma}\label{T:intspecial}
Suppose $m$ is an integer with $\left(\frac{m}{p}\right)=1$. Then the intersection of $Z(m)$ with $\mathcal{A}_{2,\kk}^{\SSS}$ consists only of superspecial points. 
\end{lemma}
\proof {}{This result is well-known. We will give a short but not so self-contained proof. Let $P$ be a supersingular point that lies on $Z(m)$. Note that $(m,p)=1$. Let $\cH_\lambda$ be the integral Hilbert modular surfaces as in Remark~\ref{rmk:specialHilbert}.}
{}{Viewing $P$ as a supersingular point on the mod $p$ reduction of an  $\mathcal{H}_{\lambda}$ as above, we can use the theory of vertex lattices (cf. \cite[\S 5]{HP17}, but we will follow the notation in \cite[\S3.1.2]{MAT}) to show that the vertex lattice $\Lambda_P$ has type $t_P=2$, and conclude by \cite[Definition 3.1.3]{MAT} that $P$ is superspecial. Let $V=L_{\lambda,\bQ}$. There is an even integer $t_{\max}$ that depends only on $V$, such that $t_P\in [2,t_{\max}]$. The quantity $t_{\max}$ can be computed via \cite[(1.2.3.1)]{HP17}: In our case $\rk V=4$ and, since $p$ splits in $\bQ[x]/(x^2-m)$, we have $\det(V_{\bQ_p})=1$, so $t_P=t_{\max}=2$. } $\hfill\square$ 

\subsection{Deformation theory}
\label{3.1}
\subsubsection{Explicit deformation theory of p-divisible groups}\label{subsub:defompdiv}
Suppose $P$ is a $\kk$ point of $\mathcal{A}_{2,\kk}$. Let $$\widehat{\mathcal{A}_2}_P:=\left({\widehat{\mathcal{A}_{2,W}}}\right)_P$$
be the formal completion of the integral model $\mathcal{A}_{2,W}$ at $P$. One is able to describe the universal Dieudonné crystal over $\widehat{\mathcal{A}_2}_P$ following \cite[\S 7]{Fal99}, \cite[\S 4]{Moo98} and \cite[\S1.4-1.5]{KM09}. Recall from \S\ref{2.1} that $\mathcal{A}_{2,(p)}$ is the Shimura variety associated to $G=\GSpin(L,Q)$. {}{It also  follows from \cite[1.3.6(2)]{KM09} that} $$(\mathbf{L}_{\text{cris},P}(W),\mathbf{Q})\simeq (L,Q)\otimes W$$
and the actions of $G$ over $L$ is compatible with the action of $G_W$ over $\mathbf{L}_{\text{cris},P}$ . 

Write $\mathbf{H}'_{\text{cris},P}$ as a 3-tuple $(\mathbf{H}'_{\text{cris},P}(W),\Fil^1\mathbf{H}'_{\text{cris},P}(W),\varphi')$ , the filtration $\Fil^1\mathbf{H}'_{\text{cris},P}(W)$ of the Dieudonné module reduces modulo $p$ to the Hodge filtration given by the Hodge cocharacter $\overline{\mu}:\mathbb{G}_{m,\mathbf{k}}\rightarrow G_{\mathbf{k}}$ corresponding to $P$. Pick $\mu:\mathbb{G}_{m,W}\rightarrow G_W$ lifting $\overline{\mu}$ (which amounts to picking a lifting of $P$ to a $W$-point of $\mathcal{A}_{2,W}$), then $\widehat{\mathcal{A}_2}_P$ is isomorphic to the completion at identity of the opposite unipotent $U_{{\mu}^{-1}}$ corresponding to the cocharacter $\mu$. Pick coordinates, so that $\mathcal{O}(U_{\mu^{-1}})=\Spf W[[x_1,x_2,x_3]]$, where the zero section corresponds to the identity of ${U_{\mu^{-1}}}$. We lift the 
Frobenius to $W[[x_1,x_2,x_3]]$ by setting $x_i\rightarrow x_i^p$, and call it $\sigma$. 

There is an explicit description of the universal Dieudonné crystal $\widehat{\mathbf{H}}'_{\text{cris},P}$ over $\widehat{\mathcal{A}_2}_P$ which we now describe. One first writes down the following 3-triple
$$(\mathbf{H}'_{\text{cris},P}(W)\otimes \mathcal{O}({U_{\mu^{-1}}}),\Fil^1\mathbf{H}'_{\text{cris},P}(W)\otimes \mathcal{O}({U_{\mu^{-1}}}),u'(\varphi'\otimes \sigma)),
$$ where $u$ is the tautological $\mathcal{O}(U_{{\mu}^{-1}})$ point of ${U_{\mu^{-1}}}$. Then there is a unique connection ${\nabla}'$ such that  
$$ \widehat{\mathbf{H}}'_{\text{cris},P}=(\mathbf{H}'_{\text{cris},P}(W)\otimes \mathcal{O}({U_{\mu^{-1}}}),\Fil^1\mathbf{H}'_{\text{cris},P}(W)\otimes \mathcal{O}({U_{\mu^{-1}}}),\nabla',u'(\varphi'\otimes \sigma)), $$
see \cite[\S 4.3]{Moo98}. We also have an explicit description of the universal K3-crystal $\widehat{\mathbf{L}}_{\cris,P}$\footnote{
The crystal $\mathbf{L}_{\text{cris}}$ is not an $F$-crystal, whereas $\mathbf{L}_{\text{cris}}(-1)$ is. However, we will still call it an $F$-crystal.} over $\widehat{\mathcal{A}_2}_P$:  
\begin{equation}\label{eq:Lhat}
\widehat{\mathbf{L}}_{\cris,P}=(\mathbf{L}_{\cris,P}(W)\otimes \mathcal{O}({U_{\mu^{-1}}}),\Fil^\bullet \mathbf{L}_{\cris,P}(W)\otimes \mathcal{O}({U_{\mu^{-1}}}),{\nabla},u(\varphi\otimes \sigma)),
\end{equation}
 where $u$ is the composition of $u$
 with the natural morphism $G_W\rightarrow \text{SO}(L,Q)_W$, and $\varphi$ is  induced from $\varphi'$ by conjugation action. See \cite[\S 4.9]{MP15} for more details.

\subsubsection{Explicit computations for $\widehat{\mathbf{L}}_{\cris,P}$} \label{expdL}
The explicit form of $\widehat{\mathbf{L}}_{\cris,P}$ was computed in  \cite[\S 3.3]{MAT}, using the theory of vertex lattice (cf. \cite{HP17}). We cite the computational results here. %\label{T:formaldef}

Suppose that $P$ is a superspecial point. There is an identification $\widehat{\mathcal{A}_2}_P\simeq U_{\mu^{-1}}\simeq \Spf W[[x,y,z]]$ such that $P$ corresponds to the ideal $(x,y,z)$, and Frobenious $\sigma$ on $W[[x,y,z]]$ given by $\sigma(x)=x^p,\sigma(y)=y^p,\sigma(z)=z^p$, and a choice of basis $\{v_1,v_2,...,v_5\}$ on $\mathbf{L}_{\cris,P}(W)$, so that under this basis, 
the cocharacter $\mu$, the
tautological point $u$
and the Frobenius $\Frob=u(\varphi\otimes \sigma)$ for $\widehat{\mathbf{L}}_{\cris,P}$, as in (\ref{eq:Lhat}), are explicitly given by
\begin{equation}
    {\mu}: \mathbb{G}_m\rightarrow \sO(L)_W,\,\,t\rightarrow \text{diag}(t^{-1},1,t,1,1),
\end{equation}
\begin{equation}\label{eq:u}   
    u=\begin{bmatrix}
        1&x&-xy-\frac{z^2}{4\epsilon}&y&z\\
        0&1&-y&0&0\\
        0&0&1&0&0\\
        0&0&-x&1&0\\
        0&0&-\frac{z}{2\epsilon}&0&1
    \end{bmatrix},
\end{equation}
\begin{equation}\label{eq:Frobv}   
\Frob=\begin{bmatrix}-\frac{1}{p}(xy+\frac{z^2}{4\epsilon})&x&p&y&z\\-\frac{y}{p}&1&0&0&0\\
\frac{1}{p}&0&0&0&0\\
-\frac{x}{p}&0&0&1&0\\
-\frac{z}{2\epsilon p}&0&0&0&1\end{bmatrix}\sigma,
\end{equation}
where $\epsilon\in \mathbb{Z}_p^*$ is a Teichmüller lifting of some nonsquare element in $\mathbb{F}_p$.

Fix $\lambda\in \mathbb{Z}_{p^2}^*$ as a square root of $\epsilon$. Then $L''_P\otimes \mathbb{Z}_p$ is spanned by $\omega_1=\lambda(pv_1-v_3),\omega_2=pv_1+v_3,\omega_3=v_2,\omega_4=v_4,\omega_5=v_5$. Under this basis, the Frobenius on $\widehat{\mathbf{L}}_{\cris,P}$ and the Gram matrix of the pairing $Q'$ on $L''_P\otimes \bZ_p$ are given by 
\begin{equation}\label{eq:FrobF}
\Frob=(I+F)\sigma,\text{ where }F=\begin{bmatrix}\frac{1}{2p}(xy+\frac{z^2}{4\epsilon})&-\frac{1}{2\lambda p}(xy+\frac{z^2}{4\epsilon})&\frac{x}{2\lambda p}&\frac{y}{2\lambda p}&\frac{z}{2\lambda p}\\
\frac{\lambda}{2p}(xy+\frac{z^2}{4\epsilon})&-\frac{1}{2p}(xy+\frac{z^2}{4\epsilon})&\frac{x}{2 p}&\frac{y}{2 p}&\frac{z}{2 p}\\
\lambda y& -y&0&0&0\\
\lambda x & -x&0&0&0\\
\frac{\lambda z}{2\epsilon} &-\frac{z}{2\epsilon}&0&0&0
\end{bmatrix}.
\end{equation}
\begin{equation}\label{eq:Pairingthing}
    Q'=\begin{bmatrix}
        -p\epsilon&0 & & & \\
        0&p&&&& \\
        &&0&1&\\
        &&1&0&\\
         &&&&\epsilon
    \end{bmatrix}.
\end{equation}
The numbers $2,4$ appeared in $xy+\frac{z^2}{4\epsilon}$ and $\frac{z}{2\epsilon}$ are Teichmüller liftings of $2,4\in \mathbb{F}_p$, but the number 2 appeared in fractions like $\frac{1}{2p}$ or $\frac{x}{2p}, \frac{x}{2\lambda p}$ is the usual $2\in\bZ$. By abuse of notation, we will use $2,4$ to denote either Teichmüller liftings or usual ones. Actually, for establishing our results, we don not really need to distinguish them. 
\subsubsection{Horizontal extension of special endomorphisms}\label{cha4}
Let $P\in \mathcal{A}_{2}(\kk)$ and $\omega\in L''_P\otimes \mathbb{Z}_p \subseteq \mathbf{L}_{\cris,P}(W)$.
By \cite[\S1.5]{KM09}, there is a unique horizontal section $\tilde{\omega}$ of $\widehat{\mathbf{L}}_{\cris,P}\otimes K=\mathbf{L}_{\cris,P}(W)\otimes K[[x,y,z]]$ extending $\omega$. With notation as in \S\ref{3.1}, we have $\tilde{\omega}=F_\infty\omega$, where  $F_\infty=\prod_{i=0}^\infty (I+F^{[i]})$ and $F^{[i]}=\sigma^i(F)$. The section $\tilde{\omega}$ controls the deformation of $\omega$. Let $R=\mathbf{k}[[x_1,...,x_n]]$, with a morphism $\iota: \Spf R\rightarrow \Spf\mathbf{k}[[x,y,z]]$. Lift ${\iota}$ to a morphism $\Spf W[[x_1,,...,x_n]] \rightarrow \Spf W[[x,y,z]]$. We still call it $\iota$.
\begin{lemma}\label{T:PDlemma} Let $I$ be an ideal of $R$, and $\tilde{I}$ be the preimage of $I$ in $ W[[x_1,...,x_n]]$. Then $\omega$ deforms to $\Spf(R/I)$ if and only if $\iota^*\tilde{\omega}\in\mathbf{L}_{\cris,P}(W)\otimes \hat{D}_{I}$, where $\hat{D}_I=\hat{D}_{W[[x_1,...,x_n]]}(\tilde{I})\subseteq K[[x_1,...,x_n]]$ is the $p$-adic completion of the PD-envelope of $\tilde{I}$.
\end{lemma}
\proof  Let $\iota^*\mathbf{L}_{\cris}(\hat{D}_I)$ be the pullback of $\widehat{\mathbf{L}}_{\cris,P}$ to $\hat{D}_I$. It can be expressed as the following 4-tuple: 
$$\iota^*\mathbf{L}_{\cris}(\hat{D}_I)=(\mathbf{L}_{\cris,P}(W)\otimes \hat{D}_I, \Fil^\bullet\mathbf{L}_{\cris,P}(W)\otimes \hat{D}_I,  {\iota}^*\nabla, F_I),$$
where $\nabla$ is the connection on $\widehat{\mathbf{L}}_{\cris,P}$ (\ref{eq:Lhat}), and the Frobenius $F_I$ can be constructed by the formula given in \cite[\S 4.3.3]{Moo98}. By \cite[\S 2.3]{DJ95}, there is an equivalence of categories between $p$-divisible groups over $\Spf(R/I)$ and Dieudonné modules over $\hat{D}_I$.  As a result, there exists an endomorphism of the $p$-divisible group $\iota^*\mathscr{A}[p^\infty]$ deforming $\omega$, if and only if there exists an endomorphism of the corresponding Dieudonné modules over $\hat{D}_I$ deforming $\omega$, if and only if there is a horizontal section of $\iota^*\mathbf{L}_{\cris}(\hat{D}_I)$ extending $\omega$. 

Recall that $\tilde{\omega}$ is constructed as the unique horizontal section $\tilde{\omega}$ of $\mathbf{L}_{\cris,P}(W)\otimes_W K[[x,y,z]]$ extending $\omega$. If $\omega$ deforms to $\Spf(R/I)$, then there is a horizontal section $\omega'$ of $\iota^*\mathbf{L}_{\cris,P}(\hat{D}_I)$ extending $\omega$.  Since ${\iota}^*\tilde{\omega}$ is also a horizontal section extending $\omega$ (in a larger ring), we must have $\omega'={\iota}^*\tilde{\omega}$ by the uniqueness of horizontal extension of $\omega$. Therefore ${\iota}^*\tilde{\omega}\in \mathbf{L}_{\cris,P}(W)\otimes \hat{D}_I$. Conversely, if 
${\iota}^*\tilde{\omega}\in \mathbf{L}_{\cris,P}(W)\otimes \hat{D}_I$, then it is must be the  horizontal section of $\iota^*\mathbf{L}_{\cris,P}(\hat{D}_I)$ extending $\omega$, so $\omega$ deforms to $\Spf(R/I)$. $\hfill\square$ 
\subsection{Formal germs of higher Newton strata} The computations in \S\ref{expdL} enable us to explicitly describe the formal germs of non-ordinary and supersingular loci. In the following we will always identify $\widehat{\mathcal{A}_2}_P$ with $\Spf W[[x,y,z]]$ as in \S\ref{expdL}. So the geometric deformation space $\mathcal{A}_{2,\kk}^{/P}$ can be identified as $\Spf \kk[[x,y,z]]$. We write $\mathcal{A}_{2,\kk}^{\AO,/P}$ and $\mathcal{A}_{2,\kk}^{\SSS,/P}$ for the formal germs of non-ordinary and supersingular strata.  
\begin{proposition}[Ogus]\label{SSstratat}
Notation as above, $\mathcal{A}_{2,\kk}^{\AO,/P}$ is cut out by $xy+\frac{z^2}{4\epsilon}=0$, inside which the supersingular locus $\mathcal{A}_{2,\kk}^{\SSS,/P}$ is furthermore cut out by $x^py+xy^p+\frac{z^{p+1}}{2\epsilon}=0$.
\end{proposition}
\proof {}{This is in the same spirit of  \cite[Lemma 3.4.1]{MAT}.
Let $\Fil^{\bullet}$ be the Hodge filtration on $\widehat{\bbL}_{\cris,P}$.} Then ${\gr_{-1}\widehat{\bbL}_{\cris,P}(\kk[[x,y,z]]})={\Span_{\kk}\{v_1\}}$ (where $v_1$ is introduced in \S\ref{expdL}). By \cite[Proposition 11]{Ogus01}, the non-ordinary locus is cut out by $(p\Frob)|_{\Span_{\kk}\{v_1\}}
$ and the supersingular locus is further cut out by $(p^2\Frob)|_{\Span_{\kk}\{v_1\}}
$. Using the explicit description (\ref{eq:Frobv}), we find that $\mathcal{A}_{2,\kk}^{\AO,/P}$ is cut out by $xy+\frac{z^2}{4\epsilon}=0$, inside which the supersingular locus $\mathcal{A}_{2,\kk}^{\SSS,/P}$ is furthermore cut out by the second equation.  $\hfill\square$

\begin{remark}\label{sslocus}
A direct computation shows that the subscheme cut out by $xy+\frac{z^2}{4\epsilon}=0$ and  $x^py+xy^p+\frac{z^{p+1}}{2\epsilon}=0$ consists of $p+1$ lines with tangent directions
$(1,0,0), (0,1,0)$ and  $(\alpha^{-1},-\alpha ,2\lambda )$, where $\alpha/\lambda\in \bF_p^*$. Each line is of multiplicity 2.
\end{remark}

\section{Compatifications and logarithmic geometry}\label{compactification} In this section we review the (integral) toroidal compactifications for $\mathcal{A}_2$, its (logarithmic) universal family, and its special divisors. 
A standard reference for the theory of toroidal compactifications for $\mathcal{A}_g$ is \cite{FC}. For the theory of integral toroidal compactifications of PEL type Shimura varieties and GSpin Shimura varieties, we refer the readers to \cite{LKW}, \cite{MP11},\cite{MP19} and \cite{BHS19} (\cite{MP11} and \cite{MP19} are on the same topic, but cover slightly different things). For the theory of log geometry and log motives, we reder the readers to \cite{Kk89},\cite{logDiu},\cite{Shi00},\cite{LogAV} and \cite{MP11}. Note that $\mathcal{A}_2$ and its special divisors can be viewed as both PEL type Shimura viarieties and GSpin Shimura varieties. Once a cone decomposition is fixed, the integral toroidal compactifications arising from different references listed above are compatible, see \cite[Introduction and Theorem 2]{MP19}. 

\subsection{Compatifications for $\mathcal{A}_{2}$} We will review the Bailey--Borel and the toroidal compactifications for $\mathcal{A}_{2}$. We won't provide a detailed introduction to the theory of compactifications, as it is an extensive topic. In stead of reading the literature listed above which are hundreds of thousands of pages long, the readers can find a concise and self-contained introduction for the toroidal compactifications of GSpin Shimura varieties in \cite{ST22}. 

\subsubsection{The Bailey--Borel compactification}
The boundary components of $\mathcal{A}_{2,\mathbb{C}}$ consists of a unique 0-dimensional cusp and a unique 1-dimensional cusp (cf. \cite[V.2]{FC} and \cite[Example 3.28]{BHS19}). As a result, $\mathcal{A}_{2,\mathbb{C}}$ admits a Bailey--Borel compatification 
\begin{equation}\label{eq:BB}
\mathcal{A}_{2,\mathbb{C}}^{\text{BB}}=\mathcal{A}_{2,\mathbb{C}}\sqcup\mathcal{A}_{1,\mathbb{C}}\sqcup\mathcal{A}_{0,\mathbb{C}},
\end{equation}
where $\mathcal{A}_{1,\mathbb{C}}$ \textit{resp}. $\mathcal{A}_{0,\mathbb{C}}$ is $j$-line \textit{resp}. a single point. The compactification extends to the integral model $\mathcal{A}_{2,(p)}$, which is also called the minimal compactification, 
see \cite[V]{FC} and \cite[Theorem 3]{MP19}. The integral model $\mathcal{A}_{2,(p)}$ admits a projective normal compatification $ \mathcal{A}_{2,(p)}^{\text{BB}}$ over $\mathbb{Z}_{(p)}$, whose generic fibre recovers $\mathcal{A}_{2,\mathbb{C}}^{\text{BB}}$. The stratification (\ref{eq:BB}) on $\mathcal{A}_{2,\mathbb{C}}^{\text{BB}}$ extends to a flat stratification of $\mathcal{A}_{2,(p)}^{\text{BB}}$:  $$ \mathcal{A}_{2,(p)}^{\text{BB}}=\mathcal{A}_{2,(p)}\sqcup\mathcal{A}_{1,(p)}\sqcup\mathcal{A}_{0,(p)}$$ where $\mathcal{A}_{1,(p)}$ \textit{resp}. $\mathcal{A}_{0,(p)}$ is the canonical integral model of $ \mathcal{A}_{1,\bC}$ \textit{resp}. $\mathcal{A}_{0,\bC}$.
\subsubsection{The toroidal compatifications}\label{subsub:toroidalcomp} According to \cite[Example 3.28]{BHS19}, we fix representatives for the cusps, i.e., fix a rank 2 integral lattice $J$ for the 1-dimensional cusp and a rank 1 integral lattice $I$ for the 0-dimensional cusp. For the 1-dimensional cusp, there is a unique cone decomposition. For the 0-dimensional cusp, we take the first cone decomposition written in \textit{loc.cit}\footnote{The authors of \textit{loc.cit} have given two cone decompositions for the 0-dimensional cusp, depending on the existence of internal rays. }. The overall cone decomposition for both cusps will be denoted $\Sigma$. The toroidal compactification $\mathcal{A}_{2,\mathbb{C}}^{\Sigma}$ and its stratification is described thoroughly in \textit{loc.cit}. Indeed, there is a natural map 
\begin{equation*}
\pi_\bC:\mathcal{A}_{2,\bC}^{\Sigma}\rightarrow \mathcal{A}_{2,\bC}^{\text{BB}} 
\end{equation*}
respecting the stratifications. The boundary strata of $\mathcal{A}_{2,\bC}^{\Sigma}$ are: 
\begin{enumerate}
    \item\label{item:1}$\mathcal{B}_{1,\bC}=\pi^{-1}(\mathcal{A}_{1,\bC})$ is
    a certain Kuga--Sato variety over $\mathcal{A}_{1,\bC}$, with fiber dimension 1. See also \cite[\S 2.2.1]{ST22}.
    \item\label{item:2}  $\mathcal{B}_{0,\bC}=\pi^{-1}(\mathcal{A}_{0,\bC})$ is a projective line over $\mathcal{A}_{0,\bC}$\footnote{In the language of \cite{ST22}, the stratum $\mathcal{B}_1$ \textit{resp}. $\mathcal{B}_0$ corresponds to the boundary component of type II \textit{resp}. III.}.
\end{enumerate} 

One can show that $\Sigma$ is furthermore an admissible complete smooth cone decomposition, and the integral model $\mathcal{A}_{2,(p)}$ admits a smooth toroidal compactification  $\mathcal{A}_{2,(p)}^{\Sigma}$ with generic fiber $\mathcal{A}_{2,\mathbb{C}}$ (See \cite{FC} and \cite[Theorem 1, 2 and 4.1.5]{MP19}). The stratification of $\mathcal{A}_{2,\mathbb{C}}^{\Sigma}$ extends to a stratification of $\mathcal{A}_{2,(p)}^{\Sigma}$ with all boundary components being flat divisors. The map $\pi_\bC$ extends to a natural map $$\pi:\mathcal{A}_{2,(p)}^{\Sigma}\rightarrow \mathcal{A}_{2,(p)}^{\text{BB}}$$
respecting the stratifications. The boundary stratum $\mathcal{B}_1=\pi^{-1}(\mathcal{A}_{1,(p)})$ \textit{resp}. $\mathcal{B}_0=\pi^{-1}(\mathcal{A}_{0,(p)})$ of $\mathcal{A}_{2,(p)}^{\Sigma}$ is a Kuga--Sato variety \textit{resp}. projective line over $\bZ_{(p)}$, with generic fibers $\mathcal{B}_{1,\bC}$ \textit{resp}. $\mathcal{B}_{0,\bC}$.

We will fix this compatification through out the paper, and call it $\mathcal{A}_{2,(p)}^{\text{tor}}$. Note that the boundary $\mathcal{B}=\mathcal{A}_{2,(p)}^{\text{tor}}\setminus\mathcal{A}_{2,(p)}$ is a single irreducible divisor of $\mathcal{A}_{2,(p)}^{\text{tor}}$ with a stratification $\mathcal{B}=\mathcal{B}_0\sqcup\mathcal{B}_1$. 

We now look at the mod $p$ fiber of  $\mathcal{A}_{2,(p)}^{\tor}$.  We write $\mathcal{A}_{2,\kk}^{\AO,\tor}$ for the closure of $\mathcal{A}_{2,\kk}^{\AO}$ in $\mathcal{A}_{2,\kk}^{\tor}$. Since $\mathcal{A}_{2,\kk}^{\SSS}$ is already proper, its closure in $\mathcal{A}_{2,\kk}^{\tor}$ is itself. It is well-known that $\mathcal{A}_{2,\kk}^{\AO,\tor}\setminus\mathcal{A}_{2,\kk}^{\AO}$ classifies extensions of supersingular elliptic curves by $\bG_m$. Such extensions, on the other hand, are 1-dimensional families in $\mathcal{B}_{1,\kk}$, as seen from the Barsotti formula $$\text{Ext}^1_{\text{alg-group}}(E,\mathbb{G}_m)\simeq \text{Hom}_{\text{group}}(\mathbb{Z},E^\vee)\simeq E^{\vee}\simeq E.$$ For our purpose, we will only use the fact that $\mathcal{A}_{2,\kk}^{\AO,\tor}\setminus\mathcal{A}_{2,\kk}^{\AO}\subseteq \mathcal{B}_{1,\kk}$.
\subsection{{}{Backgrounds on logarithmic geometry}} We will provide a brief introduction to the various log constructions used in this paper and direct readers to the relevant literature for more detailed explanations.
\subsubsection{Log schemes}
A \textit{log scheme} is a triple $(S,M,\alpha)$, where $S$ is a scheme, $M$ is a sheaf of monoids over the étale site $S_{\et}$, and a morphism $\alpha: M\rightarrow \mathcal{O}_S$ with respect to the multiplicative structure on $\mathcal{O}_S$, such that $\alpha^{-1}(\mathcal{O}_S^*)\rightarrow \mathcal{O}_S^*$ is an isomorphism. Following conventions, we will usually denote $(S,M,\alpha)$ simply by $(S,M)$ or even $S$ if the log structure is clear from the context. One can take direct and inverse images of log structures along morphisms of schemes \S\cite[(1.4)]{Kk89}. One can also make sense of the logarithmic version of (immersion, finite type, flat, smooth, étale, etc) morphisms of log schemes. We refer the readers to \cite{Kk89} or \cite{Shi00} for a detailed explanation of the theory. %There is also a notion of \textit{log formal scheme}, which is a formal scheme $\mathfrak{X}$ equipped with a log structure $M$, defined in a similar manner as log schemes. See \cite[\S 2]{Shi00} for details. 

In our paper, we will only use log schemes which are \textit{fine} in the sense of \cite[\S2]{Kk89}. This roughly means that $M$ is associated to a monoid that is finitely generated and integral. We will also use the terminology \textit{saturated}, and \textit{sharp} as defined in \cite[Definition 1.2.1.2, 1.2.1.4]{MP11}: We say $M$ is {saturated} (\textit{resp}. sharp), if the stalk of $M/\cO_{S}^*$ at any geometric point of $S$ is a {saturated} (\textit{resp}. sharp) monoid. A log structure is called \textit{fs} if $M$ is fine and saturated. Most of the schemes considered in this paper (e.g. those arise from good toroidal compactifications) are equipped with fs sharp log structures:
\begin{lemma}\label{lm:logpullbackfs}
Consider a regular scheme $S$ with log structure associated to a reduced normal crossing divisor. Let $S'\rightarrow S$ be any morphism of schemes. Then the inverse image log structure on $S'$ is fs and sharp.
%Suppose that $r\in \bN$ and $t_1,...,t_r\in m$. Equip $\Spec R$ with a log structure $\alpha: \mathcal{O}_{\Spec R}^*\oplus \bN^r\rightarrow \mathcal{O}_{\Spec R}$ sending $u\oplus (n_1,...,n_r)$ to $ut_1^{n_1}...t_{r}^{n_r}$ (e.g., the log structures induced from normal crossing divisors, cf. \cite[Example (2.5)]{Kk89}). Then it is fs and sharp. Furthermore, for any local morphism of local rings $(R,m)\rightarrow (R',m')$, the  pullback log structure on $\Spec R'$ is fs and sharp. %As a consequence, it is legitimate to talk about the inverse image log 1-motive $\mathfrak{L}_{\Spec \kk[[t]]}$ or $\mathfrak{L}_{\Spec \kk[[t]]/t^n}$ and its log Dieudonné module.
\end{lemma}
\begin{proof}
The inverse image of a fine log structure is fine (\cite[\S 2.4.1]{Kk89}). The rest of the properties can be checked over the stalks over geometric points. That is, we need to show that $M/\cO_{S}^*$ is sharp and saturated over the stalks of geometric points. Therefore it suffices to show that case where $S'\rightarrow S$ is a local morphism of strict Henselian local schemes. In this case, we can write $S=\Spec (R,\mathfrak{m})$ and $S'=\Spec (R',\mathfrak{m}')$, and it suffices to show that the global section of $M/\cO_{S}^*$ is sharp and saturated. It follows easily from \cite[Example (2.5)]{Kk89}) that the log structure $M$ on $S$ is the one associated to a pre-log structure (\cite[(1.3)]{Kk89}):$$\bN^r\rightarrow \mathfrak{m}\hookrightarrow R$$ sending generators to a regular sequence.  
The inverse image log structure $M'$ on $S'$ is then the one associated to the pre-log structure $\bN^r\rightarrow \mathfrak{m}\rightarrow \mathfrak{m}'\hookrightarrow R'$. Note that the image of this map is disjoint from ${R'}^*$, so $M'=\bN^r\oplus {R'}^*$. In other words, $M'/{R'}^*=\bN^r$, which is evidently saturated and sharp. 
\end{proof}

\subsubsection{Log 1-motives} Let $(S,M)$ be an fs log scheme. We will consider the Kummer log flat topology on $(S,M)$ (\cite[\S1.2.2]{MP11}). The Kato's log torus $\bG_m^{\log}$ is a sheaf in the Kummer log flat topos:
$$ (T,N)\rightarrow \Gamma(T,N^{\mathrm{gp}}).$$
There is a natural embedding of the usual $\bG_m$ into $\bG_m^{\log}$. Let $\cJ$ be a semi-abelian scheme over $S$ which is an extension of an abelian scheme $B$ by an isotrivial torus $T$ (\textit{isotrivial} means locally constant in the finite étale topology). Write $T=\Hom(X,\bG_m)$, where $X$ is the character lattice. We can associated to $T$ a log torus $T^{\log}=\Hom(X, \bG_m^{\log})$. Pushing $\cJ$ along the embedding $T\hookrightarrow T^{\log}$, we get a sheaf $\cJ^{\log}$ which is an extension of $B$ by $T^{\log}$. A \textit{log 1-motive} (\cite[Definition 2.2]{LogAV} or \cite[Definition 1.2.2.1]{MP11}) over $(S,M)$ is a complex $$\mathbf{M}=[Y\xrightarrow{u} \cJ^{\log}]$$ in degrees -1 and 0, where $\cJ^{\log}$ is a sheaf that arises from the above process and $Y$ is a $\bZ$-lattice locally constant in the finite étale topology. %There is a similar definition for log 1-motives over fs log formal schemes. 

Log 1-motives are relevant to us since, roughly speaking, degenerations of abelian schemes over a complete local Noetherian normal ring correspond to log 1-motives, see \cite[Proposition 1.2.4.2]{MP11}.

\subsubsection{Log $p$-divisible groups}\label{subsub:logavandpdiv} Let $(S,M)$ be a locally Noetherian fs log scheme, equipped with Kummer flat log topology. Following \cite{logDiu}, let $$(\mathrm{fin}/S)_d\subseteq (\mathrm{fin}/S)_r.$$
be full subcategories of the category of sheaves of abelian groups on $S$, which are ``stable under the Cartier \underline{d}uality" or ``\underline{r}epresentable". We call a log scheme representing a sheave in $(\mathrm{fin}/S)_r$ a \textit{log finite group scheme over $S$}. By \cite[Proposition 1.4]{logDiu}, log finite group schemes are (classically) finite, log flat, fs and Kummer type log schemes over $S$.

Following \cite[\S 4.1]{logDiu}, for $*=d,r$, we can define the category $(p\mathrm{-div}/S)_{*}$ as the full subcategory of the category of sheaves of
abelian groups on $S$ consisting of objects $G$ satisfying\begin{itemize}
    \item  $G=\bigcup_{n} G([p^n])$.
    \item  $p:G\rightarrow G$ is surjective. 
    \item For each $n$, $G([p^n])$ belongs to $(\mathrm{fin}/S)_{*}$.
\end{itemize}
We will call an object in $(\mathrm{fin}/S)_{r}$ a \textit{log $p$-divisible group}.  

Let $\mathbf{M}=[Y\xrightarrow{u} \cJ^{\log}]$ be a log 1-motive. One can functorially associate a log $p$-divisible group $\mathbf{M}[p^{\infty}]\in (p\mathrm{-div}/S)_{d}$, see \cite[Proposition 3.5, Definition 3.6]{Sheer}. In particular, its $p^n$-torsion is represented by a log finte group scheme over $S$.

%Let $(S,M)$ be a locally Noetherian fs log scheme. In \cite[\S 3]{LogAV}, the authors define the so called \textit{log abelian varieties with constant degeneration} over $S$ out of \textit{pointwise polarizable log 1-motives over $S$}, and proved an equivalence of categories between the two kinds of objects (\cite[Theorem 3.4]{LogAV}). A \textit{log abelian variety} is, very roughly speaking, pointwisely log abelian varieties with constant degeneration, see \cite[Definition 4.1]{LogAV}. For example, $\mathcal{A}_{2,(p)}^{\tor}$, with the log structure is induced by the normal crossing boundary divisor (\cite[Example (2.5)]{Kk89}), is a Noetherian fs log scheme. Then there is a universal log abelian scheme over $\mathcal{A}_{2,(p)}^{\tor}$ (\cite[\S4.4]{}) that admits $\mathscr{A}$ as the semi-abelian part (\cite[\S 9.2]{LogAV}). To a log abelian variety, one can associate a log $p$-divisible group. 

\subsubsection{Log crystals and log Dieudonné crystals}\label{subsub:logcryd}
Let $k/\mathbb{F}_p$ be a perfect field and let $W:=W(k)$. We regard $\Spec W$ as a log scheme with a fine log structure. Let $(S,M)$ be a fine log scheme over $W$. There is a notion of \textit{log crystalline site} $(S/W)_{\cris}^{\log}$ and \textit{log crystals} introduced by Kato in \cite[\S5,\S6]{Kk89} (strictly speaking, Kato only defines  $(S/W_n)_{\cris}^{\log}$ where $W_n$ is a truncation, but the definition generalizes to $S/W$). When talking about crystalline sites, we will always assume that the log structure on $W$ is trivial. 

Following \cite[\S  1.3.3]{MP19}, a \textit{log Dieudonné crystal} on an fs and sharp log scheme $(S,M)$ is a four tuple $(\mathbb{M},F,V,\Fil^\bullet)$, where $\mathbb{M}$ is a crystal on $(S/W)_{\cris}^{\log}$, $\varphi$ is  Frobenius, $V$ is Verschiebung and $\Fil^\bullet$ is a two step Hodge filtration. One can functorially attach log Dieudonné crystals to log 1-motives over locally $p$-nilpotent fs sharp log schemes (see \cite[\S 1.3.3]{MP11}).  

A crystalline Dieudonné functor for log $p$-divisible groups over a certain $p$-nilpotent log scheme $S$ is also constructed in \cite[\S 5]{logDiu}. In summary, we have the following theories \begin{center}
    
\begin{tikzcd}
\text{log 1 motives} \arrow[d,"\text{\cite{Sheer}}"] \arrow[rd,"\text{\cite{MP11}}"] &                              \\
\text{log $p$-divisible groups} \arrow[r,"\text{\cite{logDiu}}"] & \text{log Dieudonné crystals}
\end{tikzcd}
\end{center}

Though not explicitly written down in the literature, it is expected that the above constructions are compatible. Indeed, the construction in \cite{MP11} traces back to Kato. Despite this, we won't use the compatibilities of these constructions. When talking about the crystalline Dieudonné functor, we only use the one constructed in \cite{MP11}.

\subsection{{}{Special divisors near a boundary point}}\label{subsub:log1motive} By construction, the scheme $\mathcal{A}_{2,(p)}^{\tor}$ with the log structure induced from the normal crossing boundary divisors is fs and sharp. Let $k$ be a finite field or $\kk$. Let $P$ be a $k$-point of the boundary, and let $R$ be the complete local ring of $\mathcal{A}_{2,W}^{\tor}$ at $P$, where $W=W(k)$. Equip $\Spec R$ with the inverse image log structure from $\mathcal{A}_{2,W}^{\tor}$, which is fs and sharp by Lemma~\ref{lm:logpullbackfs}. Then the universal family over $\Spec R$ can be thought of as a log 1-motive. In fact, by \cite[Proposition 1.2.4.2]{MP11}, the abelian scheme $\mathscr{A}_{(\Spec R) -D}$ gives rise to a positive polarized log 1-motive \begin{equation}\label{eq:motive1}
    \mathfrak{L}=[Y\xrightarrow{u} \cJ^{\log}]
\end{equation} over $R$. Here $D$ is the boundary divisor in $\Spec R$, $Y$ is a constant sheaf of free abelian groups (of rank at most 2), $\cJ$ is an extension of an abelian variety $B$ (of dimension at most 1) with a torus $T$. The polarization and endomorphisms on $\mathscr{A}_{(\Spec R) -D}$ also have their avatars on $\mathfrak{L}_{R}$.

One can attach a log Dieudonné crystal $\bD(\mathfrak{L})$ over $R$ to the log 1-motive $\mathfrak{L}$ (\cite[Proposition 2.4.1.1]{MP11} or \cite[Proposition 1.3.5]{MP19}). The crystal $\bD(\mathfrak{L})$ is further equipped with a three step weight filtration $W_\bullet \bD(\mathfrak{L})$, such that  $$W_{-1}=0,\;W_0=\underline{\Hom}(Y,\mathbf{1}),\; W_1=\bD(\cJ),\; W_2= \bD(\mathfrak{L});\; \gr_{1}^W=\bD(B),\;\gr_2^W=\bD(T).$$ Here $\mathbf{1}:=\bD(\bQ_p/\bZ_p)$. Important $p$-adic Hodge theoretical properties of log crystalline Diuedonné functors  can be found in \cite[Proposition 1.4.10]{MP19} and \cite[Proposition 2.4.1.1]{MP11} (the integral version is via a logarithmic version of Breuil--Kisin modules (\cite[\S 2.2.4]{MP11})).% If $S\rightarrow \Spec R$ is a morphism of fs log formal schemes over $W$, let $\bD(\mathfrak{L}_S)$ be the pullback of $\bD(\mathfrak{L})$ to $S$.

% Note that we  have defined $\bbH_{\cris}'$ in \S\ref{subsubsec:SA} as the  Dieudonné crystal for $\mathscr{A}$. We have slightly absused the notation here since $\bbH_{\cris,R}'$ is the natural extension of the original one to the boundary.

%or as a log abelian scheme (\cite[Definition 4.1]{LogAV}). The two notions are related, but are different: roughly speaking, a log abelian scheme is ``a quotient'' of the log 1-motive, in the sense that Tate curve is a quotient of $\bG_m(K)$ by $q^{\bZ}$. Our main reference here is \cite{LogAV} (the authors have a series of papers on log abelian schemes from Vol I to Vol IV) and Madapusi's Thesis \cite{MP11}. We will not bother the readers with the precise definitions in the log geometry, since we will only use them as a black box.

%The restriction of the category of log $F$-crystals over $ \mathcal{A}_{2,(p)}^{\tor}$ to $\mathcal{A}_{2,(p)}$ (considered as a log scheme with trivial log structure) is equivalent to the  category of $F$-crystals over $\mathcal{A}_{2,(p)}$ in a natural sense. of $\bbH_{\cris}'$ to $\mathcal{A}_{2,(p)}$ (considered as a log scheme with trivial log structure) is of course just the $F$-crystal $\bbH_{\cris}'$ we defined in \S\ref{subsubsec:SA}. That is the reason we abuse the notation and still use $\bbH_{\cris}'$ to denote the log $F$-crystal. 
\subsubsection{{}{The compactified special divisors}}\label{subsub:specialdivisorscompact}
In the following we will assume that $m$ is coprime to $p$. Let ${\mathcal{H}_\lambda}$ be the set of integral Hilbert modular surfaces that unions to $^{\mathrm{O}}\cZ(m)$  as in Remark~\ref{rmk:specialHilbert}. From \cite[\S 2.1.28]{MP19}, each  $\mathcal{H}_\lambda$ admits a toroidal compactification $\mathcal{H}_\lambda^{\tor}$ smooth over $\bZ_{(p)}$ which is induced from the chosen cusp representative and the cone decomposition $\Sigma$, and the map $\mathcal{H}_\lambda\rightarrow \mathcal{A}_{2,(p)}$ extends to a map $\mathcal{H}_\lambda^{\tor}\rightarrow \mathcal{A}_{2,(p)}^{\tor}$. The compactification of $^{\mathrm{O}}\mathcal{Z}(m)$ is the union of $\mathcal{H}_\lambda^{\tor}$, which is smooth and will  be denoted by $^{\mathrm{O}}\overline{\mathcal{Z}(m)}$. The image of  $^{\mathrm{O}}\overline{\mathcal{Z}(m)}$ in $\mathcal{A}_{2,(p)}^{\tor}$ is then the Zariski closure of  
$\mathcal{Z}(m)$ in $\mathcal{A}_{2,(p)}^{\tor}$, which will be denoted by $\overline{\mathcal{Z}(m)}$.
%It is tempting to use the notation $\mathcal{Z}(m)^{\tor}$ to denote the compactification. However, this notation is reserved for a class of more important divisors of modularity significance, see \S\ref{subsub:modularspecial}. 
If not otherwise specified, the notation $\overline{Z(m)}$ will be used to denote $\overline{\cZ(m)}_{\kk}$.  \\

\begin{definition}\label{def:logspecialendo}
Let $k$ be a finite field or $\kk$, and $W=W(k)$. Let $P$ be a $k$-point in the boundary of $\mathcal{A}^{\tor}_{2,(p)}$ and let $R$ and $\mathfrak{L}$ be as in \S\ref{subsub:log1motive}. Suppose that $R'$ is a Noetherian local algebra over $W$ such that $p$ is nilpotent, and $S=\Spec R'\rightarrow \Spec R$ is a local morphism of fs sharp local 
 log schemes over $W$. We call an $s\in \End(\mathfrak{L}_{S})$ \textit{resp}. $\End(\mathfrak{L}_{S}[p^\infty])$ \textit{resp}. $ \End(\bD(\mathfrak{L}_{S}))$ 
a \textit{special endomorphism}\footnote{Ideally, one should define elements in $\End(\mathfrak{L}_{S})$ and $\End(\mathfrak{L}_{S}[p^\infty])$ to be special, if their crystalline realizations in  $ \End(\bD(\mathfrak{L}_{S}))$ satisfies $\Tr{s}=0$ and $s^\dagger=s$. But let's keep this ad hoc definition for now.}, if $\Tr{s}=0$ and $s^\dagger=s$. \\  
\end{definition}

By functoriality, a special endomorphism $s\in \End(\mathfrak{L}_{S})$ induces special endomorphisms of    $\mathfrak{L}_{S}[p^\infty]$ and $\bD(\mathfrak{L}_{S})$, which are usually denoted by the same symbol. Moreover, let $S'\rightarrow S$, then the canonical map $\End(\mathfrak{L}_{S})\rightarrow \End(\mathfrak{L}_{S'})$ is an embedding. To see this, we can use the fact that the base changes of  map $u$ in (\ref{eq:motive1}) to $S$ and $S'$ are injective (due to the existence of a polarization, see \cite[Remark 2.8.4]{LogAV}) to reduce to  the injectivity of $\End(\mathcal{J}_{S}^{\log})\rightarrow \End(\mathcal{J}_{S'}^{\log})$, which further reduces to the injectivity of $\End(\mathcal{J}_{S})\rightarrow \End(\mathcal{J}_{S'})$ by \cite[Proposition 2.5]{LogAV}.

As in the classical case, for a special endomorphism $s\in \End(\mathfrak{L}_{S})$, the self product $s\comp s$ is a scalar. To see this, one uses the injectivity $\End(\mathfrak{L}_{S})\hookrightarrow \End(\mathfrak{L}_{S'})$ to reduce to the case where $S$ is the log point $P$. Let \begin{equation}\label{eq:M_0m}
   M_0=\bD(\mathfrak{L}_P)(W),
\end{equation}where $W$ is equipped with a log structure induced from the log structure on $P$ (see \cite[\S 2.4.1]{MP11}). Note that $M_0$ is a rank four $W$-module with a polarization $\psi_0$. In the following, let $\End(M_0)$ be the endomorphism algebra of $M_0$ as a $W$-module, which is the matrix algebra $\mathrm{M}_4(W)$. If we represent the image of $s$ in $\End(M_0)_{\bC_p}$ by a $4\times 4$ matrix $N$, then $\Tr N=0$, and $N$ is invariant under a certain involution $x\rightarrow x'$ on $\End(M_0)_{\bC_p}$ induced by $\psi_0$. It suffices to show that $N^2$ is scalar. This is purely a linear algebra question, which can be approached in the same manner as in the classical setting. In fact, under suitable change of basis, we can identify $\End(M_0)_{\bC_p}$ with $\mathrm{M}_2(B)$, where $B=\mathrm{M}_2(\bC_p)$, 
and the involution is identified as $$x'=\prescript{t}{}x{^\iota}$$ 
where $^t(-)$ is the matrix transposition of $\mathrm{M}_2(B)$, and $(-)^\iota$ is the involution on $B$ given by $$\begin{bmatrix}
    a & b \\
    c & d
\end{bmatrix}^\iota=\begin{bmatrix}
    d & -b \\
    -c & a
\end{bmatrix}.$$ 
By \cite[Appendix A.3]{KR00}, if $x\in \mathrm{M}_2(B)$ is such that $x'=x$ and $\Tr x=0$, then $x^2=xx'$ is a scalar. This shows that our $N^2$ is a scalar matrix, hence $s\comp s$ is a scalar. 

In the following, let $(L''_P,Q')$ be the $\bZ$-lattice of special endomorphisms over $\mathfrak{L}_{P}$, where $Q'$ is defined by $s\comp s=[Q'(s)]$. This generalizes the lattice $L''_P$ in Definition~\ref{ddd1} to the boundary. 
\begin{lemma}\label{lm:boundaryrank}
    Suppose that $P$ is a $\kk$-point in  $\mathcal{A}_{2,\kk}^{\AO,\tor}\setminus\mathcal{A}_{2,\kk}^{\AO}$. Then $\rk L''_P\leq 1$. %The equality holds if and only if $\mathfrak{L}_P=[Y_P\xrightarrow{u} \cJ^{\log}_P]$ splits into a supersingular elliptic curve $[0\rightarrow E]$ and a ``Tate's curve'' $[\underline{\bZ}\rightarrow \bG_m^{\log}]$ up to isogeny (in the sense that there is an isogeny $\cJ_P\rightarrow \cJ_P'$ and $[Y_P\rightarrow {\cJ'}^{\log}_P]$ splits as above). In this case, every element $s\in L''_P$ satisfies $Q'(s)=n^2$ for some $n\in \bZ$. 
\end{lemma}
\begin{proof}
Our assumption on $P$ implies that $Y_P=\underline{\bZ}$, and $\cJ_P$ is an extension of a supersingular elliptic curve $E$ by $\bG_m$.
Suppose that there exists a nonzero $s\in L''_P$. We write it into the following form:
\begin{equation*}
 \begin{tikzcd}
\underline{\bZ} \arrow[r,"u"] \arrow[d, "s_{-1}"'] & \cJ_P^{\log} \arrow[d, "s_0"] \\
\underline{\bZ} \arrow[r,"u"]                      & \cJ_P^{\log}                 
\end{tikzcd}   
\end{equation*}
Since $s^\dagger=s$, $s_0$ uniquely determines $s$ via polarization. By \cite[Proposition 2.5, \S6]{LogAV}, restricting an endomorphism to $\cJ_P$ induces an isomorphism $$\End(\cJ_P^{\log})\xrightarrow{\sim}\End(\cJ_P).$$
Now we have $s_0\comp s_0=[Q'(s)]$, but $s_0$ cannot be the scaling map. The only case that this can happen is that $\cJ_P$ splits as $\bG_m\times E$ up to isogeny, and $s_0$ acts as $(n,-n)$ in the rational endomorphism ring $\End(\cJ_{P})_{\bQ}=\End(\bG_m)_{\bQ}\times \End(E)_{\bQ}$ (in fact, up to isogeny $\mathfrak{L}_P$ splits as the sum of an elliptic curve $[0\rightarrow E]$ and a ``Tate's curve'' $[\underline{\bZ}\rightarrow \bG_{m}^{\log}]$. We don't need this). Given another nonzero $s'\in L_P''$, the same reasoning shows that the corresponding $s_0'=(n',-n')\in \End(\bG_m)_{\bQ}\times \End(E)_{\bQ}$. Identify $\End(\bG_m)_{\bQ}=\bQ$.  There is some rational number $r$ such that $rn=n'$. So $rs_0=s_0'$. Via polarization this implies that $rs=s'$. Therefore $L''_P= 1$ in this case.
%then determines a splitting $\mathfrak{L}_P=[0\rightarrow E]\oplus [\underline{\bZ}\rightarrow \bG_m^{\log}]$ up to isogeny where $s$ acts as $(n,-n)$. Note that this already implies that $L''_P\leq 1$.
%Conversely, if $\mathfrak{L}_P$ splits up to isogeny, it is easy to reverse the above argument to see that any special endomorhism is of form $(n,-n)$.  
\end{proof}
\subsection{Deformation theory for log 1-motives} We will need an explicit deformation theory for log 1-motives in the spirit of \S\ref{subsub:defompdiv}, which is worked out in \cite{MP11}. This can be seen as a logarithmic generalization of the explicit deformation theory of $p$-divisible group developed in \cite[\S 7]{Fal99}, \cite[\S 4]{Moo98} and \cite[\S1.4-1.5]{KM09}.

Let $P$ be a $k$-point in the boundary of $\mathcal{A}^{\tor}_{2,(p)}$, with the inverse image log structure. Let $R$ and $\mathfrak{L}$ be as in \S\ref{subsub:log1motive}. The ring $R$ is the universal deformation ring of the polarized motive $\mathfrak{L}_P$ in the sense of \cite[\S 3.1]{MP11}.

Now starting from $M_0=\bD(\mathfrak{L}_P)(W)$ and a reductive group $G$ fixing certain Frobenius-invariant tensors $\{s_{\alpha,0}\}\subseteq M_0^{\otimes}$, a map $R\rightarrow R_G$ in explicit coordinates is constructed in \cite[\S3.2, \S 3.3]{MP11}, which is the normalization of a continuous quotient. This explicit map is the logarithmic generalization of the explicit coordinates constructed from opposite unipotents as per \cite[\S1.5.4]{KM09}, and is used to analyze the formal neighborhood at a $k$-point of a toroidal compactification of the integral canonical model of a Shimura variety of Hodge type (see \cite[Chapter 4]{MP11}), which closely follows the framework of Kisin's work.

\subsubsection{Deformation of special endomorphisms}\label{subsub:dese} The goal is to deduce an (extremely weak) ``analytic local moduli interpretation'' of $\overline{\cZ(m)}$ in terms of special endomorphisms of log 1-motives, which is enough for use in \S\ref{subsub:locinterbadloci}. We will assume that $m$ is coprime to $p$. For our sake, we will assume that the base change of $P$ to $\kk$ lies in  $\mathcal{A}_{2,\kk}^{\AO,\tor}\setminus\mathcal{A}_{2,\kk}^{\AO}$. 

Let $M_0=\bD(\mathfrak{L}_P)(W)$ be as in (\ref{eq:M_0m}), which is equipped with a polarization $\psi_0$ and a three step weight filtration $W_\bullet M_0$. Let $s_0$ be a special endomorphism of $\mathfrak{L}_P$ with $s_0\comp s_0=m$. We again use $s_0$ for its crystalline realization in $ M_0^{\otimes}$. Let $G$ be the reductive subgroup of $\GSp(M_0,\psi_0)$ fixing $s_0$.

We also make the assumption that $P$ lifts to a point $P'$ of $^{\mathrm{O}}\overline{\cZ(m)}$. Let $s$ be the tautological special endomorphism of $^{\mathrm{O}}\cZ(m)$, and apply \cite[Proposition 1.2.4.2]{MP11} to $^{\mathrm{O}}\overline{\cZ(m)}^{/P'}_W$, we see that $s$ gives rise to a special endomorphism of $\mathfrak{L}_P$ which, by Lemma~\ref{lm:boundaryrank}, must be $s_0$ or $-s_0$.

Following \cite{MP11}, we will run the construction for $R\rightarrow R_G$ in our setting.
\begin{enumerate}
\item Let $(Q,\lambda)$ \textit{resp}. $(Q_0,\lambda_0)$ be the log 1-motive $\mathfrak{L}$ \textit{resp}. $\mathfrak{L}_P$.
    \item Let $M_0^{\mathrm{sab}}=M_0/W_0M_0$; $\cP_{\mathrm{wt}}$ be the parabolic that fixes $W_\bullet M_0$; $U_{\mathrm{wt}}\subseteq \cP_{\mathrm{wt}}$ be the unipotent radical; $U_{\mathrm{wt}}^{-2}\subseteq U_{\mathrm{wt}}$ be the subgroup that acts trivially on  $M_0^{\mathrm{sab}}$; $\mathcal{P}_F$ be the parabolic that fixes the Hodge filtration.
  \item Let $\cP_{\mathrm{wt},G}$, $U_{\mathrm{wt},G}$, $U_{\mathrm{wt},G}^{-2}$, $\mathcal{P}_{F,G}$ be the intersections of $G$ with  $\cP_{\mathrm{wt}}$, $U_{\mathrm{wt}}$, $U_{\mathrm{wt}}^{-2}$, $\mathcal{P}_F$.
     \item Choose a cocharacter $\mu_0:\bG_{m}\otimes k\rightarrow \cP_{\mathrm{wt},G}\otimes k$ that splits the Hodge filtration $\Fil^\bullet(M_0\otimes k)$, which exists by \cite[Proposition 2.3.2.3]{MP11}(5) and the fact that $P$ lifts to $^{\mathrm{O}}\overline{\cZ(m)}$ (one can choose a $W$-point $x$ of $^{\mathrm{O}}\overline{\cZ(m)}$ lifting $P$ such that, over the generic point of $x$, one has an étale tensor $s_{\et,x}$ that arise from the
tautological special endomorphism $s$ on $^{\mathrm{O}}{\cZ(m)}$. Then $s_{\et,x}$ gives rise to $s_0$ or $-s_0$ via log $p$-adic Hodge theory). Let $\mu:\bG_m\rightarrow P_{\mathrm{wt},G}$ be a lift of $\mu_0$ splitting the Hodge filtration and satisfying some other conditions (see \cite[\S 3.2.3]{MP11}).
     \item  Let ${U}_F^{\mathrm{op}}\subseteq \GSp(M_0,\psi_0)$ and ${U}_{F,G}^{\mathrm{op}}\subseteq G$ be the opposite unipotents associated to $\mu$. 
     \item Let $U^+_F= {U}_F^{\mathrm{op}}/U_{\mathrm{wt}}^{-2}$ and $U^+_{F,G}= {U}_{F,G}^{\mathrm{op}}/U_{\mathrm{wt},G}^{-2}$. Let $\hat{U}^+$ \textit{resp}. $\hat{U}_{G}^+$ be the completion of $U^+_F$ \textit{resp}. ${U}_{F,G}^+$ at the identity section, it is the formal spectrum of a formally smooth $W$-algebra $R^+$ \textit{resp}. $R_G^+$. 
     \item Let $\sigma$ and $\sigma_G$ be certain polyhedral cones that depend only on $Q$ and $G$. In addition, let $\mathbf{S}_{Q,\sigma}$ and $\mathbf{S}_{\sigma_G}$ be certain monoids associated to these polyhedral cones. Furthermore, let $R_{\sigma}$ and $R_{\sigma_G}$ be completions of $\Spec W[\mathbf{S}_{Q,\sigma}]$ and $W[\mathbf{S}_{\sigma_G}]$ along certain sections, see \cite[\S 3.2.4,\S3.3.1]{MP11} for more details. 
  \item $R\rightarrow R_G$ is identified with $R^+\hat{\otimes}R_{\sigma}\rightarrow R^+_G\hat{\otimes}R_{\sigma_G}$.
\end{enumerate}
Since $P$ lifts to a point $P'$ of $^{\mathrm{O}}\overline{\cZ(m)}$ (which is smooth over $\bZ_{(p)}$ as we have noted in \S\ref{subsub:specialdivisorscompact}), the same proof of \cite[Theorem 4.2.3.1]{MP11} shows that $\Spf R_G\rightarrow \Spf R$ can be identified with \begin{equation}\label{eq:embcompl}
^{\mathrm{O}}\overline{\cZ(m)}_W^{/P'}\rightarrow \mathcal{A}_{2,W}^{\tor,/P}.\end{equation}
In particular, $R_G$ is smooth of dimension 2. Since
$^{\mathrm{O}}\overline{\cZ(m)}_W$ is étale locally a divisor, (\ref{eq:embcompl}) is an embedding, so $R_G$ is a quotient of $R$. This construction depends on the lifting $P'$ (and the choice of the $W$-point $x$), but we have the following:   
\begin{lemma}\label{lm:isoformalformal}
Suppose that $P$ lies in  $\mathcal{A}_{2,\kk}^{\AO,\tor}\setminus\mathcal{A}_{2,\kk}^{\AO}$, which lifts to a point $P'$ of $^{\mathrm{O}}\overline{\cZ(m)}$. Then for any lift of $P''$ of $P$ to $^{\mathrm{O}}\overline{\cZ(m)}$, the closed embedding \begin{equation}\label{eq:embcomp2}
    ^{\mathrm{O}}\overline{\cZ(m)}_W^{/P''}\hookrightarrow \mathcal{A}_{2,W}^{\tor,/P}=\Spf R
\end{equation} has image exactly $\Spf R_G$. 
\end{lemma}
\begin{proof}
Instead of working over $\kk$, let's work over sufficiently large finite field $k$. So we assume $P$ is defined over $k$ and $W=W(k)$, and prove the lemma in the new setting. Applying \cite[Proposition 1.2.4.2]{MP11} to $^{\mathrm{O}}\overline{\cZ(m)}^{/P''}_W$, we see that the tautological special endomorphism $s$ gives rise to a special endomorphism of $\mathfrak{L}_P$ which, by Lemma~\ref{lm:boundaryrank}, must be $\pm s_0$. For any $y\in {^{\mathrm{O}}\overline{\cZ(m)}}_W^{/P''}(W)$, the existence of the Hodge cycle $\pm s$ verifies the main conditions of \cite[Proposition 3.3.4.2]{MP11} (for readers convenience: this is an analogue of Kisin's result \cite[Proposition 1.5.8]{KM09}). The rationality assumption \textit{loc.cit} (3.3.3.5) can be verified by \textit{loc.cit} Corollary 4.2.2.5 (In order for the proof of  Corollary 4.2.2.5 to work, we need an $\mathcal{O}_L$-point of $^{\mathrm{O}}\overline{\cZ(m)}$ whose special fiber maps to $P$, this is guaranteed by the assumption that $P$ lifts to a point of $^{\mathrm{O}}\overline{\cZ(m)}$). Finally, the subcondition ``$y_{\sigma,0}:R_\sigma\rightarrow W/p=k$ factors through $R_{\sigma_G}$'' in \cite[Proposition 3.3.4.2]{MP11} is satisfied since the same thing holds for $P$. So $y$ factors through $\Spf R_G$. As a result, (\ref{eq:embcomp2}) factors through $\Spf R_G$. Since they both have dimension 2, the image of  (\ref{eq:embcomp2}) is exactly $\Spf R_G$.
\end{proof}
\begin{corollary}\label{cor:atmostliftings}
    Suppose that $P$ lies in  $\mathcal{A}_{2,\kk}^{\AO,\tor}\setminus\mathcal{A}_{2,\kk}^{\AO}$. Then there are at most two liftings of $P$ to $^{\mathrm{O}}\overline{\cZ(m)}$ (this can be seen as a very weak moduli interpretation: two liftings of $P$ correspond to two special endomorphisms $\pm s_0$ over $\mathfrak{L}_P$).
\end{corollary}
\begin{proof}
 If there is no lifting, then there is nothing to prove. Suppose that $P$ admits one lifting $P'$, so we can construct a deformation ring $R_G$ as above. By Lemma~\ref{eq:embcomp2}, if there are $n$ distinct liftings of $P$ to $^{\mathrm{O}}\overline{\cZ(m)}$, then there are $n$ distinct liftings of $R_G$ to $^{\mathrm{O}}\overline{\cZ(m)}$. Let $D$ be the boundary divisor of $R_G$. Then there are $n$ distinct liftings of $(\Spec R_G) - D $ to $^{\mathrm{O}}\overline{\cZ(m)}$. This means that $\mathscr{A}_{(\Spec R_G) - D}$ admits at least $n$ distinct special endomorphisms squaring to $m$. On the other hand, the lattice of special endomorphisms on $\mathscr{A}_{(\Spec R_G) - D}$ can only inject into the set of special endomorphisms on $\mathfrak{L}_{P}$: if it was not an injection, say some special endomorphism $s'$ maps to 0, then $0<m'=s'\comp s'$ maps to 0, which is absurd.  
Since $\rk L''_P=1$, we must have $n\leq 2$.
\end{proof}

\section{Arithmetic intersection theory}\label{sec:AIT}
This section is a summary of the constructions and results on arithmetic intersection theory that will be used in this paper. 
 \subsection{Global intersections and modularity}\label{C:intersection} 
 In this section we review the modularity result, which is the key input in estimating the intersection of $\overline{Z(m)}$ with a curve. 
 
 \subsubsection{Modularity of special divisors in $\mathcal{A}_{2,\kk}$} We follow \cite{AGHMP17}. For $\mu\in L^\vee/L, m\in \mathbb{Q}_{>0}$, let $\mathcal{Z}(m,\mu)$ be the corresponding special divisor over $\mathbb{Z}$ in $\mathcal{A}_2$ {}{defined in \cite[\S4.5, Definition 4.5.6]{AGHMP17}, which is flat over $\mathbb{Z}_{p}$ since our $L$ is self-dual at $p$} (cf. \cite[Proposition 5.2.1]{MP15} and \cite[Proposition 4.5.8]{AGHMP17}). {}{Roughly speaking}, $\mathcal{Z}(m,\mu)$ is the substack of $\mathcal{A}_2$ parametrizing abelian surfaces $A$ with a special quasi-endomorphism $s$ such that $Q(s)=m$, and the $l$-adic relization (\textit{resp}. {}{crystalline realization, which is only relevant when $A$ is defined in characteristic $p$}) of $s$ lies in the image of $(\mu+L)\otimes \mathbb{Z}_l$ (\textit{resp}. $(\mu+L)\otimes \mathbb{Z}_p$) in $\End(H^1_{\et}(A,\bQ_l))$ (\textit{resp}. $\End(\bD(A)\otimes\bQ)$, where $\bD$ stands for the Dieudonné module). By definition, we have $\mathcal{Z}(m,0)=\mathcal{Z}(m)$. 
 
 Let $Z(m,\mu):=\mathcal{Z}(m,\mu)_{\kk}$. Let $\mathcal{L}\in \text{Pic}(\mathcal{A}_{2,\kk})_\mathbb{Q}$ be the Hodge line bundle, which is the line bundle of weight 1 modular forms. Let $\mathfrak{e}_\mu$ be the standard basis of $\mathbb{C}[L^\vee/L]$. Let $\tau\in \bH$ and consider the following generating series: 
 $$\Phi_L:=\mathcal{L}^{-1}\mathfrak{e}_0+\sum_{m>0,\mu\in L^\vee/L}Z(m,\mu) q^m\mathfrak{e}_\mu
,\text{ where }q=e^{2\pi i\tau}.$$
Let $\rho_L$ denotes the Weil representation on $\mathbb{C}[L^\vee/L]$ and $M_{5/2}(\rho_L)$ be the space of vector valued modular forms of $\text{Mp}_2(\mathbb{Z})$ with respect to $\rho_L$ of weight $5/2$. The following proposition is a direct consequence of \cite[Theorem B]{HBK17}:
\begin{theorem}
$\Phi_L\in M_{5/2}(\rho_L)\otimes \Pic(\mathcal{A}_{2,\kk})_\mathbb{Q}$.
\end{theorem}
\subsubsection{Modularity of special divisors in $\mathcal{A}^{\tor}_{2,\kk}$}\label{subsub:modularspecial} One can extend the modularity results to the toroidal compactification of $\mathcal{A}_{2,\kk}$. We 
will use the notation $\overline{\mathcal{Z}(\mu, m)}$ to denote the closure of the special divisor $\mathcal{Z}(\mu, m)$ in $\mathcal{A}^{\tor}_{2,\kk}$. To every $\overline{\mathcal{Z}(\mu, m)}$, there is an associated divisor $$\mathcal{Z}^{\text{tor}}(\mu, m)= \overline{\mathcal{Z}(\mu, m)}+ \beta(\mu,m)\mathcal{B},$$
where $\beta(\mu_,m)\in \mathbb{R}$ is certain number defined in \cite{BHS19}, and $\mathcal{B}$ is the boundary divisor. Let ${Z}^{\text{tor}}(\mu, m)=\mathcal{Z}^{\text{tor}}(\mu, m)_{\kk}$, consider the following generating series 
$$\Phi^{\text{tor}}_L:=\mathcal{L}^{-1}\mathfrak{e}_0+\sum_{\mu\in L^\vee/L}\sum_{m\in Q(\mu)+\mathbb{Z}}{Z}^{\text{tor}}(\mu, m) q^m\mathfrak{e}_\mu
,\text{ where }q=e^{2\pi i\tau}.$$
The following proposition is a direct consequence of \cite[Theorem 3.1]{ST22}:
\begin{theorem}
$\Phi^{\tor}_L\in M_{5/2}(\rho_L)\otimes \Pic(\mathcal{A}^{\tor}_{2,\kk})_\mathbb{Q}$.
\end{theorem}
The theorem will not be used directly in this paper. But note that it is used in prove the result \cite[Proposition 4.10]{ST22} which we cited here as Lemma~\ref{intersEisen}. \subsubsection{Global intersections and Eisenstein series}\label{subsub;gine} Let $E_0(\tau)\in M_{5/2}(\rho_L)$ be the vector valued Eisenstein series with constant term $\mathfrak{e}_0$ as in \cite[\S 2.1]{Bru}(by taking $\kappa=5/2$). Denote the Fourier expansion of $E_0(\tau)$ as
$$E_0(\tau)= \sum_{m\geq 0,\mu\in \mathbb{Z}+Q(\mu)}q_L(m,\mu) q^m\mathfrak{e}_\mu
,\text{ where }q=e^{2\pi i\tau}.$$
Let $q_L(m):=q_L(m,0)$, the quantity $q_L(m)$ has the asymptotic growth $q_L(m)\sim m^{\frac{3}{2}}$ (\cite[\S 4.3.1]{MAT}. This estimation plays an important role in counting the intersection number $\overline{Z(m)}\cdot C$. The following is a direct consequence of \cite[Proposition 4.10]{ST22}: 
\begin{lemma}\label{intersEisen}
Let $C$ be a smooth projective curve mapping to $\mathcal{A}^{\tor}_{2,\kk}$, whose generic point lies in $\mathcal{A}_{2,\kk}$. Then for $m\gg 1$, one has $C\cdot \overline{Z(m)}\sim m^{\frac{3}{2}}$, and $C\cdot \overline{Z(m)}=|q_L(m)|(\mathcal{L}\cdot C)+o(m^{\frac{3}{2}}).$
\end{lemma}
The theorem \cite[Theorem 17]{BK01} provides an explicit expression of $q_L(m)$. Before stating it, we first recall the notion of local density, which will be used throughout the paper:
\begin{definition}
For a quadratic $\mathbb{Z}$-lattice $(L,Q)$, a prime $l$ and an integer $m$, the local density of $L$ representing $m$ over $\mathbb{Z}_l$ is defined as $$\delta(l,L,m)= \lim_{a\rightarrow \infty} l^{a(1-\rk\,L)}\#\{v\in L/l^aL|Q(v)\equiv m\mod l^a\}.$$
\end{definition}
\begin{theorem}[Bruinier--Kuss]\label{BrunierW}
Write $m=m_0f^2$ where $\gcd(f,2\det L)=1$ and $v_l(m_0)\in\{0,1\}$ for all $l\nmid2\det L$, then $$q_L(m)=\frac{-16\sqrt{2}\pi^2m^{\frac{3}{2}}L(2,\chi_D)}{2\sqrt{|L^\vee/L|}\zeta(4)}\left(\sum_{d|f} \mu(d)\chi_D(d)d^{-2}\sigma_{-3}(f/d)\right)\prod_{l|2\det L}\frac{\delta(l,L,m)}{1-l^{-4}},$$
where $\mu$ is the Möbius function and $D=-2m_0\det L$.
\end{theorem}
\subsection{Local intersections}\label{sub:localintersections} 
Let $C$ be a smooth projective curve and $\tau:C\rightarrow \mathcal{A}^{\tor}_{2,\kk}$ be a non-constant map, whose generic point lies in the almost ordinary locus of $\mathcal{A}_{2,\kk}$. Define 
$$l_P(m):=(C\cdot \overline{Z(m)})_P$$
to be the local intersection number of $\overline{{Z}(m)}$ with $C$ at a point $P\in C(\kk)$. We will equip $C$ with the inverse image log structure from $\mathcal{A}^{\tor}_{2,\kk}$. The goal of this section is to reduce the computation of $l_P(m)$ to a point-counting problem in quadratic lattices.  

Since the computation of $l_P(m)$ is a local question, we will pick a uniformizer $t$ of $C$ at $P$, write $C^{/P}=\kk[[t]]$ and work completely with $\kk[[t]]$. We remind the readers that $\Spf\kk[[t]]$ is a formal scheme over $\mathcal{A}^{\tor,/P}_{2,\kk}$ via $\tau^{/P}$.

\subsubsection{Local intersection in good locus}
Suppose that $P\in \mathcal{A}_{2,\kk}(\kk)$. For each $n\geq 1$ we define a sub-lattice $L''_{n,P}\subseteq L''_P$ by \begin{equation}\label{eq:LprprnP}
    L''_{n,P}:=\left\{s\in L''_P\,|\,s \text{ lifts to a special endomorphism of } \mathscr{A}_{\Spec\kk[[t]]/(t^n)}\right\}.
\end{equation}
The moduli interpretation of $\cZ(m)$ in Definition~\ref{ddd1} implies that  \begin{equation}\label{eq:l_Platticecounting}
    l_P(m)=\sum_{n\geq 1}\#\left\{s\in L''_{n,P}|Q'(s)=m\right\}.
\end{equation} 
Therefore, to compute $l_P(m)$ is to a count point in various quadratic lattices $(L''_{n,P},Q')$. 
\subsubsection{Local intersection in bad locus}\label{subsub:locinterbadloci} Let $\mathfrak{L}$ be the log 1-motive over $R=\mathcal{A}_{2,(p)}^{\tor,/P}$  (\S\ref{subsub:log1motive}). Suppose that $P\in \mathcal{A}_{2,\kk}^{\AO,\tor}(\kk)\setminus\mathcal{A}_{2,\kk}(\kk)$, considered as a log point with the inverse image log structure. We have defined the lattice of special endomorphisms $(L''_P,Q')$ over $\mathfrak{L}_P$ (\S\ref{subsub:specialdivisorscompact}), which is of rank at most 1 by Lemma~\ref{lm:boundaryrank}. By Lemma~\ref{lm:logpullbackfs}, the inverse image log structures on $\kk[[t]]$ and all $\kk[[t]]/t^n$ are fs and sharp, so it is legitimate to talk about the pullback of $\mathfrak{L}$ and its log Dieudonné crystal over these local rings. For each $n\geq 1$ we define a sub-lattice $L''_{n,P}\subseteq L''_P$ by \begin{equation}\label{eq:LprprnP2}
    L''_{n,P}:=\left\{s\in L''_P\,|\,s \text{ lifts to a special endomorphism of } \mathfrak{L}_{\Spec\kk[[t]]/(t^n)}\right\}.
\end{equation}
\begin{lemma}\label{lm:specialboundarycount} Notation as above. Suppose that $p\nmid m$, then 
    $$l_P(m)\leq \sum_{n\geq 1}\#\left\{s\in L''_{n,P}|Q'(s)=m\right\}.$$
\end{lemma}
\begin{proof}
  If $P$ does not lift to $^{\mathrm{O}}\overline{\cZ(m)}$, then $l_P(m)=0$, and there is nothing to prove. Suppose that $P$ lifts to a point $P'$ in $^{\mathrm{O}}\overline{\cZ(m)}$, then one can construct an explicit deformation ring $R_G$ as per \S\ref{subsub:dese}. By Corollary~\ref{cor:atmostliftings}, there can be at most two liftings of $P$ to $^{\mathrm{O}}\overline{\cZ(m)}$  and by Lemma~\ref{lm:isoformalformal}, the formal neighborhood of any lifting can be identified with $\Spf R_G$.  
  
Let $N$ be the largest integer such that $\Spec\kk[[t]]/t^N$ factors through $\Spf R_G$. Let $P''$ be a lifting of $P$ to $^{\mathrm{O}}\overline{\cZ(m)}$, then $N$ is the largest integer such that $\Spec \kk[[t]]/t^N$  factors through $^{\mathrm{O}}\overline{\cZ(m)}^{/P''}$. Since there are at most two liftings, it follows that $l_P(m)\leq 2N$. 

Lift $\Spec\kk[[t]]/t^N$ to $^{\mathrm{O}}\overline{\cZ(m)}^{/P''}$. The tautological special endomorphism $s_{\mathrm{taut}}$ on $^{\mathrm{O}}\overline{\cZ(m)}^{/P''}$ gives rise via \cite[Proposition 1.2.4.2]{MP11} a special endomorphism $s$ on $\mathfrak{L}_{\Spec\kk[[t]]/t^N}$. That is, $s\in L''_{N,P}$. Since $\rk L''_P=1$, this means $\#\left\{s\in L''_{n,P}|Q'(s)=m\right\}=2$ for each $n\leq N$. Therefore
$$N=\frac{1}{2}\sum_{n=1}^N \#\left\{s\in L''_{n,P}|Q'(s)=m\right\}\leq  \frac{1}{2}\sum_{n\geq 1} \#\left\{s\in L''_{n,P}|Q'(s)=m\right\}.$$
Combining $l_P(m)\leq 2N$ gives the estimate. 
\end{proof}

\subsection{Eisenstein series and quadratic forms}\label{sub:EsQl}
In \S\ref{sub:localintersections} we have reduced the computation of $l_P(m)$ to a counting problem in quadratic lattices. Counting problems of this kind are related to the theory of Eisenstein series. 

Let's assume that $P$ is a superspecial point, which is the case that is the most relevant to us. Recall that we have the lattice $(L''_P,Q')$. When the choice of $P$ is clear from the context, we abbreviate $L''_P$ as $L''$. Following \cite[\S2.3]{MAT}, let $L'$ be a lattice such that $L''\subseteq L'\subseteq L''\otimes \mathbb{Q}$, $L'\otimes \mathbb{Z}_p=L''\otimes \mathbb{Z}_p$ and is maximal over places $l\neq p$. The pairing $Q'$ on $L''$ extends naturally to $L''\otimes \mathbb{Q}$, hence to $L'$. It is known from \textit{loc.cit} Lemma 2.3.2 that $$(L'\otimes \mathbb{Z}_l,Q')\simeq (L\otimes\mathbb{Z}_l,Q)$$ for $l\neq p$. In particular, $L'$ and $L''$ have rank 5.

In the setup of \S\ref{sub:localintersections},  in order to estimate the quantity $l_P(m)$, we need to further consider {}{a certain} sublattice $L'''\subset L'$ such that $L'''\otimes \mathbb{Z}_l=L'\otimes \mathbb{Z}_l$ for $l\neq p$. {}{ The lattice $L'''$ will be made precise later and is related to $L''_{n,P}$. But for now, let's just regard it as a sublattice that does not carry a specific meaning.} 

Let $\theta_{L'''}(q)$ be the theta series of the positive definite lattice $L'''$, which is a modular form of weight $5/2$. By definition, the $m$-th coefficient of $\theta_{L'''}(q)$ computes the number of $v\in L'''$ with $Q'(v)=m$. Decompose the theta series into $E_{L'''}(q)+G_{L'''}(q)$, where $E_{L'''}$ is the Eisenstein series and $G_{L'''}$ the cusp form, and let $q_{L'''}(m)$ be the Fourier coefficient of $E_{L'''}(q)$. The proof of \cite[Theorem 17]{BK01} directly applies to the lattice $L'''$ and gives the following
\begin{theorem}\label{siegelmass}
 Write $m=m_0f^2$ where $\gcd(f,2\det L')=1$ and $v_l(m_0)\in\{0,1\}$ for all $l\nmid2\det L$, then $$q_{L'''}(m)=\frac{-16\sqrt{2}\pi^2m^{\frac{3}{2}}L(2,\chi_{D'})}{2\sqrt{|{L'''}^\vee/L'''|}\zeta(4)}\left(\sum_{d|f} \mu(d)\chi_{D'}(d)d^{-2}\sigma_{-3}(f/d)\right)\prod_{l|2\det L'}\frac{\delta(l,L''',m)}{1-l^{-4}},$$
where $\mu$ is the M\"obius function, $D'=-2m_0\det L'$ and $\delta(l,L''',m)$ is the local density of $L'''$ representing $m$ over $\mathbb{Z}_l$.
\end{theorem}
\begin{corollary}
\label{T:BKformula} For $L'''$ as above, we have 
$$\frac{q_{L'''}(m)}{-q_L(m)}=\frac{\delta(p,L''',m)\left[1-\left(\frac{m}{p}\right)p^{-2}\right]}{\sqrt{|(L'''\otimes \mathbb{Z}_p)^\vee/(L'''\otimes \mathbb{Z}_p)|}(1-p^{-4})}.$$
\end{corollary}
\proof As observed before, we have $(L'\otimes \mathbb{Z}_l,Q')\simeq (L\otimes\mathbb{Z}_l,Q)$. Since $L'''\otimes \bZ_l=L'\otimes \bZ_l$, we have $(L'''\otimes \mathbb{Z}_l,Q')\simeq (L\otimes\mathbb{Z}_l,Q)$. It follows that $\delta(l,L''',m)=\delta(l,L,m)$ for $l\neq p$. Since $L$ is self dual at $p$, we have $p\nmid \det L$. One then combines Theorem~\ref{BrunierW} and Theorem~\ref{siegelmass}.   $\hfill\square$\\[10pt]
We will also need the following result, which is essentially \cite[Proposition 9.1.5]{MAT}, for estimating the Fourier coefficient $q_G(m)$ of the cusp form $G_{L'''}(q)$:
\begin{proposition}\label{T:DWformula}
(Duke, Waibel) Let $S$ be a fixed finite set of primes and $\theta$ be the theta series attached to a positive definite quadratic lattice $(L''',Q')$ which is of rank $5$ with discriminant $D_\theta$, such that all prime factors of $D_\theta$ are in $S$. Then there exist absolutely bounded constants $N_0$ and $C_0$ such that $q_G(m)\leq C_0D_{\theta}^{N_0}m^{\frac{5}{4}}$.
\end{proposition} 
Finally we need a result by Hanke (\cite[Remark 3.4.1, Lemma 3.2]{Han04}), which computes local densities:
\begin{lemma}\label{T:Hanke}
If $(p,m)=1$ then $\delta(p,L''',m)=p^{1-\rk  L'''}\#\left\{v\in \frac{L'''}{pL'''}\bigg|Q'(v)\equiv m \mod p \right\}$.
\end{lemma}

\section{Algebraization theorems} \label{algebraization} In this section we prove Theorem~\ref{TT:4}. Recall the setup: let $U$ be a connected smooth (not necessarily projective) curve with a nonconstant morphism $U\rightarrow \mathcal{A}_{2,\overline{\mathbb{F}}_p}$, and $P$ being a $\overline{\mathbb{F}}_p$ point on $U$. We first use an idea of \cite{Zin01} to establish a general extension result of $p$-divisible groups, which enables us to descend the formal special endomorphism of $\mathscr{A}_{K(U^{/P})}[p^\infty]$ to an algebraic special endomorphism of $\mathscr{A}_{K(U)}[p^\infty]$. We then use $p$-adic isogeny theorem over function field (\cite[Theorem 2.6]{DJ98}) to pass from endomorphism rings of $p$-divisible groups to endomorphism rings of abelian surfaces, hence proving that $\mathscr{A}_{K(U)}$ admits a special endomoprhism.

\subsection{A result on $p$-divisible groups}
We call a ring extension $A\hookrightarrow B$ \textit{ $p$-integrally closed}, if $
\forall b\in B, b^p\in A\Rightarrow b\in A$. Note that if $A$ is reduced, then $B$ is also reduced. 
\begin{lemma}\label{cond}
If $A$ is a finite type connected smooth  $\mathbf{k}$-algebra, and $B$ is the completion of $A$ at a prime ideal, then $B$ is a normal domain. Furthermore, the extensions $A\hookrightarrow B$ and $\Frac(A)\hookrightarrow\Frac(B)$ are $p$-integrally closed.
\end{lemma}
\proof The first assertion follows from \cite[Scholie 7.8.3(v)]{EGAIV2}.
Now we prove that $A\hookrightarrow B$ is $p$-integrally closed. Denote the Frobenius of $A$ by $\sigma_A$, and the prime ideal in question by $\mathfrak{p}$, so $B=\varprojlim_n {A}/{\mathfrak{p}^n}$. A result of Kunz tells that $\sigma_A$ is flat, cf. \cite[Tag 0EC0]{stacks-project}. Since the systems of ideals $\{\mathfrak{p}^n\}_{n\geq 1}$ and $\{\sigma_A(\mathfrak{p})^n\}_{n\geq 1}$ are cofinal to each other, we have \begin{equation}\label{eq:cani}
B\otimes_{\sigma_A} A = \left(\varprojlim_n\frac{A}{\mathfrak{p}^n}\right)\otimes_{\sigma_A}A\simeq \varprojlim_n \left(\frac{A}{\mathfrak{p}^n}\otimes_{\sigma_A} A\right)\simeq B,\end{equation}
where the middle isomorphism, which tells that taking completion commutes with $\sigma_A$, follows from the fact that $A$ is a finite presented over itself via $\sigma_A$. Now suppose that $B$ is not $p$-integrally closed over $A$, i.e., 
$\exists a\in A, b\in B-A$ such that $b^p=a$. Then there is an injection $A[x]/(x^p-a)\hookrightarrow B$ which sends $x$ to $b$. By the flatness of $\sigma_A$, we have $$A[x]/(x-a)^p\simeq A[x]/(x^p-a)\otimes_{\sigma_A}A\hookrightarrow B\otimes_{\sigma_A}A\simeq B.$$
This implies that $B$ is nonreduced, contradicting the fact that $B$ is integral. 

Now we prove that the pair $\text{Frac}(A)\hookrightarrow\text{Frac}(B)$ is $p$-integrally closed. 
Observe first that $\text{Frac}(A)\hookrightarrow (A-\{0\})^{-1}B$ is $p$-integrally closed, then notice that $(A-\{0\})^{-1}B$ is integrally closed in $\text{Frac}(B)$. 
If $b\in \text{Frac}(B)$ such that $b^p\in \text{Frac}(A)$, then $b^p\in (A-\{0\})^{-1}B$. This implies that $b\in (A-\{0\})^{-1}B$, so $b\in \text{Frac}(A)$. $\hfill\square$

\begin{proposition}
\label{injj}
Suppose that $A\hookrightarrow B$ is a $p$-integrally closed extension of reduced $\bF_p$-algebras. Let $\mathscr{G}_1,\mathscr{G}_2$ be isoclinic p-divisible groups over $A$ with slopes $\mu_1>\mu_2$, and let $\mathscr{G}_i(n)$ be their truncations. Then the following natural morphisms are injective:\begin{enumerate}
    \item\label{injj1}$\Ext_{ A}^1(\mathscr{G}_2(n),\mathscr{G}_1(n)) \hookrightarrow \Ext_{B}^1( \mathscr{G}_{2,B}(n),\mathscr{G}_{1,B}(n))$,
    \item \label{injj2}$\Ext_{A}^1(\mathscr{G}_2,\mathscr{G}_1) \hookrightarrow \Ext_{B}^1( \mathscr{G}_{2,B},\mathscr{G}_{1,B})$,
\end{enumerate}
where $\mathscr{G}_{i,B}$, $\mathscr{G}_{i,B}(n)$ are the base changes of the respective objects to $B$. 
\end{proposition}
\proof  The idea is essentially from \cite{Zin01}, but we will make the proof as self-contained as possible.  Let $\sigma_A,\sigma_B$ be the Frobenius endomorphism over $A$,$B$, respectively. Let  $\mu_i=\frac{r_i}{s}$, $r_i,s\in \bZ$.

Since (\ref{injj2}) follows from (\ref{injj1}) by taking limit over various $n$, it suffices to show (\ref{injj2}). In the following, we fix an $n$. It suffices to show that if the extension splits over $ B$ then it splits over $A$. Let $\mathscr{G}(n)$ be an extension of $\mathscr{G}_1(n)$ by $\mathscr{G}_2(n)$, i.e., there is an exact sequence 
\begin{equation}\label{eq:splittingG_1G_2}
    0\rightarrow \mathscr{G}_1(n)\rightarrow \mathscr{G}(n)\rightarrow \mathscr{G}_2(n)\rightarrow0.
\end{equation}

Let $\cM$ (\textit{resp}. $\cL$) be the structure bi-algebra of $\mathscr{G}(n)$  (\textit{resp}. $\mathscr{G}_2(n)$) over $A$. There is a canonical embedding of bi-algebras $\cL\subseteq \cM$ induced from $\mathscr{G}(n)\rightarrow\mathscr{G}_2(n)$, which makes $\mathcal{L}$ locally a direct summand of $\cM$ (as $A$-modules). Let $\Phi$ be the $\sigma_A^s$-semilinear endomorphism on $
\cM$ induced by the isogeny $p^{-r_2}\mathrm{Fr}_{\mathscr{G}(n)}^{s}:\mathscr{G}(n)\rightarrow \mathscr{G}(n)^{(p^s)}$, which is well-defined since $r_1>r_2$. Then $\mathcal{L}$ is furthermore $\Phi$-invariant. Note that $\Phi$ is a $\sigma_A^s$-linear isomorhism on $\cL$ and is nilpotent on $\cM/\cL$ (and this uniquely characterizes $\cL$ -- we won't need this). Let $\cM^{\text{nil}}$ be the $A$-submodule of $\cM$ which is nilpotent under $\Phi$. By Lemma~\ref{lm:M+L} below, the splitting of  
(\ref{eq:splittingG_1G_2}) is equivalent to the following identity on $A$-modules
\begin{equation}\label{eq:sujnil}
   \mathcal{M}= \cM^{\text{nil}}+ \cL.
\end{equation}

Now consider $\mathcal{M}\otimes B$ (\textit{resp}. $\mathcal{L}\otimes B$), which is the structure bi-algebra of $\mathscr{G}_B(n)$ (\textit{resp}. $\mathscr{G}_{2,B}(n)$). Running the same set of constructions as above, we have a $\sigma_B^s$-semilinear endomorphism $\Phi\otimes \sigma_B^s$ on $\mathcal{M}\otimes B$ which, by abuse of notation, will still be denoted by $\Phi$. The map $\Phi$ is an $\sigma_B^s$-linear isomorphism on $\cL\otimes B$ and is nilpotent on $(\mathcal{M}\otimes B)/(\mathcal{L}\otimes B)\simeq (\cM/\cL)\otimes B$. Let $(\cM\otimes B)^{\text{nil}}$ be the $B$-submodule of $\cM$ which is nilpotent under $\Phi$. Possibly enlarge $N$, the map $\Phi^N$ annihilates $(\cM\otimes B)^{\text{nil}}$. Apply Lemma~\ref{lm:M+L} to $B$, the splitting of (\ref{eq:splittingG_1G_2}) over $B$ implies that \begin{equation}\label{eq:sujnil2}
    \cM\otimes B=(\cM\otimes B)^{\text{nil}}+(\cL\otimes B).
\end{equation}
We now use (\ref{eq:sujnil2}) as well as the $p$-integrality assumption to show that (\ref{eq:sujnil}) holds. This can be done locally. Possibly localizing $A$, we can assume that $\cM,\cL,\cM/\cL$ are all free\footnote{We may also need to localize $B$ in order to apply (\ref{eq:sujnil2}) to the localized $A$. This is totally fine since $p$-integrally closed extension is preserved by localization.}. Let $e_1,...,e_s$ be an $A$-basis of $\cL$. Let $e_1,...,e_s,e_{s+1},... e_t$ be an $A$-basis of $\cM$.
Since all modules involved are finite, there is a some sufficiently large $N$ such that $\Phi^N$ annihilates $\cM/\cL$, $\cM^{\text{nil}} $ and $ (\cM\otimes B)^{\text{nil}}$. It is easy to see that $\Phi^N(e_1),...,\Phi^N(e_s)$ is again an $A$-basis of $\cL$, while $\Phi^N(e_{s+1}),... \Phi^N(e_t)\in \cL$. In particular, expressing an $m\in \cM$ into a linear combination of basis elements, we have $\Phi^N(m)\in \cL$. 

From (\ref{eq:sujnil2}), given an element $m\in \cM$, there are $b_1,...,b_s\in B$ and $g\in (\cM\otimes B)^{\text{nil}}$ such that \begin{equation}\label{eq:m=g+l}
    m\otimes 1=g+\sum_{i=1}^s e_i\otimes b_i.
\end{equation} Take $\Phi^N$ on both sides. Since $g$ is killed, we find that $$\sum_{i=1}^s \Phi^N(e_i)\otimes b_i^{p^N}=\Phi^N(m)\otimes 1\in\cL.$$  
Since $\Phi^N(e_1),...,\Phi^N(e_s)$ is an $A$-basis of $\cL$, and $\Phi^N(e_1)\otimes 1,...,\Phi^N(e_s)\otimes 1$ is a $B$-basis of $\cL\otimes B$, we find that $b_i^{p^N}\in A$, hence $b_i\in A$ by $p$-integrality. Plugging this back to (\ref{eq:m=g+l}), we find that $g\in \cM$. Since $\Phi^N$ kills $g$, we have $g\in \cM^{\text{nil}}$. Again plugging this back to (\ref{eq:m=g+l}), we find that $m\in \cL+\cM^{\text{nil}}$,  verifying (\ref{eq:sujnil}).  $\hfill\square$ 
\begin{lemma}\label{lm:M+L}
In the context of the proof of Proposition~\ref{injj}, the splitting of (\ref{eq:splittingG_1G_2}) is equivalent to $ \mathcal{M}= \cM^{\text{nil}}+ \cL$. 
  \end{lemma}
 \begin{proof}
The readers are recommended to first read Example~\ref{ex:splitex} before proceed. We will use the same notation as in the proof of Proposition~\ref{injj}.

If (\ref{eq:splittingG_1G_2}) splits, write the structure algebra for $\mathscr{G}_1(n)$ as $\cN$, we have $\cM\simeq \cL\otimes \cN$. Let $e:\cN\rightarrow A$ be the co-unit, then $1\otimes \ker e$ is nilpotent under $\Phi$. Therefore $\cM^{\text{nil}}\supseteq \cL\otimes \ker e$. We immediately have $\mathcal{M}= \cM^{\text{nil}}+ \cL$ as $A$-modules. 

Conversely, suppose that $\mathcal{M}= \cM^{\text{nil}}+ \cL$ as $A$-modules. We always have $\cM^{\text{nil}}\cap \cL =\mathbf{0}$ (without reducedness this is not true, for example, the nilradical of $A$ will be contained in the intersection). It follows that $\mathcal{M}= \cM^{\text{nil}}\oplus \cL$ as $A$-modules. Since $\Phi$ is compatible with the $A$-algebra structure on $\mathcal{M}$, it is easy to see that $\cM^{\text{nil}}$ is an ideal of $\cM$. So it is legitimate to talk about the following morphism of $A$-algebras $$\pi:\mathcal{M}=\cM^{\text{nil}}\oplus \cL\rightarrow \cL, \,\,m=(m^{\text{nil}},m^{\text{iso}})\rightarrow m^{\text{iso}}.$$
Now we show that $\pi$ is a morphism of bi-algebras. First, note that $\cL\subseteq \cM$ is already a sub-bi-algebra that is carried identically onto its target (as bi-algebras). Second, note that the formation of $\Phi$ preserves the bi-algebra structure (since it comes from an isogeny of group schemes). Let $N$ be sufficiently large so $\Phi^N$ annihilates $\cM^{\text{nil}}$. To show that $\pi$ preserves the co-units, it suffices to show that for any $m^{\text{nil}}\in \cM^{\text{nil}}$, $e_{\cM}(m^{\text{nil}})=0$. This follows since $\sigma_A^Ne_{\cM}(m^{\text{nil}})=e_{\cM}(\Phi^{N}m^{\text{nil}})=0$ and $A$ is reduced. Let $\fm_{\cM}$ and $\fm_{\cL}$ be the co-multiplication laws on $\cM$ and $\cL$, respectively. To show that $\pi$ preserves co-multiplication laws, 
it suffices to show that for any $m^{\text{nil}}\in \cM^{\text{nil}}$, $\fm_{\cM}(m^{\text{nil}})\in \cM\otimes_A \cM$ maps to 0 in $\cL\otimes_A \cL$. Let the image of $\fm_{\cM}(m^{\text{nil}})$ in $\cL\otimes_A \cL$ be $\alpha$. We have \begin{align*}
    \Phi^N\otimes_A \Phi^N(\alpha)=(\pi\otimes_A \pi) \comp (\Phi^N\otimes_A \Phi^N) (\fm_{\cM}(m^{\text{nil}}))=\pi\otimes_A \pi (\fm_{\cM}(\Phi^N(m^{\text{nil}})))=0
\end{align*}
Note that $\Phi\otimes_A \Phi$ is a $\sigma_A^{s}\otimes_A \sigma_A^{s}$-linear isomorphism on $\cL\otimes_A \cL$. Since $\cL$ is a locally free $A$-module. By explicitly writing out the local basis, and using the reducedness of $A$, one finds that $\alpha=0$. This finishes the proof that $\pi$ is a morphism of bi-algebras. Therefore there is a morphism of group schemes $\mathscr{G}_2(n)\rightarrow \mathscr{G}(n)$  splitting (\ref{eq:splittingG_1G_2}).
\end{proof}
 \begin{example}\label{ex:splitex}
Let the setting be as in Lemma~\ref{lm:M+L}. Let's give two examples to show how $\cM^{\text{nil}}$ governs the splitting. Suppose that $n=1$.  

First, consider the split $A$-group scheme $\mathscr{G}(1)=\mathscr{G}_1(1)\times \mathscr{G}_2(1)$, where $\mathscr{G}_1(1)=\mu_{p,A}$ is the constant multiplicative group scheme over $A$, and $\mathscr{G}_2(1)=\underline{\bF_p}_A$ is its dual. The structure bi-algebra $\cM$ of $\mathscr{G}(1)$ is $\prod_{i \in \bF_p}A[x_i]/(x^p_i-1)$, the structure bi-algebra $\cL$ of $\mathscr{G}_2(1)$ is $\prod_{i \in \bF_p}A$. The module $\cM^{\text{nil}}$ is $\prod_{i\in \bF_p}(x_i-1)A[x_i]/(x^p_i-1)$. We have $\cM=\cM^{\text{nil}}\oplus \cL$ as $A$-modules. Furthermore, $\cM^{\text{nil}}$ is an ideal of $\cM$, such that $\cM/\cM^{\text{nil}}\simeq \cL$. 

Second, let $A=k=\bF_p(t)$. Let $\mathscr{G}_1(1)=\alpha_{p,k}$ be the local-local group scheme $\Spec k[x]/(x^p)$, and $\mathscr{G}_2(1)=\underline{\bF_{p}}_k$. Let $k'=\bF_p(t^{\frac{1}{p}})$ be an inseparable extension. For $i\in \bF_p$, let $\cM_i= k[x_i]/(x^p_i-it)$, $\cM=\prod_{i\in \bF_p} \cM_i$, and $\mathscr{G}(1)=\Spec  \cM$ with group structure $$ \cM_{i+j}\rightarrow\cM_i\otimes \cM_j,\,\,x_{i+j}\rightarrow x_i\otimes 1+1\otimes x_j.$$
Then $\mathscr{G}(1)$ is a nontrivial extension of $\mathscr{G}_1(1)$ by $\mathscr{G}_2(1)$ since $\cM_0$ is non-reduced, but $\cM_i= k'$ for $i\neq 0$. In this case, $\mathcal{M}^{\text{nil}}$ is only $(x_0)\cM_0$, and we have $\cM^{\text{nil}}+\cL\subsetneq \cM$.
 \end{example}
\subsection{Algbraization of formal special endomorphisms}
 \begin{theorem}
 \label{T: algebraicC}
Suppose that $U$ is a connected smooth (affine) curve with a nonconstant morphism $U\rightarrow \mathcal{A}_{2,\kk}$. Let $K(U)$ be the function field of $U$ and $\widehat{K}(U)$ be its completion at a place. If the image of $U$ is generically almost ordinary and $\mathscr{A}_{\widehat{K}(U)}[p^\infty]$ admits a special endomorphism, then the image of $U$ lies in a special divisor.
 \end{theorem}
\proof Let $D=K(U)$ and $E=\widehat{K}(U)=\hat{D}$. By shrinking $U$ one can assume that $\mathbf{L}_{\cris,U}(-1)$ has constant Newton polygon. $\mathbf{L}_{\cris,U}(1)$ admits a slope filtration in the category of $F$-crystals over $U$:
$$0\subseteq \Fil_{1/2}\subseteq \Fil_1\subseteq \mathbf{L}_{\cris,U}(1),$$
where $\Fil_{1/2}$ is the slope $1/2$ part and $ \Fil_1$ is the slope $\leq 1$ part.
Since $\Fil_1$ is actually a Dieudonné crystal over $U$, the main result of \cite{DJ95} implies that $\Fil_1=\bbD(\mathcal{G})$ where $\bbD$ is the covariant Dieudonne crystal functor and $\mathcal{G}$ is a $p$-divisible group over $U$. We have $\Fil_{1/2}=\bbD(\mathcal{G}_{\text{con}})$, where $\mathcal{G}_{\text{con}}$ and $\mathcal{G}$ sit inside the classical connected-etale sequence \begin{equation}\label{ecs}
    1\rightarrow \mathcal{G}_{\text{con}}\rightarrow \mathcal{G}\rightarrow \mathcal{G}_{\text{et}} \rightarrow 1,
\end{equation}
which corresponds via $\bbD$ to \begin{equation}\label{ecf}
    0\rightarrow \Fil_{1/2}\rightarrow \Fil_{1}\rightarrow \Fil_{1}/\Fil_{1/2}\rightarrow 0. \end{equation}
Let $\tilde{\omega}$ be the formal special endomorphism of $\mathscr{A}_{E}[p^\infty]$, which gives rise to a horizontal section of slope 1 of the crystal $\bbL_{\cris,\Spec E}(-1)$. Therefore (\ref{ecf}) splits over $\Spec E$. \cite[Theorem 4.6]{DJ99} then implies that (\ref{ecs}) splits over $\Spec E$. It follows from Proposition~\ref{injj} and Lemma~\ref{cond} that (\ref{ecs}), hence (\ref{ecf})
 splits over $\Spec D$. Therefore $\tilde{\omega}$ comes from the base change of a special endomorphism $\omega$ of $\mathscr{A}_D[p^\infty]$. Now \cite[Theorem 2.6]{DJ98} implies that $$ \omega\in \End_{D}\left(\mathscr{A}_D[p^\infty]\right)\simeq \End_D\left(\mathscr{A}_D\right)\otimes \mathbb{Z}_p.$$ 
Therefore $\Spec D\rightarrow \mathcal{A}_{2,\kk}$ factors through some $Z(m)$, hence $U$ factors though $Z(m)$. $\hfill\square$\\

\begin{remark}\label{remarkOth}
The proof of Theorem~\ref{T: algebraicC} generalizes to height $1/n$ strata of GSpin Shimura varieties of signature $(2,2n-1)$. The corresponding statement of the theorem would be the following:\begin{thm}
    Let $X$ be a smooth subscheme of $\mathcal{S}_{\kk}$ generically lying in the height $1/n$ strata. If a special endomorphism $\omega$ of the pullback p-divisible group at a point $P\in X$ extends to the formal neighbourhood $X^{/P}$, then $X\subseteq Z(m)$.\end{thm}
To prove such a theorem, note first that Proposition~\ref{injj} and Lemma~\ref{cond} are general. Then modify the proof of Theorem~\ref{T: algebraicC} by using Kuga--Satake construction before applying $p$-adic isogeny theorem.  \\
\end{remark}

\section{Estimation of local and global intersections}\label{cha3}
The goal of this section is to establish the main decay estimation needed for proving \S\ref{T:interMain}. 
The main results will be stated in \S\ref{sec:overview}. %while the rest of the chapter will be occupied by technical computations which can be skipped upon first reading. 

\subsection{Overview} \label{sec:overview}
Let the setup be the same as \S\ref{sub:localintersections}. That is, $\tau:C\rightarrow \mathcal{A}_{2,\kk}^{\tor}$ is a non-constant morphism, whose generic point lies in the almost ordinary locus of $\mathcal{A}_{2,\kk}$. For a $P\in C(\kk)$, we choose a uniformizer $t$ so that $C^{/P}\simeq \kk[[t]]$.  The point $P$ is called almost-ordinary (\textit{resp}. supersingular, \textit{resp}. superspecial \textit{resp}. of bad reduction) if its image is so.

The main result of this chapter is the following key lemma, which is an analogue of \cite[Proposition 7.2.5]{MAT}:
\begin{lemma}\label{decayresult}
If the image of $C$ is generically almost ordinary and does not lie on any special divisor, and $m$ is a positive integer with $\left(\frac{m}{p}\right)=1$. Then there exists an absolute constant $0<\alpha<1$, depending only on $C$, such that 
$$\sum_{P\,\,\mathrm{  supersingular} }l_P(m)\leq \alpha C\cdot \overline{Z(m)}+o(m^{\frac{3}{2}}).$$
\end{lemma}
Note that the all supersingular intersections of $C$ and $\overline{Z(m)}$ are superspecial by Lemma~\ref{T:intspecial}, so we only need to bound the local intersections at superspecial points. In the rest of this section, $P$ will be a superspecial point. Recall from \S\ref{sub:localintersections}, we have defined $$
    L''_{n,P}:=\left\{s\in L''_P\,|\,s \text{ lifts to } \Spf\kk[[t]]/(t^n)\right\},
$$
and showed that
$$
    l_P(m)=\sum_{n\geq 1}\#\left\{s\in L''_{n,P}|Q'(s)=m\right\}.
$$
Let $A$ be the intersection number of $C^{/P}$ with the supersingular locus at $P$, which is finite if $C$ is generically almost ordinary. Furthermore, let
$$A_n=\left\{\begin{aligned}
    &0,\,\,&n=-1,\\
    &A(1+p^2+...+p^{2n}),\,\,&n\geq 0.
\end{aligned}\right.$$
Unless otherwise specified, we will always assume that $P$ is superspecial. When $P$ is clear from the context, we abbreviate $L''_P$ and $L''_{n,P}$ as $L''$ and $L''_n$. \\%We The following definitions concern the different patterns of decays and the possible positions of $C^{/P}$.

\begin{definition}\label{def:3} Let $\omega$ denote a special endomorphism of $\mathscr{A}_P[p^\infty]$, $\Gamma$
denote a sub-lattice of $L''\otimes \mathbb{Z}_p$.
\begin{enumerate}
    \item $\omega$ is said to \textbf{lift to order} $N$, if it lifts to $\Spf\kk[[t]]/(t^N)$. $\omega$ is said to \textbf{lift exactly to order} $N$, if it lifts to $\Spf\kk[[t]]/(t^N)$, but not $\Spf\kk[[t]]/(t^{N+1})$.
    \item $\omega$ is said to \textbf{decay rapidly}, if $p^n\omega$ does not lift to $A_{n}+1$. It is said to \textbf{decay very rapidly}, if $p^n\omega$ does not lift to $A_{n-1}+\frac{A}{p+1}p^{2n}+1$.
    \item $\Gamma$ is said to \textbf{decay rapidly}, if every primitive element $\omega\in \Gamma$ decays rapidly. It is said to \textbf{decay very rapidly}, if every primitive element decays very rapidly.
    \item $L''$ is said to \textbf{decay completely}, if for every $n\geq 0$, there is an integer $N_n$ such that for every primitive $\omega\in L''$,  $p^n\omega$ does not lift to $N_n+1$. In other words, 
$L''_{N_n+1}\subseteq p^{n+1}L''$.
    \item $L''$ is said to \textbf{decay very rapidly at first step} with respect to $p$-splitting integers, if for any positive integer $m$ with $\left(\frac{m}{p}\right)=1$ and any $\omega\in L''$ with $Q'(\omega)=m$, $\omega$ does not lift to $\frac{A}{p+1}+1$.\\
\end{enumerate}
\end{definition} 

\iffalse\begin{definition}\label{def:position}
Let $\sigma_0$ be the absolute Frobenious on $\mathcal{A}_{2,\kk}^{/P}$. Regard $C^{/P}$ as a formal scheme which is finite over $\mathcal{A}_{2,\kk}^{/P}$ via $\tau^{/P}$. We say that $C^{/P}$ is in 
\begin{enumerate}
    \item \textbf{supersingular position}, if the tangent cones of
$C^{/P}$ and $\sigma_0(C^{/P})$ have the same induced reduced structure.
\item \textbf{special position}, if the tangent cones of
$C^{/P}$ and $\sigma_0^2(C^{/P})$ have the same induced reduced structure.
\item \textbf{general position}, if it is not in special position.
\end{enumerate}
\end{definition}
\fi

In order to do concrete computations, we use the coordinate in \S\ref{3.1} to make the identification $\widehat{\mathcal{A}_2}_P=W[[x,y,z]]$ and $\mathcal{A}_{2,\kk}^{/P}=\kk[[x,y,z]]$. Recall from Proposition~\ref{SSstratat} that in $\kk[[x,y,z]]$, the germ of the almost ordinary locus is given by $xy+\frac{z^2}{4\epsilon}=0$. So we can always assume that  $\tau^{/P}$ is given by
\begin{equation}\label{eq:curve}
    \tau^{/P}:\left\{\begin{aligned}
             &x= \beta^{-1}t^a+o(t^a),\\
             &y= -\beta t^b+o(t^b),\\
             &z=2{\lambda}t^{c}+o(t^c),\\
    \end{aligned}\right.
\end{equation}
with $\beta\neq 0$ and $a+b=2c$. {}{We will say that $C^{/P}$ is of \textit{supersingular/special/general position} according to the numerical condition listed in the following table (the definition can be made intrinsic to $C^{/P}$, but we won't need this)}: 

%Note that $C^{/P}$ is in supersingular position if and only if its tangent cone with induced reduced structure is one of the $p+1$ tangent lines given in Remark~\ref{sslocus}. It is also easy to see that $C^{/P}$ is in special position if and only if its tangent cone with induced reduced structure is one of the $p^2+1$ tangent lines with tangent directions $(1,0,0), (0,1,0)$ and $(\alpha^{-1} ,-\alpha ,2\lambda )$ where $\alpha\in \bF_{p^2}^*$. Therefore, the positions of $C^{/P}$ translate to numerical conditions on $a,b,c$ and $\beta$, which can be summarized in the following Table~\ref{table:1}:
\begin{table}[!h]
\centering
\begin{tabular}{ |p{2.2cm}||p{6cm}|}
 \hline
Positions & Numerical conditions \\
 \hline
Supersingular   & $a,b,c$ distinct or  $a=b=c$ and  $\beta/\lambda\in \bF_p^*$. \\
\hline
\hfil Special &  $a,b,c$ distinct or $a=b=c$ and $\beta\in \mathbb{F}_{p^2}^{*}$  \\\hline
\hfil General   &$a=b=c$ and $\beta\notin \mathbb{F}_{p^2}$\\
 \hline
\end{tabular}
\caption{}
\label{table:1}
\end{table}
Note that if $C^{/P}$ is in the supersingular position, then it is also in the special position. 

\begin{theorem}\label{Thm:localdensity} Suppose that $C$ does not lie on any $\overline{Z(m)}$, then
\begin{enumerate}
    \item\label{Thm:localdensity1}  $L''$ decays completely.
    \item \label{Thm:localdensity2}   
$C^{/P}$ is in special position but not in supersingular position, then $L''\otimes \mathbb{Z}_p$ admits a saturated rank 4 submodule $\Gamma_0$ which decays rapidly, inside which there is a saturated rank 2 submodule $\Gamma_1$ which decays very rapidly. 
    \item\label{Thm:localdensity3}   For all other positions of $C^{/P}$, the whole lattice $L''\otimes \mathbb{Z}_p$ decays rapidly, and there is a saturated rank 2 submodule $\Gamma_1$ which decays very rapidly. Furthermore, $L''$ decays very rapidly at first step with respect to $p$-splitting integers.
\end{enumerate}
\end{theorem}
%The statement for the submodules which decay rapidly and very rapidly in (\ref{Thm:localdensity3}) will not be used, but we still state and prove it for completeness (and also in compatible with previous literature \cite{MAT} and \cite{MST22}). 

\begin{remark}
The proof of Lemma~\ref{decayresult} mainly relies on Theorem~\ref{Thm:localdensity}. Besides that, some relatively easy local density computations are also needed.
\end{remark}

\subsection{Description of $F_\infty$ and $\tilde{\omega}$ over the non-ordinary stratum}\label{alphabeta} In this section we deduce an explicit description of the matrix $F_\infty$ (over the non-ordinary stratum), which we have defined in \S\ref{cha4}. We begin by fixing the notation:
\begin{itemize}
    \item The symbols $F, \omega_i, \tilde{\omega}, \epsilon, \lambda$, \textit{etc} will be consistent with that of \S\ref{expdL} and \S\ref{cha4}. 
    \item The symbols  $A$, $C$, $\tau$, $C^{/P}$ and $\tau^{/P}$ will be consistent with \S\ref{sec:overview}.
    \item {}{The notation $x^{[i]},y^{[i]},z^{[i]}$ and $M^{[i]}$ for a matrix $M$ is explained in \S\ref{subsec:Notation}: for example $x^{[i]}=x^{p^i}$.}
    \item {}{For a series $f\in W[[x,y,z]]$ or $\kk[[x,y,z]]$, we write $f\cdot C^{/P}$ for the degree of $\tau^{/P*}(f)$ as a series in $t$. It is the intersection number of $C^{/P}$ with the vanishing locus of $f$. }
\end{itemize}
\begin{definition}
Let $n\geq 0$. An \textbf{admissible index set} is an ordered set $\mathbf{I}=\{i_1<...<i_{2n}\}$ with $2n$ non-negative integers, such that $i_{2k}+i_{2k+1}\equiv 1 (\mod 2)$ for $k=1,2,...,n-1$. An \textbf{extended admissible index set} is one of the following:\begin{enumerate}
    \item  a 2-tuple $(i,\mathbf{I})$ such that $i\in \bZ_{\geq 0}$, and $\mathbf{I}$ is either empty or an non-empty admissible index set such that $\min \mathbf{I}>i$ and $\min \mathbf{I}+i\equiv1(\mod 2)$, 
    \item a 2-tuple $(\mathbf{I},j)$ such that $j\in \bZ_{\geq 0}$, and $\mathbf{I}$ is either empty or an non-empty admissible index set such that $\max \mathbf{I}<j$ and $\min \mathbf{I}+j\equiv1(\mod 2)$, \item a 3-tuple $(i,\mathbf{I},j)$ such that $i,j\in \bZ_{\geq 0}$ and $i<j$, such that $(i,\mathbf{I})$ and $(\mathbf{I},j)$ are both extended admissible index sets. \end{enumerate} 
\end{definition}
%\begin{definition}The \textbf{length} of an admissible index set or an extended admissible index set is defined as: $|\mathbf{I}|=|(i,\mathbf{I})|=\#\mathbf{I}/2$ and $|(\mathbf{I},j)|=|(i,\mathbf{I},j)|=1+ \#\mathbf{I}/2$.\end{definition}
For a pair of natural numbers $(i,j)$, we define $\alpha_{i,j}\in W[[x,y,z]]$ to be  $\alpha_{i,j}=x^{[i]}y^{[j]}+x^{[j]}y^{[i]}+\frac{z^{[i]}z^{[j]}}{2\epsilon}$. %We also let $\beta_0\in W[[x,y,z]]$ be $xy+\frac{z^2}{4\epsilon}$. 
We have $A=\alpha_{0,1}\cdot C^{/P}=\deg\tau^{/P*}({\alpha_{0,1}})$, {}{which is the intersection number of $C^{/P}$ with the vanishing locus of $\alpha_{0,1}$}. More generally, for an admissible index set $\mathbf{I}=\{i_1<...<i_{2n}\}$, we define $\alpha_{\mathbf{I}}=\alpha_{i_1,i_2}...\alpha_{i_{2n-1},i_{2n}}$.
Besides, we set
 $$D_{i,j}= \begin{bmatrix}
  (-1)^{i+j} & \frac{(-1)^{i+1}}{\lambda} \\
   (-1)^j\lambda & -1 \\
  \end{bmatrix}, B= \begin{bmatrix}
  \lambda y & -y \\ 
    \lambda x& -x \\
    \frac{\lambda z}{2\epsilon}& -\frac{z}{2\epsilon}
  \end{bmatrix} 
  ,\,\,E=\begin{bmatrix}
   \frac{x}{\lambda} & \frac{y}{\lambda} & \frac{z}{\lambda} \\
   x & y & z \\
  \end{bmatrix},$$
 and define $V_{\mathbf{I}}=E^{[i_1]}B^{[i_2]}...E^{[i_{n-1}]}B^{[i_n]}$ when $\mathbf{I}\neq \emptyset$, and $V_{\mathbf{I}}=I$ when $\mathbf{I}=\emptyset$. For $n\geq 0$, let \begin{equation}
     \label{eq:XY} X_n=\sum_{\substack{\mathbf{I} \text{ admissible}\\ |\mathbf{I}|=2n}}V_{\mathbf{I}}, \,\,Z_n=\sum_{\substack{(i,\mathbf{I}) \text{ admissible}\\ |\mathbf{I}|=2n}}B^{[i]}V_{\mathbf{I}} \end{equation}
     be the $2\times 2$ and $3\times 2$ matrices. For $n\geq 1$, let 
    \begin{equation}\label{eq:ZW} 
Y_n=\sum_{\substack{(\mathbf{I},j) \text{ admissible}\\ |\mathbf{I}|=2n-2}}V_{\mathbf{I}}E^{[j]},\,\,W_n=\sum_{\substack{(i,\mathbf{I},j) \text{ admissible}\\ |\mathbf{I}|=2n-2}}B^{[i]}V_{\mathbf{I}}E^{[j]}.
    \end{equation}
 be the $2\times 3, 3\times 3$ matrices. We also define $Y_0=\mathrm{O}_{2\times 3}, W_0=I_{3\times 3}$.  
\begin{proposition}
\label{T:1} Over the non-ordinary stratum, we have the following:\begin{enumerate}
 \item \label{T:12}If $\mathbf{I}\neq \emptyset$ is admissible, then $V_{\mathbf{I}}={(-2)^{|\mathbf{I}|/2-1}}\alpha_{\mathbf{I}}D_{\min\mathbf{I}, \max\mathbf{I}}$, 
    \item\label{T:11} For $n\geq 0$, let $F_n=2^{-n}\begin{bmatrix}
X_n & Y_n\\
Z_n & W_n
\end{bmatrix}$, then  $F_\infty=\sum_{n\geq 0}p^{-n}F_n$.
\end{enumerate} 
 
\end{proposition} 
\proof Over the non-ordinary stratum, we have $$F_\infty=\prod_{i\geq 0}(I+F^{[i]})= \prod_{i\geq 0}\left(I+\begin{bmatrix}0& \frac{1}{2p}E^{[i]}\\ B^{[i]}&0\end{bmatrix}\right).$$
 It is clear that $F_\infty$ is the sum over all possible matrices of the form
 $R_{i_1,...,i_n}:= \prod_{k=1 }^n\begin{bmatrix}0& \frac{1}{2p}E^{[i_k]}\\
 B^{[i_k]}&0\end{bmatrix}$,
 where $0\leq i_{1}<i_2<...<i_{n}$. Straightforward computation yields
\begin{equation*}
 R_{i_1,...,i_n}=\left\{
             \begin{array}{lr}
             \begin{bmatrix}\frac{1}{(2p)^k}E^{[i_1]}B^{[i_2]}...E^{[i_{n-1}]}B^{[i_n]}& 0\\ 0&\frac{1}{(2p)^k}B^{[i_1]}E^{[i_2]}...B^{[i_{n-1}]}E^{[i_n]}\end{bmatrix} , n=2k  &  \\
             \\
             \\
             \begin{bmatrix}0& \frac{1}{(2p)^{k+1}}E^{[i_1]}B^{[i_2]}...B^{[i_{n-1}]}E^{[i_n]}\\ \frac{1}{(2p)^{k}}B^{[i_1]}E^{[i_2]}...E^{[i_{n-1}]}B^{[i_{n}]}& 0\end{bmatrix}, n=2k+1 &  
             \end{array}
\right. ,
\end{equation*}
Note that $D_{i,j}D_{k,l}= 0$ if $j+k\equiv 0 (\mod 2)$, otherwise it is $-2D_{il}$. So we have 
\begin{equation*}
 E^{[i_1]}B^{[i_2]}...E^{[i_{2n-1}]}B^{[i_{2n}]}=\left\{
             \begin{array}{lr}
             (-2)^{n-1}\alpha_{i_1,i_2}...\alpha_{i_{2n-1},i_{2n}} D_{i_1,i_{2n}}, \,\,\, i_{2k}+i_{2k+1}\equiv 1(\mod 2)\\
             \\
             \\
             0, \,\,\,\text{otherwise}&  
             \end{array}
\right.
\end{equation*}
Therefore (\ref{T:12}) and (\ref{T:11}) holds. 
$\hfill\square$

\subsubsection{}Let $\omega\in L''\otimes\mathbb{Z}_p$ be a special endomorphism. Under basis $\{\omega_j\}_{j=1}^5$, we can write $\omega=\sum_{j=1}^5 a_j\omega_j$ with $a_j\in \mathbb{Z}_p$. Express $a_j=\sum_{k\geq 0}a_{j,k}p^k$, where $a_{j,k}$ are Teichmüller liftings of elements in $\mathbb{F}_p$. For $n\geq 1$, let $F_n$ be the matrix in Proposition~\ref{T:1}(\ref{T:11}).  Let $G^i_{n}(\omega)\in W[[x,y,z]]$ be defined as  
\begin{equation}\label{eq:vector11}
  \begin{bmatrix}
     G^1_{n}(\omega)\\
     G^2_{n}(\omega)\\
     \vdots\\
     G^5_{n}(\omega)
  \end{bmatrix}:=  \sum_{k\geq 0} F_{n+k}\begin{bmatrix}
     \alpha_{1,k}\\
     \alpha_{2,k}\\
     \vdots\\
    \alpha_{5,k}
  \end{bmatrix}.
\end{equation}
For $\star=$"$\even$" or "$\odd$" and $n\geq 1$, we also define  $X_{n,\star},Z_{n,\star}$ (\textit{resp}. $Y_{n,\star}, W_{n,\star}$) by adding one more restriction ``$\min \mathbf{I}$ is $\star$'' (\textit{resp}. ``$i$ is $\star$'') to the summation in (\ref{eq:XY}) and (\ref{eq:ZW}). We then let $F_{n,\star}:=2^{-n}\begin{bmatrix}X_{n,\star}&Y_{n,\star}\\
Z_{n,\star}&W_{n,\star}
\end{bmatrix}$. Similarly, define $G^i_{n,\star}(\omega)\in W[[x,y,z]]$ by replacing $F_{n+k}$ in (\ref{eq:vector11}) by $F_{n+k,\star}$. Note that $G^i_{n,\star}(\omega_j)$ is just the $(i,j)$-th entry of $F_{n,\star}$, and $ G_{n,\star}^i(\omega)=\sum_{j=1}^5\sum_{k\geq 0}a_{j,k}G_{n+k,\star}^i\left(\omega_j\right)$. 

If we express $\tilde{\omega}$ as $\sum_{j=1}^5 \tilde{a}_j\omega_j$, then $\tilde{a}_j= \sum_{n\geq 1}p^{-n}G_{n}^j(\omega)+O(1)$ ({}{Here the meaning of $O(1)$ is explained in \S\ref{subsec:Notation}}). Therefore $G_{n}^j(\omega)$ controls the deformation of $\omega$. The most important property of $G_{n}^j(\omega)$ is its recursive behavior, from which we realize that the deformation of $\omega$ is already controlled by $G_{1,\text{even}}^1(\omega)$.

\begin{proposition}[Recursive formula]\label{T:gdecay}~\begin{enumerate}
 \item \label{it:recurs1}$\begin{bmatrix}
    G_{1}^1(\omega)\\
    G_{1}^2(\omega)
    \end{bmatrix}=G_{1,\even}^1(\omega)\begin{bmatrix}
   1\\
   \lambda
   \end{bmatrix}  + G_{1,\even}^1(\omega)^{[1]}\begin{bmatrix}
   1\\
   -\lambda
   \end{bmatrix}  $, $\begin{bmatrix}
G_{1}^3(\omega)\\
     G_{1}^4(\omega)\\
      G_{1}^5(\omega)
   \end{bmatrix} =2\lambda\sum_{i\geq 0}(-1)^{i}G_{1,\even}^1(\omega)^{[i+1]}\begin{bmatrix}
     y^{[i]}\\
     x^{[i]}\\
     \frac{z^{[i]}}{2\epsilon}   \end{bmatrix}$.
  %   \item For $n>1$, we have \begin{equation*} \begin{bmatrix}
 %   G_{n}^1(\omega)\\
  %  G_{n}^2(\omega)
  %  \end{bmatrix}=(-1)^{n-1}\sum_{{|\mathbf{I}|=2n-2}}(-1)^{\max\mathbf{I}+1} \alpha_\mathbf{I}\begin{bmatrix}
   %  (-1)^{\min\mathbf{I}}\\
   % \lambda         \end{bmatrix}G^1_{1,\even}(\omega)^{[\max\mathbf{I}+1]}.\\
%\end{equation*}\begin{equation*}
 %%  \begin{bmatrix}
   % G_{n}^3(\omega)\\
    % G_{n}^4(\omega)\\
     % G_{n}^5(\omega)\end{bmatrix}=2\lambda(-1)^{n} \sum_{|\mathbf{I}|=2n-2}(-1)^{\max\mathbf{I}+1} \alpha_\mathbf{I}\begin{bmatrix}    y^{[i]}\\    x^{[i]}\\     \frac{z^{[i]}}{2\epsilon}   \end{bmatrix}G^1_{1,\even}(\omega)^{[\max\mathbf{I}+1]}.\end{equation*}
\item\label{it:recurs2} For $n>1$, we have  $$\begin{bmatrix}
    G_{n}^1(\omega)\\
    G_{n}^2(\omega)\\
    G_{n}^3(\omega)\\
    G_{n}^4(\omega)\\
    G_{n}^5(\omega)
    \end{bmatrix}=\sum_{|\mathbf{I}|=2n-2}(-1)^{\max\mathbf{I}+n}\alpha_{\mathbf{I}}G^1_{1,\even}(\omega)^{[\max\mathbf{I}+1]}\begin{bmatrix}
     (-1)^{\min\mathbf{I}}\\
    \lambda        \\
     -2\lambda y^{[i]}\\
     -2\lambda  x^{[i]}\\
     -z^{[i]}/\lambda
    \end{bmatrix}    $$
 \end{enumerate} 
\end{proposition}
\proof  We fix an integer $m\geq 1$. Let $\mathbf{U}_n$ and $\mathbf{V}_n$ be the vectors representing the first row of $X_{n,\even}$ and $Y_{n,\even}$. Clearly, $X_{n,\even}=\begin{bmatrix}
    1\\
    \lambda
    \end{bmatrix}\mathbf{U}_n$, $Y_{n,\even}=\begin{bmatrix}
    1\\
    \lambda
    \end{bmatrix}\mathbf{V}_n$, and 
    $
    (X_n,Y_n)= \begin{bmatrix}
     1\\
     \lambda
     \end{bmatrix}(\mathbf{U}_n,\mathbf{V}_n)+ \begin{bmatrix}
     1\\
     -\lambda
     \end{bmatrix}(\mathbf{U}_n,\mathbf{V}_n)^{[1]}$. This implies the first identity of (\ref{it:recurs1}). For $n>m$, we have      
\begin{align*}
    X_n&= \sum_{|\mathbf{I}|=2n}V_{\mathbf{I}}=\sum_{|\mathbf{I}|=2n-2m}V_{\mathbf{I}}\sum_{\substack{|\mathbf{I}'|=2m\\\min\mathbf{I}'\text{ is odd}} } V_{\mathbf{I}'}^{[\max \mathbf{I}]}=\sum_{|\mathbf{I}|=2n-2}V_{\mathbf{I}}\sum_{\substack{|\mathbf{I}''|=2m\\\min\mathbf{I}''\text{ is even}} } V_{\mathbf{I}''}^{[\max \mathbf{I}+1]}\\&= \sum_{|\mathbf{I}|=2n-2m}V_{\mathbf{I}} X_{m,\text{even}}^{[\max \mathbf{I}+1]}=\sum_{|\mathbf{I}|=2n-2m}V_{\mathbf{I}}\left(\begin{bmatrix}
    1\\
    \lambda 
    \end{bmatrix}\mathbf{U}_{m}\right)^{\left[\max\mathbf{I}+1\right]}\\&=
    (-2)^{n-m}\sum_{|\mathbf{I}|=2n-2m}(-1)^{\max \mathbf{I}+1}\alpha_\mathbf{I}\begin{bmatrix}
    (-1)^{\min\mathbf{I}}\\
    \lambda
    \end{bmatrix}\mathbf{U}_m^{[\max\mathbf{I}+1]}.
\end{align*}
Similarly, 
\begin{align*}
    Y_n=   (-2)^{n-m}\sum_{|\mathbf{I}|=2n-2m}(-1)^{\max \mathbf{I}+1}\alpha_\mathbf{I}\begin{bmatrix}
    (-1)^{\min\mathbf{I}}\\
    \lambda
    \end{bmatrix}\mathbf{V}_m^{[\max\mathbf{I}+1]}.
    \end{align*}
Now we define $\max(i,\mathbf{I})=\max \mathbf{I}$ if $\mathbf{I}\neq \emptyset$ and $\max(i,\mathbf{I})=i$ if $ \mathbf{I}= \emptyset$. For $n\geq m$, similar argument as above shows that 
    \begin{align*}
    Z_n&= \sum_{|\mathbf{I}|=2n}B^{[i]}V_{\mathbf{I}}=\sum_{|\mathbf{I}|=2n-2m}B^{[i]}V_{\mathbf{I}}\left(\begin{bmatrix}
    1\\
    \lambda 
    \end{bmatrix}\mathbf{U}_{m}\right)^{\left[\max(i,\mathbf{I})+1\right]}\\&=\sum_{|\mathbf{I}|=2n-2m}\begin{bmatrix}
    y\\
    x\\
    \frac{z}{2\epsilon}
    \end{bmatrix}^{[i]}\begin{bmatrix}
    (-1)^i\lambda&-1
    \end{bmatrix}V_{\mathbf{I}}\begin{bmatrix}
    1\\
    (-1)^{\max(i,\mathbf{I})+1}\lambda 
    \end{bmatrix}\mathbf{U}_{m}^{\left[\max(i,\mathbf{I})+1\right]}\\&= 
  (-2)^{n-m+1}\lambda\sum_{|\mathbf{I}|=2n-2m}(-1)^{\max(i,\mathbf{I})+1} \alpha_\mathbf{I}\begin{bmatrix}
    y\\
    x\\
    \frac{z}{2\epsilon}
    \end{bmatrix}^{[i]}\mathbf{U}_{m}^{\left[\max(i,\mathbf{I})+1\right]},\\
    W_n&=
  (-2)^{n-m+1}\lambda\sum_{|\mathbf{I}|=2n-2m}(-1)^{\max(i,\mathbf{I})+1} \alpha_\mathbf{I}\begin{bmatrix}
    y\\
    x\\
    \frac{z}{2\epsilon}
    \end{bmatrix}^{[i]}\mathbf{V}_{m}^{\left[\max(i,\mathbf{I})+1\right]}.
\end{align*}
Then the reduction formula of $G^i_n(\omega)$ follows directly from the reduction of the corresponding matrices and the definition of $G^i_n(\omega)$. For example, for $n>1$,
\begin{align*}
   {2^{-n}}(X_n,Y_n)&= (-1)^{n-m}2^{-m}\sum_{|\mathbf{I}|=2n-2m}(-1)^{\max\mathbf{I}+1} \alpha_\mathbf{I}\begin{bmatrix}
   (-1)^{\min\mathbf{I}}\\
   \lambda
   \end{bmatrix}(\mathbf{U_m},\mathbf{V}_m)^{[\max\mathbf{I}+1]}\\&=(-1)^{n-m}\sum_{|\mathbf{I}|=2n-2m}(-1)^{\max\mathbf{I}+1} \alpha_\mathbf{I}\begin{bmatrix}
   (-1)^{\min\mathbf{I}}\\
   \lambda
   \end{bmatrix}[G_{m,\even}^1(\omega_1),...,G_{m,\even}^1(\omega_5)]^{[\max\mathbf{I}+1]}.
\end{align*}
Let $n>1$. For any $k\geq 0$, take $(n,m)$ to be $(n+k,1+k)$ in the above formula. We get
\begin{align*}
G_n^1(\omega)&=\sum_{j=1}^5\sum_{k\geq 0}a_{j,k}G_{n+k}^1(\omega_j)\\&=\sum_{j=1}^5\sum_{k\geq 0}a_{j,k}(-1)^{n-1}\sum_{|\mathbf{I}|=2n-2}(-1)^{\max\mathbf{I}+\min\mathbf{I}+1} \alpha_\mathbf{I}G_{1+k,\even}^1(\omega_j)^{[\max\mathbf{I}+1]}\\&=(-1)^{n-1}\sum_{|\mathbf{I}|=2n-2}(-1)^{\max\mathbf{I}+\min\mathbf{I}+1} \alpha_\mathbf{I}\sum_{j=1}^5\sum_{k\geq 0}a_{j,k}G_{1+k,\even}^1(\omega_j)^{[\max\mathbf{I}+1]}\\&=\sum_{|\mathbf{I}|=2n-2}(-1)^{\max\mathbf{I}+\min\mathbf{I}+n} \alpha_\mathbf{I}{G^1_{1,\even}(\omega)}^{[\max\mathbf{I}+1]}.
\end{align*}
This proves the first row of (\ref{it:recurs2}). Other formulas are obtained in a similar way, and is left to the readers. 
$\hfill\square$
\subsection{{Deformation along $C^{/P}$}}\label{5.2}
For $\omega\in L''\otimes\mathbb{Z}_p$, we define $\mathcal{D}(\omega):= G_{1,\text{even}}^1(\omega)(\mod p)\in \kk[[x,y,z]]$. Let $\deg(\omega,C^{/P}):=\mathcal{D}(\omega)\cdot C^{/P}= \deg(\tau^{/P*}\mathcal{D}(\omega))$, {}{which is the intersection number of $C^{/P}$ with $\mathcal{D}(\omega)$}. When $C^{/P}$ is clear from the context, we simply abbreviate it as {}{$\deg(\omega)$}. As noted before, $\mathcal{D}(\omega)$ controls the deformation of $\omega$: 
\begin{proposition}
\label{T:Decay} $p^{n}\omega$ lifts exactly to $A(1+p^2+... +p^{2n-2})+\deg(\omega)p^{2n}$.
%\begin{enumerate}   \item\label{T:Decay1} If $\hat{C}$ is generically supersingular, then $p\omega$ never decays.     \item\label{T:Decay2} If $\hat{C}$ is generically almost ordinary, then $   \hat{C}\cdot\mathcal{D}_{n+1}=A(1+p^2+... +p^{2n-2})+dp^{2n}$.   \item\label{T:Decay3}  If $\hat{C}$ is generically ordinary, and if $H_e\leq (1+p)B<H_{e+1}$, then    $$\hat{C}.\mathcal{D}_{n+1}=\left\{\begin{aligned}     &A(1+p^2+...+p^{2n-2})+dp^{2n},\,\,\,\,&n\leq e-1,\\&    B(1+p+...+p^{n-e-1})+Dp^{n-e},\,\,\,\,&n\geq e.\end{aligned}\right.$$Moreover, we always have $D\geq A(1+...+p^{2e-2})+dp^{2e}$. If furthermore $(1+p)B>H_e$, then $D=A(1+...+p^{2e-2})+dp^{2e}$.\end{enumerate}
\end{proposition}
We first need a lemma:
\begin{lemma}\label{less}
The minimal value of the set $\{ \alpha_{0,i}\cdot C^{/P}\}_{i\geq 1}$ is $A=\alpha_{0,1}\cdot C^{/P}$. 
\end{lemma}
%\begin{example}It is not always true that $A=\hat{C}\cdot\alpha_{0,1}< \hat{C}\cdot\alpha_{0,i}$ for $i> 1$. For example, take $\hat{C}$ to be $x=(\lambda-\frac{1}{4\epsilon})t, y=t,z=t$, then $A=\infty$, while $\hat{C}\cdot\alpha_{0,2}=p+1$. However, in this particular case, we do have $B=\hat{C}\cdot\beta_{0}< \hat{C}\cdot\alpha_{0,i}$.    \end{example}
\proof  This lemma is a consequence of \cite[Theorem 19]{Ogus01}, but we also give a direct computational proof.
If $A=\min\{(x^py)\cdot C^{/P},(xy^p)\cdot C^{/P},(z^{p+1})\cdot C^{/P}\}$, the lemma is obvious. Otherwise, we can write ${\tau}^{/P}$ as \begin{equation}\label{eq:4.12}
    \left\{\begin{array}{lr}
             x= \delta^{-1}(t^a+ f_{1}t^{a+1}+f_2t^{a+2}+...), &\\
             y= -\delta \epsilon (t^a+g_{1}t^{a+1}+g_2t^{a+2}+...),& \\
             z=2{\epsilon}(t^{a}+ h_{1}t^{a+1}+h_2t^{a+2}+...),\\
             \end{array}\right.
\end{equation}
for some $\delta\in \mathbb{F}_p^*$. We use bump induction: Set $f_0=g_0=h_0=1$. For system of elements $\{f_i,g_i,h_i\}_{i\geq 0}$, and $l,k\geq 0$, we define the following system of assertions \begin{itemize}
    \item $\mathscr{W}_k:2h_k-f_k-g_k=0$,
    \item $\mathscr{X}_k:f_k-g_k=0$,
    \item $\mathscr{A}^i_l:2\sum_{k=0}^{[\frac{l}{p^i}]} h_k^{p^i}h_{l-kp^i}-\sum_{k=0}^{[\frac{l}{p^i}]} f_k^{p^i}g_{l-kp^i}-\sum_{k=0}^{[\frac{l}{p^i}]} g_k^{p^i}f_{l-kp^i}=0$ 
    \item $ \mathscr{B}_l:\sum_{k=0}^l h_kh_{l-k}-\sum_{k=0}^l f_kg_{l-k}=0$.
\end{itemize} 

Furthermore, $ \mathscr{W}_{\leq n}$, $\mathscr{W}_{<\infty}$ are defined to be $\bigwedge_{k\leq n}\mathscr{W}_{k}$, $\bigwedge_{k\geq 1}\mathscr{W}_{k}$ (wedge means being simultaneously true). We introduce the similar notation for $\mathscr{X}$, $\mathscr{A}^i$, and $\mathscr{B}$. Since $C^{/P}$ is generically almost ordinary, $\mathscr{B}_{< \infty}$ is true. Moreover, $\mathscr{W}_0,\mathscr{X}_0,\mathscr{A}^i_0,\mathscr{B}_0$ are vacuously true, and one easily verifies that $\mathscr{A}^i_1\Leftrightarrow \mathscr{B}_1\Leftrightarrow \mathscr{W}_1$.  

Now the coefficient of the term $t^{ap^i+a+l}$ in $\alpha_{0,i}\cdot C^{/P}$ is 0 iff $\mathscr{A}^i_l$ is true. In particular, $\mathscr{A}^1_{\leq A-{a(p+1)}-1}$ is true. Now Lemma~
\ref{induction} (\ref{induction4}) shows that 
$\mathscr{A}^i_{\leq A-{a(p+1)}-1}$ is true for all $i\geq 1$. Therefore, for any $i>1$, we have $\alpha_{0,i}\cdot C^{/P}\geq ap^i+a+A-{a(p+1)}>A$. This concludes the lemma. $\hfill\square$
\begin{lemma}\label{induction}
Notation as in the proof of Lemma~\ref{less}. If $\mathscr{B}_{<\infty}$ is true, then the following hold:\begin{enumerate}
    \item \label{induction1}$\mathscr{W}_{\leq 2n}\Rightarrow \mathscr{X}_{\leq n}$,
    \item  \label{induction2}$ \mathscr{W}_{\leq 2n-1}\Rightarrow(\mathscr{W}_{\leq 2n}\Leftrightarrow \mathscr{X}_{\leq n}) $,
    \item \label{induction3}$\mathscr{W}_{\leq 2n}\Rightarrow \mathscr{W}_{\leq 2n+1}$,
    \item  \label{induction4}$\mathscr{A}^i_{\leq n}\Leftrightarrow \mathscr{W}_{\leq n}$.
\end{enumerate} 
\end{lemma}
\proof~\begin{enumerate}
    \item  We verify by direct computation that $\mathscr{W}_{\leq 2n}\wedge\mathscr{X}_{\leq n-1}\Rightarrow \mathscr{X}_n $. Then do induction.
\item If $\mathscr{W}_{\leq 2n-1}$ is true, then by (\ref{induction1}), we have $\mathscr{X}_{\leq n-1}$. We then verify by direct computation that $(\mathscr{W}_{2n}\Leftrightarrow \mathscr{X}_n)$.
\item If $\mathscr{W}_{\leq 2n}$ is true, then by (\ref{induction1}), $\mathscr{X}_{\leq n}$ is also true. It then follows that $\mathscr{W}_{2n+1}$ is true. \item We verify by direct computation that for any $0\leq r<p^i$, \begin{equation}\label{logic}
    \mathscr{W}_{\leq l}\wedge \mathscr{X}_{\leq l}\wedge(\bigwedge_{k<l} \mathscr{W}_{{k}p^i+r})\Rightarrow(\mathscr{A}^i_{lp^i+r}\Leftrightarrow \mathscr{W}_{lp^i+r}).
\end{equation}
We then do induction. If $n<p^i$, the statement is obvious. Now suppose that  (\ref{induction4}) is true for all $n<kp^i$, it suffices to show that it is true for all $kp^i+r, 0\leq r<p^i$. 
To this end, suppose that $\mathscr{A}^i_{\leq kp^i+r}$ is true for all $0\leq r<p^i$. By the induction hypothesis, $\mathscr{W}_{\leq kp^i-1}$ is true, hence (\ref{logic}) implies that  $\mathscr{W}_{kp^i+r}$ is true. Conversely, if $\mathscr{W}_{\leq kp^i+r}$ is true, then by (\ref{induction1}), $\mathscr{X}_{\leq [\frac{kp^i+r}{2}]}$ is true. In particular, $\mathscr{X}_{\leq k}$ is true, hence (\ref{logic}) implies $\mathscr{A}^i_{kp^i+r}$ for $0\leq r<p^i$.
$\hfill\square$\\\end{enumerate}

\noindent\textit{Proof of Proposition~\ref{T:Decay}}. Let $\mathcal{D}_n'(\omega):= G_{n}^1(\omega)(\mod p)\in \kk[[x,y,z]]$. By definition, we have $\mathcal{D}_{n+1}'(\omega)=\mathcal{D}'(p^n\omega)$. We can then use Proposition~\ref{T:gdecay} to express $\mathcal{D}_{n+1}'(\omega)$ in terms of $\mathcal{D}(\omega)$. As a consequence of Lemma~\ref{less}, there is a unique term in the set $\left\{\alpha_{\mathbf{I}}G^1_{1,\even}(\omega)^{[\max\mathbf{I}+1]}\big|\,\,|\mathbf{I}|=2n\right\}$ which has the minimal intersection number with $C^{/P}$. This minimal term is achieved exactly when $\mathbf{I}=\{0,1,...,2n-2,2n-1\}$, and the minimal intersection number is $N:=A(1+p^2+...+p^{2n-2})+\deg(\omega)p^{2n}$. Therefore, the intersection numbers of $\mathcal{D}'(p^n\omega)$ and $\mathcal{D}(p^n\omega)$ with $C^{/P}$ are both $N$. We then apply  Proposition~\ref{T:gdecay} and Lemma~\ref{T:PDlemma} to show that $p^n\omega$ lifts exactly to $N$.$\hfill\square$

\subsection{Proof of Theorem~\ref{Thm:localdensity}}Recall that we have $A=\alpha_{0,1}\cdot C^{/P}=\deg {\tau^{/P}}^*(\alpha_{0,1})$. There is an apparent dependence of the quantity $A$ upon $a,b,c$ and $\beta$:
\begin{lemma}\label{lm:Acase}The following are true:\begin{enumerate}
    \item\label{lm:Acase1} if $a,b,c$ are distinct, then $A=\max\{a,b\}+\min\{a,b\}p$.
    \item\label{lm:Acase2} if $a=b=c$ and $\beta/\lambda\notin \mathbb{F}_p^*$ , then $A=a(p+1)$.
    \item\label{lm:Acase3} if $a=b=c$ and $\beta/\lambda\in \mathbb{F}_p^*$, then $A>a(p+1)$. 
\end{enumerate}
\end{lemma} 
\proof Since $A=\deg{\tau^{/P*}}(\alpha_{0,1})$, easy computation gives (\ref{lm:Acase1}). 
In the case (\ref{lm:Acase2}) and (\ref{lm:Acase3}), we 
at least have $A\geq a(p+1)$. Moreover, $A>a(p+1)$ if and only if the coefficient of  $t^{a(p+1)}$ in ${\tau^{/P*}}(\alpha_{0,1})$ is 0, i.e., if and only if 
\begin{equation}\label{eq:todm}
    -\beta^{p-1}- \beta^{1-p}+\frac{(2\lambda)^{p+1}}{2\epsilon}=0.
\end{equation}
An easy computation then shows that (\ref{eq:todm}) holds exactly when $\beta/\lambda\in\bF_p^*$. 
$\hfill\square$

\begin{lemma}
\label{T:som}The following are true: 
\begin{enumerate}
    \item \label{T:som1} $p^n\omega_1$ and $p^n\omega_2$ lift exactly to $A_{n}$,
    \item \label{T:som2}$p^n\omega_3,p^n\omega_4$ and $p^n\omega_5$ lift exactly to $A_{n-1}+ap^{2n}$,$A_{n-1}+bp^{2n}$ and $A_{n-1}+cp^{2n}$ respectively,
    \item\label{T:som3} $L''_n=L''$ for $n\leq \min\{a,b\}$.
\end{enumerate}
\end{lemma}
\proof Proposition~\ref{T:Decay} implies that $\omega$ lifts exactly to $\deg(\omega)$. It is easy to see that $\deg(\omega_1)=\deg(\omega_2)=A$, $\deg(\omega_3)=a,\deg(\omega_4)=b$, and $\deg(\omega_5)=c$. So (\ref{T:som1}) and (\ref{T:som2}) follows. By (\ref{T:som1}) and (\ref{T:som2}), $\omega_i$ at least lifts to $\kk[[t]]/(t^{\min\{a,b\}})$, and this implies (\ref{T:som3}). $\hfill\square$\\[10pt]
We now begin a case by case study of the patterns of decay according to Table~\ref{table:1}. According to Lemma~\ref{lm:Acase}, 
the case where $a,b,c$ are distinct, and the case where $a=b=c$ and $\beta\lambda\in \bF_p^*$ assert different effects on the quantity $A$. We will further divide the supersingular case into two subcases:
\begin{enumerate}
    \item Supersingular case I: $a,b,c$ are disctinct,
    \item Supersingular case II: $a=b=c$ and $\beta/\lambda\in \bF_p^*$.
\end{enumerate}
\subsubsection{Supersingular case I and general case} Write $\omega$ as $\sum_{j=1 }^5 a_jp^{n_j}\omega_j$ where $a_j\in \mathbb{Z}_p^*$, then
\begin{lemma}\label{SSS}
$\deg(\omega)=\min _{j=1,2,...,5}\{\deg(p^{n_{j}}\omega_{j})\}$.
\end{lemma}
\proof By Lemma~\ref{T:som}, $\deg(p^{n_j}\omega_j)$ is $A_{n_j}$ for $j=1,2$ and $A_{n_j-1}+\gamma p^{2n_j}$, where $\gamma=a,b,c$ for $i=3,4,5$. Note that $A_{n_j}$ can never be equal to $A_{n_j-1}+\gamma p^{2n_j}$. So looking at corresponding lowest degrees of $\tau^{/P*}\mathcal{D}_1$, we have $\deg(\omega)=\min\{\deg(\sum_{j=1}^2a_jp^{n_j}\omega_j),\deg(\sum_{j=3}^5a_jp^{n_j}\omega_j) \}$. 

On one hand,  $\deg(\sum_{j=1}^2a_jp^{n_j}\omega_j)=\min\{ \deg(p^{n_1}\omega_1),\deg(p^{n_2}\omega_2)\}$, since even if $n_1=n_2$, the terms of lowest degrees from $p^{n_1}\omega_1$ and $p^{n_2}\omega_2$ can not cancel because $\lambda\notin \mathbb{Z}_p$. On the other hand, we claim  $\deg(\sum_{j=3}^5a_jp^{n_j}\omega_j)=\min\{ \deg(p^{n_3}\omega_1), \deg(p^{n_4}\omega_4),\deg(p^{n_5}\omega_5)\}$.  This is trivially true when $a,b,c$ are distinct. But for $a=b=c$ and $n_3=n_4=n_5$, if $a_3,a_4,a_5\in \mathbb{Z}_p^*\cup\{0\}$, not simultaneously 0, are such that $\deg(\sum_{j=3}^5a_jp^{n_j}\omega_j)> \deg(p^{n_3}\omega_1)$, then $$\overline{a_3}\beta^{-1}-\overline{a_4}\beta+2\overline{a_5}\lambda=0$$
and this shows that $\beta\in \mathbb{F}_{p^2}$. Therefore the claim is also true in the general case. $\hfill\square$

\subsubsection{Supersingular case II}
Let $\beta/\lambda= \delta$, we can rewrite $\tau^{/P}$ in the form of (\ref{eq:4.12}). Let $\omega_3'=\delta^2\epsilon\omega_3-\omega_4-\delta\omega_5$ and  $\omega_4'=\delta^2\epsilon\omega_3+\omega_4$. Let $e=A-a(p+1)$. 
\begin{lemma}\label{WW}
Notation being the same as in Lemma~\ref{less}, $e\geq 2$ is even. Moreover we have $\mathscr{W}_{\leq e-1}$, $\neg\mathscr{W}_{e}$, $\mathscr{X}_{\leq e/2-1}$ and $\neg\mathscr{X}_{e/2}$.
\end{lemma}
\proof Note that we have $\mathscr{A}^1_{\leq e-1}$ but $\neg \mathscr{A}^1_{e}$. If $e=2n+1$ is odd, then we have $\mathscr{A}^1_{\leq 2n}$, by Lemma~\ref{induction} (\ref{induction4}) we see $\mathscr{W}_{\leq 2n}$, then by Lemma~\ref{induction} (\ref{induction3}) we have $\mathscr{W}_{\leq 2n+1}$. By (\ref{induction4}) again we have $\mathscr{A}^1_{2n+1}$, a contradiction.  $\mathscr{W}_{\leq e-1}$, $\neg\mathscr{W}_{e}$ is a consequence of (\ref{induction4}), while $\mathscr{X}_{\leq e/2-1}$ is a consequence of (\ref{induction1}). If we also have $\mathscr{X}_{e/2}$ then we will have $\mathscr{W}_{e}$ by (\ref{induction2}), so  $\neg\mathscr{X}_{e/2}$.$\hfill\square$\\[5pt]
The following two lemmas are immediate from Lemma~\ref{WW}. Moreover, the proof of Lemma~\ref{WWWW} is similar to Lemma~\ref{SSS}. We left the proof to the readers. 

\begin{lemma}\label{WWW}
$p^n\omega_3'$ lifts exactly to $A_{n-1}+(a+e)p^{2n}$, $p^n\omega_4'$ lifts exactly to $A_{n-1}+(a+e/2)p^{2n}$.
\end{lemma}

\begin{lemma}\label{WWWW}
 Write $\omega=a_1p^{n_1}\omega_1+a_2p^{n_2}\omega_2+a_3p^{n_3}\omega'_3+a_4p^{n_4}\omega_4'+a_5p^{n_5}\omega_5$ where $a_j\in \mathbb{Z}_p^*$, then $$\deg(\omega)=\min _{j=1,2,...,5}\{\deg(p^{n_{1}}\omega_{1}),\deg(p^{n_{2}}\omega_{2}),\deg(p^{n_{3}}\omega'_{3}),\deg(p^{n_{4}}\omega_{4}'),\deg(p^{n_{5}}\omega_{5})\}.$$
\end{lemma}
\subsubsection{Special but non-supersingular case}\label{sub:caseS}
In this case, there exists $\alpha \in \bF_p^*$ and $\gamma\in \bF_p$ with $\alpha+\gamma^2\epsilon\in (\bF_p^*)^2$, such that $-\beta+ \alpha\beta^{-1}= 2\lambda\gamma$. Let $\eta=\alpha\omega_3+\omega_4+\gamma\omega_5$. $\eta$ may resist decay, but we at least have expected decay for the submodule $\Gamma_0=\mathbb{Z}_p\omega_1\oplus \mathbb{Z}_p\omega_2\oplus\mathbb{Z}_p\omega_3\oplus\mathbb{Z}_p\omega_5$. The following lemma is similar to Lemma~\ref{SSS}, whose proof is left to the readers: 
\begin{lemma}\label{SST}
Write $\omega\in \Gamma_0$ as $a_1p^{n_1}\omega_1+a_2p^{n_2}\omega_2+a_3p^{n_3}\omega_3+ a_5p^{n_5}\omega_5$ where $a_j\in \mathbb{Z}_p^*$, then $$\deg(\omega)=\min _{j=1,2,3,5}\{\deg(p^{n_{j}}\omega_{j})\}.$$
\end{lemma}\,\\
\textit{Proof of Theorem~\ref{Thm:localdensity}}. We first prove (\ref{Thm:localdensity2}) and (\ref{Thm:localdensity3}), and finally prove (\ref{Thm:localdensity1}). 

For (\ref{Thm:localdensity2}), we use the observations made in  \S\ref{sub:caseS}. Let $\Gamma_0$ be the submodule defined there, and $\Gamma_1=\mathbb{Z}_p\omega_3\oplus \mathbb{Z}_p\omega_5$, then Lemma~\ref{T:som} and Lemma~\ref{SST} implies that $\Gamma_0$ decays rapidly and $\Gamma_1$ decays very rapidly. 

For (\ref{Thm:localdensity3}), we separate cases. We first treat supersingular case I and general case. Furthermore, in supersingular case I we suppose without loss of generality that $a\leq b$. Let $\Gamma_1=\mathbb{Z}_p\omega_3\oplus \mathbb{Z}_p\omega_5$. Then Lemma~\ref{T:som} and Lemma~\ref{SSS} implies that $L''\otimes \bZ_p$ decays rapidly and $\Gamma_1$ decays very rapidly.  In the following, let $m$ be such that $\left(\frac{m}{p}\right)=1$, and let $x_3$ and $x_5$ be the component of $\omega$ at $\omega_3$ and $\omega_5$. Suppose $n>a$, 
Lemma~\ref{SSS} again implies that for any $\omega\in L''_n\otimes \mathbb{Z}_p$, one has $p|x_3$. Now by (\ref{eq:Pairingthing}), the reduction modulo $p$ of the equation $Q(\omega)=m$ can be explicitly written as $\epsilon \overline{x_5}^2=\overline{m}$. Since $\epsilon$ is not a square in $\mathbb{F}_p$ while $m$ is, the equation admits no solution. This implies that $L''$ decays very rapidly at first step with respect to $p$-splitting integers. 

For supersingular case II, let $\Gamma_1=\mathbb{Z}_p\omega_4'\oplus \mathbb{Z}_p\omega_5$, then Lemma~\ref{T:som} and Lemma~\ref{WWWW} implies that $L''\otimes \bZ_p$ decays rapidly and $\Gamma_1$ decays very rapidly. The proof of the statement about $L''$ decaying very rapidly at first step with respect to $p$-splitting integers is similar to the supersingular case I and the general case. 

Finally we treat (\ref{Thm:localdensity1}). In the supersingular or general case, we already see that $L''\otimes \bZ_p$ decays rapidly. This implies that $L''$ decays completely. For the special but non-supersingular case,  we shall prove the existence of an $N_n$ for each $n\geq 0$ such that $L''_{N_n+1}\subseteq p^{n+1}L''$. {}{Suppose for some $n$ there is no such $N_n$, so for every $j$ there exists an $e_j\in L''_{j+1}-p^{n+1}L''$. Possibly scaling $e_j$, there exists an $\eta$ such that for every $j$, there exists a $v_j\in p^{n+1}\Gamma_0$ such that $e_j=p^n\eta+ v_j$. However, this implies that for any $j'>j$, $v_j-v_{j'}\in \Gamma_0\cap L''_{j+1}\otimes\bZ_p$. Since $\Gamma_0$ decays rapidly, $\{v_j\}$ converges in the $p$-adic topology. Therefore $e_j$ converges to an element $e\in L''\otimes \bZ_p$ that never decays.} By Theorem~\ref{T: algebraicC}, the existence of $e$ implies that $C$ lies on a special divisor. Contradiction. $\hfill\square$
\subsection{Proof of Lemma~\ref{decayresult}} Throughout $C$ will be a smooth curve mapping to $\mathcal{A}_{2,\kk}^{\tor}$ whose image lies generically in the almost ordinary strutum, but does not lie on any special divisor. The letter $m$ will always be a positive integer with $\left(\frac{m}{p}\right)=1$. Recall from \S\ref{C:intersection} that $q_L(m)$ is defined to be $q_L(m,0)$, the Fourier coefficients of the Eisenstein series attached to $L$. We define the global intersection of $\overline{Z(m)}$ and $C$ at a supersingular point $P$ as
$$g_P(m):= \frac{A}{p^2-1}|q_L(m)|.$$
The notation is clearly chosen to contrast the local intersections $l_P(m)$. However, note that the local intersection $l_P(m)$ is defined at any point, while the global intersection is only defined at supersingular points. The global intersection is related to the intersection number $C\cdot\overline{{Z}(m)}$ in the following sense: \begin{proposition}\label{T:intersection}
Let $C$ be generic almost ordinary, then
$$\sum_{P\in C\,\,\mathrm{supersingular}}g_P(m)=C\cdot\overline{{Z}(m)}+o(m^{\frac{3}{2}}).$$
\end{proposition}
\proof Let $\mathcal{L}$ be the Hodge line bundle. By Lemma~\ref{intersEisen} we have $$C\cdot \overline{Z(m)}=|q_L(m)|(C\cdot\mathcal{L})+o(m^{\frac{3}{2}}).$$
By \cite[Theorem 4.5.4, Theorem 6.2.3]{Box15} the supersingular locus is cut out by a section of $\mathcal{L}^{p^{2}-1}$ inside the almost ordinary locus.$\hfill\square$\\[10pt]
\textit{Proof of Lemma~\ref{decayresult}}: 
We first fix the supersepcial point $P$. {}{Following \S\ref{sub:EsQl}, for each $n$, let $L'_n$ be a $\bZ$-sublattice of $L'$ such that  $L''_n\subseteq L'_n\subseteq L''_n\otimes \bQ$, $L'_n\otimes \bZ_p =L''_n\otimes \bZ_p$, and $L'_n\otimes \bZ_l=L'\otimes\bZ_l$ for $l\neq p$.} We have 
\begin{equation}\label{eq:sdfsdf123}
    l_P(m)\leq \sum_{n\geq 1}\#\left\{s\in L'_{n}|Q'(s)=m\right\}.
\end{equation}

 To compute the number on the right hand side of (\ref{eq:sdfsdf123}), we will take $L'''$ in \S\ref{sub:EsQl} to be various $L'_n$ and use Eisenstein series. Decompose the theta series of $L'_n$ as $E+G$ as \textit{loc.cit}.  We first compute $\delta(p,L'_n,m)$ for $\left(\frac{m}{p}\right)=1$.
Following \S\ref{expdL}, {}{we write an element of $v
\in L''_n\otimes \bZ_p=L'_n\otimes \bZ_p$ as $v=u'\omega_1+w'\omega_2+ x'\omega_3+y'\omega_3+z'\omega_5$. Then the explicit form (\ref{eq:Pairingthing}) gives $Q'(v)=x'y'+\epsilon {z'}^2+p{w'}^2-p\epsilon {u'}^2$.}  By Lemma~\ref{T:Hanke} it suffices to count the number of elements in the set $\{v\in L'_n/pL'_n|Q'(v)\equiv m \mod p\}$, which is the same as counting the number of the solutions of $$x'y'+\epsilon {z'}^2+p{w'}^2-p\epsilon {u'}^2=m$$
in $\mathbb{F}_p$.  If we are in supersingular or general case, then Theorem~\ref{Thm:localdensity}(\ref{Thm:localdensity3}) tells that 
$$\delta(p,L'_n,m)\leq \left\{\begin{aligned}
   & 1-\frac{1}{p},\,\,&n\leq \frac{A}{p+1},\\
   & 0,\,\,&n>\frac{A}{p+1}.
\end{aligned}\right. $$
If we are in special but nonsupersingular case, then we only have
$$\delta(p,L'_n,m)\leq \left\{\begin{aligned}
   & 1-\frac{1}{p},\,\,&n\leq a,\\
   & 2,\,\,&n>a.
\end{aligned}\right. $$

We now use (\ref{eq:sdfsdf123}) to give a bound on $l_P(m)$. Proposition~\ref{T:DWformula} 
implies that for $n\leq \frac{A}{p+1}$ one has $q_G(m)=o(m^{\frac{3}{2}})$. Therefore, by Corollary~\ref{T:BKformula} and the fact that $|{L'}^{\vee}/L'|=p^2$, in supersingular or general case, we have 
\begin{equation*}
    l_P(m)\leq \frac{A(1-\frac{1}{p})}{p(p+1)(1+p^{-2})}|q_L(m)|+o(m^{\frac{3}{2}})=\frac{(p-1)^2}{p^2+1}g_P(m)+o(m^{\frac{3}{2}}).
\end{equation*}\\

Now we treat the special but nonsupersingular case. By Theorem~\ref{Thm:localdensity}(\ref{Thm:localdensity1}), there is an $N$ such that $\#\{s\in L''_n|Q'(s)=m\}=0$ when $n>N$. Therefore, instead of (\ref{eq:sdfsdf123}), we have a better bound 
\begin{equation}\label{eq:sdfsdf124}
    l_P(m)\leq \sum_{n=1}^{N}\#\left\{s\in L'_{n}|Q'(s)=m\right\}.  
\end{equation} %Denote by $L_{0,1}, L_{r,1} $ and $ L_{r,2}$ the lattices $L''_{1}, L''_{A_{r-1}+1} $ and $ L''_{A_{r-1}+ap^{2r}+1}$. 
Let $L_{0,1}',L_{r,1}' $ and $L_{r,2}'$ be the lattices $L', L'_{A_{r-1}+1}$ and $L'_{A_{r-1}+ap^{2r}+1}$, and let $t_{r,i}(m)=\#\{s\in L'_{r,i}|Q'(s)=m\}$. Using (\ref{eq:sdfsdf124}), we get an estimate 
\begin{equation}\label{eq:11}
        l_P(m) \leq  at_{0,1}(m)+(A-a)t_{0,2}(m)+\sum_{r=1}^{r_0}(at_{r,1}(m)+(A-a)t_{r,2}(m))p^{2r},
\end{equation}
where $r_0$ is a constant that depends only on $N$. Let $E_{r,i}, G_{r,i}$ be the Eisenstein series and cusp forms attached to $L_{r,i}'$ (\S\ref{sub:EsQl}), and let$$E=aE_{0,1}(m)+(A-a)E_{0,2}(m)+\sum_{r=1}^{r_0}(aE_{r,1}(m)+(A-a)E_{r,2}(m))p^{2r},$$  
$$G=aG_{0,1}(m)+(A-a)G_{0,2}(m)+\sum_{r=1}^{r_0}(aG_{r,1}(m)+(A-a)G_{r,2}(m))p^{2r}.$$
Let $q_E(m)$ and $q_G(m)$ be the Fourier coefficient of $E$ and $G$. By Proposition~\ref{T:DWformula}, for $r\leq r_0$ we have  $q_{G_{r,i}}(m)=o(m^{3/2})$, so $q_G(m)=o(m^{3/2})$. For the Eisenstein series, we will apply Corollary~\ref{T:BKformula}. Since the local densities $\delta(p,L'_n,m)$ are already computed, it suffices to figure out $|{L'_{r,i}}^{\vee}/L'_{r,i}|$.
Since $|{L'}^{\vee}/L'|=p^2$, we have $|{L'_{r,i}}^{\vee}/L'_{r,i}|=(p[L':L'_{r,i}])^2$. By Theorem~\ref{Thm:localdensity}(\ref{Thm:localdensity2}), we have $[L':L'_{r,1}]\geq p^{4r}$ and $[L_{r,1}':L'_{r,2}]\geq p^2$. It follows that \begin{align*}
   \frac{q_E(m)}{|q_L(m)|}& \leq   
    \frac{a(1-\frac{1}{p})}{p(1+p^{-2})}+\frac{2(A-a)}{p(p^2+1)}+2\sum_{r=1}^{r_0} \frac{ap^{2r}}{p^{4r-1}(p^2+1)}+\frac{(A-a)p^{2r}}{p^{4r+1}(p^2+1)}\\
    &<\alpha_P\left( \frac{a(p-1)}{p^{2}+1}+\frac{2(A-a)}{p(p^2+1)}+2\sum_{r=1}^{\infty} \frac{a}{p^{2r-1}(p^2+1)}+\frac{A-a}{p^{2r+1}(p^2+1)}\right)\\
    &=  \frac{\alpha_P A}{p^2-1}
\end{align*}
for some $0<\alpha_P<1$. As a result, we have $l_P(m)\leq  \alpha_Pg_P(m)+o(m^{\frac{3}{2}})$. 

Finally we vary the superspecial point $P$. Since there are only finitely many superspecial points, there exists an absolute constant $0<\alpha<1$ such that $l_P(m)\leq  \alpha g_P(m)+o(m^{\frac{3}{2}})$ for all superspecial $P$. By Lemma~\ref{T:intspecial}, the only supersingular intersections of $C$ with $\overline{Z(m)}$ are superspecial. Summing up over all superspecial points and use Proposition~\ref{T:intersection}, we obtain the desired bound. $\hfill\square$ \\

\begin{remark}\label{lastrmk}
If there exists a special endomorphism $\omega\in L''\otimes \mathbb{Z}_p$ that never decays, then it is possible that $\delta(p,L'_n,m)=2$ for all large $n$, so one only gets
$$\frac{q_E(m)}{|q_L(m)|} \leq   
    \frac{a(1-\frac{1}{p})}{p(1+p^{-2})}+\frac{2(A-a)}{p(p^2+1)}+2\sum_{r=1}^{\infty} \frac{ap^{2r}}{p^{4r-1}(p^2+1)}+\frac{(A-a)p^{2r}}{p^{4r+1}(p^2+1)}=\frac{A}{p^2-1},$$
which only gives that $l_P(m)\leq g_P(m)+o(m^{\frac{3}{2}})$. In this case, one is unable to deduce the unboundedness of the number of intersection points. Of course, the failure of the method in this particular situation meets our expectation, as discussed in the introduction.
\end{remark}

\section{Proof of the main theorems}\label{reduction}
 Let $C$ be a smooth connected projective curve with a nonconstant map $\tau:C\rightarrow \mathcal{A}^{\tor}_{2,\kk}$ sending the generic point to the almost ordinary stratum. The main goal of the whole section is to prove the following
\begin{theorem}\label{T:interMain}
Suppose that the image of $C$ does not lie on any special divisor. Let $\Delta$ be an infinite set of positive integers $m$ such that  $\left(\frac{m}{p}\right)=1$, then there are infinitely many $\kk$-points of $C$ lying on $\bigcup_{m\in \Delta} Z(m)$.
\end{theorem}
The proof of Theorem~\ref{T:interMain} will occupy the rest of the paper. We first deduce two main consequence from it:
\begin{theorem}
\label{T:maintheorem}
Suppose $U$ is a connected smooth curve with a nonconstant morphism $U\rightarrow \mathcal{A}_{2,\kk}$ whose image lies generically in the almost ordinary Newton stratum. Let $l$ be a positive integer such that $\left(\frac{l}{p}\right)=1$. If the image of $U$ does not lie on any special divisor, then:
\begin{enumerate}
    \item\label{T:maintheorem1} there are infinitely many $\kk$-points on $U$ admitting real multiplication by $\mathbb{Q}[x]/(x^2-l)$.
    \item \label{T:maintheorem2}there are infinitely many $\kk$-points on $U$ which is non-simple. 
\end{enumerate}
\end{theorem}
\proof Compactify the map $U\rightarrow \mathcal{A}_{2,\kk}$ to $C\rightarrow \mathcal{A}^{\tor}_{2,\kk}$ and let $\Delta=\{lk^2|k\in\mathbb{Z}\}$ in Theorem~\ref{T:interMain}. We obtain (\ref{T:maintheorem1}) by Lemma~\ref{T:non-simple}. (\ref{T:maintheorem2}) is a special case of (\ref{T:maintheorem1}) by letting $l=1$. $\hfill\square$
\begin{theorem}
\label{T:unlikelyint}
Suppose $U$ is a connected smooth (affine) curve with a nonconstant map $U\rightarrow \mathcal{A}_{2,\kk}$ whose image lies generically in the almost ordinary locus, and not in any special divisors. Let $S\subseteq U(\kk)$ be a finite set containing all the points mapping to the supersingular loci, then: 
\begin{enumerate}
    \item \label{T:unlikelyint1}all special divisors $Z(m)$ such that $\left(\frac{m}{p}\right)=-1$ are disjoint from $U-S$,
    \item\label{T:unlikelyint2} only finitely many special divisors $Z(m)$ such that $\left(\frac{m}{p}\right)=1$ are disjoint from $U-S$.
\end{enumerate}
\end{theorem}
\proof  
(\ref{T:unlikelyint1}) follows from the fact that when  $\left(\frac{m}{p}\right)=-1$, $Z(m)$ does not admit almost ordinary point, hence all intersections of $Z(m)$ with $U$ lies on the supersingular loci. (\ref{T:unlikelyint2}) follows by taking compatification and use Theorem~\ref{T:interMain}. $\hfill\square$\\[10pt]
Before proving Theorem~\ref{T:interMain}, we need several lemmas:
\begin{lemma}\label{lm:neverdecayimpliesZ(m)}
Let $C$ be as above. Let $P\in C(\kk)$ be an almost ordinary point or a bad reduction point. Suppose that $L''_P$ contains an element $s\neq 0$ such that $s\in L_{n,P}''\otimes \bZ_p$ for any $n>0$. Then the image of $C$ lies on a special divisor $Z(m)$.
\end{lemma}
\begin{proof}
Let $t$ be a uniformizer of $C$ at $P$. Let $\widehat{K}(C)=\kk((t))$ be the completion of the function field $K(C)$ at $P$.

If $P$ is almost ordinary, then the condition implies that $s$ lifts to a special endomorphism in $\mathscr{A}_{C^{/P}}[p^{\infty}]$ by \cite[Lemma 2.4.4]{DJ95}, hence $\mathscr{A}_{\widehat{K}(C)}[p^{\infty}]$ admits a special endomorphism. Theorem~\ref{T: algebraicC} then implies that $C$ lies on some $Z(m)$.

Now suppose that $P$ is of bad reduction. Let $\mathfrak{L}$ be the log 1-motive introduced in \S\ref{subsub:log1motive}. By Lemma~\ref{lm:logpullbackfs}, the inverse image log structures on $\Spec \kk[[t]]$ and $\Spec \kk[[t]]/t^n$ are all fs and sharp, so we can talk about the pullback of $\mathfrak{L}$ and its $p$-divisible group over these rings. The following is a logarithmic analogue of \cite[Lemma 2.4.4]{DJ95}:\\ 

\begin{claim}\label{cl:1}
$s$ lifts to a special endomorphism of the log $p$-divisible group $\mathfrak{L}_{\Spec \kk[[t]]}[p^\infty]$.\\
\end{claim}

Once the claim is proved, we can restrict $\mathfrak{L}_{\Spec \kk[[t]]}[p^\infty]$ to $\widehat{K}(C)$. Since $\widehat{K}(C)$ has trivial log structure, the restriction is nothing other than $\mathscr{A}_{\widehat{K}(C)}[p^\infty]$ (more precisely, the restriction is the 1-motive corresponding to $\mathscr{A}_{\widehat{K}(C)}$, whose $p$-divisible group coincides with the $p$-divisible group of $\mathscr{A}_{\widehat{K}(C)}$). We deduce that $\mathscr{A}_{\widehat{K}(C)}[p^{\infty}]$ admits a special endomorphism. Theorem~\ref{T: algebraicC} then implies that $C$ lies on some $Z(m)$.

It remains to prove the claim. Let's use the notation $ \Spec \kk[[t]]/t^{\infty}:=\Spec \kk[[t]]$. For $m\leq n\leq \infty$, let $\iota_{m,n}:\Spec \kk[[t]]/t^{m} \hookrightarrow  \Spec \kk[[t]]/t^{n}$ be closed immersions of fs log schemes. The conditions in the lemma imply that there is a compatible collection $(s_n)_{n\in \bN}$ such that (a) $s_1=s$ and (b) for $n<\infty$, each $s_n$ is a special endomorphism of $\mathfrak{L}_{\Spec \kk[[t]]/t^n}[p^\infty]$ and (c) for $m\leq n<\infty$, $\iota_{m,n}^* s_{n}=s_m$. 

Fix one $i>0$. For $n\leq \infty$, let $G_{n}=\mathfrak{L}_{\Spec \kk[[t]]/t^n}[p^\infty]$, and let $G_{n}[p^i]=(\Spec B_{i,n},M_{i,n},\alpha_{i,n})$ (i.e., a log scheme). Since $G_{n}[p^i]= \iota_{n,\infty}^*G_\infty[p^i]$, we have $$(\Spec B_{i,n},M_{i,n},\alpha_{i,n})= (\Spec B_{i,\infty}/t^n,\iota_{n,\infty}^*(M_{i,\infty},\alpha_{i,\infty})).$$ 
By \cite[Proposition 1.4]{logDiu}, $B_{i,\infty}$ is finite over $\kk[[t]]$. This has two consequences: 
\begin{enumerate}
 \item\label{it:qqq1} It implies that $B_{i,\infty}$ is $t$-adically complete, so $ B_{i,\infty}= \varprojlim_n B_{i,n}$ as algebras. 
    \item\label{it:qqq2} Since $\kk[[t]]$ is a strictly Henselian local ring, so by \cite[Tag 04GH]{stacks-project}, all $B_{i,n}$'s are products of strict Henselian local rings. Furthermore, $\dim B_{i,n}=0$ for all $n<\infty$, while $\dim B_{i,\infty}\leq 1$.  
\end{enumerate}
The collection $(s_n)_{n\in \bN}$ gives rise to a compatible family of morphisms of log schemes 
\begin{align*}
    s_{i,n}&:\Spec B_{i,n}\rightarrow \Spec B_{i,n},\\
    \theta_{i,n}&: s_{i,n}^{-1}M_{i,n}\rightarrow M_{i,n}.
\end{align*}
Taking direct limit of $s_{i,n}$, this immediately gives rise to a morphism of schemes ${s}_{i,\infty}: \Spec B_{i,\infty}\rightarrow \Spec B_{i,\infty}$. We need to show that it upgrades to a morphism of log schemes. In concrete terms, it means that one need to construct a map $\theta_{i,\infty}:{s}_{i,\infty}^{-1}M_{i,\infty}\rightarrow M_{i,\infty}$ that is compatible with ${s}_{i,\infty}$ and various $s_{i,n},\theta_{i,n}$ for $n<\infty$.

One need \textit{a priori} construct $\theta_{i,\infty}$ as a morphism of étale sheaves. But we can make several simplifications. By (\ref{it:qqq2}), the étale site over $\Spec B_{i,\infty}$ contains no more information than the Zariski site. So it suffices to construct $\theta_{i,\infty}$ in the Zariski site. Now the Zariski topology for $\Spec B_{i,\infty}$ is extremely simple: there are essentially two open sets $\Spec B_{i,\infty}$ (i.e. the ambient space) and $\Spec B_{i,\infty}[\frac{1}{t}]$ (i.e. the generic points). But $\Spec B_{i,\infty}[\frac{1}{t}]$ has trivial log structure (since $\Spec \kk[[t]]$ has trivial log structure over the generic point). Therefore, to construct $\theta_{i,\infty}$, it suffices to work with global sections of the relevant sheaves. By abuse of notation, in the following we will use $M_{i,n}$ to denote the global section of the corresponding sheaf on $\Spec B_{i,n}$, and use $\alpha_{i,n}$ to denote the map $M_{i,n}\rightarrow B_{i,n}$. For $n<\infty$, we have a compatible family of commuting squares 
\begin{equation}\label{eq:nlessinftycommute}
    \begin{tikzcd}
{M_{i,n}} \arrow[r, "\theta_{i,n}"] \arrow[d, "{\alpha_{i,n}}"'] & {M_{i,n}} \arrow[d, "{\alpha_{i,n}}"] \\
{B_{i,n}} \arrow[r, "{{s}_{i,n}}"']              & {B_{i,n}}                                 
\end{tikzcd}
\end{equation}
It suffices to construct a map $\theta_{i,\infty}: M_{i,\infty}\rightarrow M_{i,\infty}$ that fits in the square (\ref{eq:nlessinftycommute}) with $n=\infty$, compatible with all squares with $n<\infty$. 

Note that $B_{i,\infty}^*\rightarrow B^*_{i,n}$ is surjective, with kernel 
$K_{i,n}=1+(t^{n})B_{i,\infty}$. Embed $B_{i,\infty}^*$ as a submonoid of $M_{i,\infty}$ via the isomorphism $\alpha_{i,n}:\alpha_{i,n}^{-1}(B_{i,\infty}^*)\xrightarrow{\sim} B_{i,\infty}^*$. Then each $K_{i,n}$ is a submonoid of $M_{i,\infty}$. The composition of maps $M_{i,\infty}\xrightarrow{\alpha_{i,\infty}} B_{i,\infty}\rightarrow B_{i,n}$ factors through $K_{i,n}$, and identifies 
$\alpha_{i,n}:M_{i,n}\rightarrow B_{i,n}$ with 
$M_{i,\infty}/K_{i,n}\rightarrow B_{i,n}$. It then follows that $M_{i,\infty}=\varprojlim_n M_{i,\infty}/K_{i,n}=\varprojlim_n M_{i,n}$. So we can take limit of $\theta_{i,n}$ to get a $\theta_{i,\infty}: M_{i,\infty}\rightarrow M_{i,\infty}$, which satisfies the desired functorial property. 

This finishes the proof that $(s_{i,\infty},\theta_{i,\infty})$ is a morphism of log schemes. To ease notation, we will just denote the morphism by $s_{i,\infty}$. It is easily checked that $s_{i,\infty}$ is an endomorphism of the log finite group scheme $G_\infty[p^i]$. Let $i$ vary, we have a compatible collection $(s_{i,\infty})_{i\in \bN^+}$, which defines an endomorphism $s_{\infty}$ of the log $p$-divisible group $G_{\infty}$. This is the special endomorphism that we want and finishes the proof of the claim.\end{proof}

\begin{proposition}
\label{T:intindep} Let $C$ be as above. Assume that the image of $C$ does not lie on any special divisor. Let $P\in C(\kk)$ be a point whose image is an almost ordinary point or a bad reduction point. Then there exists a number $N_P$ only depending on $P$ such that for any $m$ coprime to $p$, $(C\cdot \overline{Z(m)})_P\leq N_P$. 
\end{proposition}
\proof Let $n$ denote a positive integer. According to Lemma~\ref{T:Rank0} (in the almost ordinary case) and Lemma~\ref{lm:boundaryrank} (in the bad reduction case), the quadratic form $Q'$ over $L''_{n,P}$ is represented by $a_{n}x^2$ for some integer $a_{n}$. Note that $a_{n}|a_{n+1}$. 

We claim that there exists $n$ such that $p|a_{n}$. Suppose this was not the case, then for any $n$, we have  $L_{n,P}''\otimes \bZ_p= L''_P\otimes \bZ_p$. Lemma~\ref{lm:neverdecayimpliesZ(m)} implies that $C$ lies on a special divisor $Z(m)$, contradicting our assumption. As a result, there exists a constant $N_P$  only depending on $P$, such that when $n>N_P/2$, the coefficient $a_{n}$ is divisible by $p$. In particular, this implies that when $n>N_P/2$ and $m$ is coprime to $p$, $\#\{x\in L''_{n,P} |Q'(x)=m \}=0$.
It follows from (\ref{eq:l_Platticecounting}) in the almost ordinary case, and Lemma~\ref{lm:specialboundarycount} in the bad reduction case, that for any $m$ coprime to $p$, we have $$l_P(m)\leq \sum_{n=1}^{[N_P/2]}\#\{x\in \bZ|a_{n}x^2=m \} \leq N_P.$$This proves the proposition.  $\hfill\square$\\[10pt]
\textit{Proof of Theorem~\ref{T:interMain}}. Let $\Delta$ be an infinite set containing positive integers $m$ such that $\left(\frac{m}{p}\right)=1$. Let $S_{\mathrm{ss}}$ (\textit{resp}. $S_{b}$, \textit{resp}. $S_{\mathrm{aos}}$) be the set of supersingular points (\textit{resp}. bad reduction points, \textit{resp} almost ordinary points that lie on $Z(m)$ for some $m\in \Delta$) on $C(\kk)$. Note that $S_{\mathrm{ss}}$ and $S_b$ are all finite. We need to show that $S_{\mathrm{aos}}$ is infinite. 

Suppose towards contradiction $|S_{\mathrm{aos}}|<\infty$. In the following, let $m\in \Delta$. Then by Proposition~\ref{T:intindep}, we have 
$$\sum_{P\in S_{\mathrm{aos}}\cup S_b} l_P(m) \leq B$$
for some constant $B$ that does not depend on $m$. On the other hand, by Lemma~\ref{decayresult}, we have 
$$\sum_{P\in S_{\mathrm{ss}}}l_P(m)\leq \alpha C\cdot \overline{Z(m)}+o(m^{\frac{3}{2}})$$  
for some $\alpha<1$ (recall that $C\cdot \overline{Z(m)}\sim m^{\frac{3}{2}}$ from Lemma~\ref{intersEisen}). It follows that $$C\cdot \overline{Z(m)} =\sum_{P\in S_{\mathrm{aos}}\cup S_b\cup S_{\mathrm{ss}}}l_P(m)\leq \alpha C\cdot \overline{Z(m)}+B+o(m^{\frac{3}{2}}).$$
Let $m\rightarrow \infty$, we arrive at a contradiction. $\hfill\square$
%Therefore we have 
%This implies that  such that $L''_{n,R}$
%According to \ref{T:Rank1}, the quadratic form $Q'$ restrict to the lattice of special endomorphisms of the semi-abelian scheme
%$\mathscr{A}_{\Spec \kk[[t_R]]/(t_R^n)}$ is represented by $a_{R,n}x^2$ for some integer $a_{R,n}$. \ref{T: algebraicC}  then 
%Since the lattice of special endomorphisms of
%$\mathscr{A}_{\Spec \kk[[t_R]]/(t_R^{n+1})}$ is contained in that of $\mathscr{A}_{\Spec \kk[[t_R]]/(t_R^{n})}$, we see that $a_{R,n}|a_{R,n+1}$. Since $\mathscr{A}_C$ does not admit any global special endomorphisms, we see that $\lim_{n\rightarrow \infty}a_{R,n}=\infty$. Therefore there exists a large number $N$, only depending on $C$ and $l$, and a subset $S_l\subseteq T_l$, such that:\\[5pt]
%(1) $\#S_l=\infty$. \\
%(2) An element $m\in S_l$ is represented as $m=a_{R,n}x^2$ for some $R\in \mathfrak{R}$ only when $n\leq N$.\\[5pt]

\bibliography{ref}

\begin{bibdiv}
\begin{biblist}

\bib{AGHMP17}{article}{
      author={Andreatta, Fabrizio},
      author={Goren, Eyal~Z.},
      author={Howard, Benjamin},
      author={Madapusi~Pera, Keerthi},
       title={Faltings heights of abelian varieties with complex multiplication},
        date={2018},
        ISSN={0003-486X},
     journal={Ann. of Math. (2)},
      volume={187},
      number={2},
       pages={391\ndash 531},
         url={https://doi-org.ezproxy.library.wisc.edu/10.4007/annals.2018.187.2.3},
}

\bib{BR08}{article}{
      author={Baker, Matthew},
      author={Ih, Su-ion},
      author={Rumely, Robert},
       title={A finiteness property of torsion points},
        date={2008},
        ISSN={1937-0652},
     journal={Algebra Number Theory},
      volume={2},
      number={2},
       pages={217\ndash 248},
         url={https://doi-org.ezproxy.library.wisc.edu/10.2140/ant.2008.2.217},
}

\bib{BK01}{article}{
      author={Bruinier, Jan},
      author={Kuss, Michael},
       title={Eisenstein series attached to lattices and modular forms on orthogonal groups},
        date={2000},
     journal={manuscripta mathematica},
      volume={106},
}

\bib{Box15}{book}{
      author={Boxer, George~Andrew},
       title={Torsion in the {C}oherent {C}ohomology of {S}himura {V}arieties and {G}alois {R}epresentations},
   publisher={ProQuest LLC, Ann Arbor, MI},
        date={2015},
        ISBN={978-1339-29190-1},
         url={http://gateway.proquest.com.ezproxy.library.wisc.edu/openurl?url_ver=Z39.88-2004&rft_val_fmt=info:ofi/fmt:kev:mtx:dissertation&res_dat=xri:pqm&rft_dat=xri:pqdiss:3738711},
        note={Thesis (Ph.D.)--Harvard University},
}

\bib{Bru}{article}{
      author={Bruinier, Jan~Hendrik},
       title={Borcherds products with prescribed divisor},
        date={2017},
     journal={Bulletin of the London Mathematical Society},
      volume={49},
      number={6},
       pages={979\ndash 987},
         url={https://doi.org/10.1112},
}

\bib{BHS19}{article}{
      author={Bruinier, Jan~Hendrik},
      author={Zemel, Shaul},
       title={Special cycles on toroidal compactifications of orthogonal {S}himura varieties},
        date={2022},
        ISSN={0025-5831},
     journal={Math. Ann.},
      volume={384},
      number={1-2},
       pages={309\ndash 371},
         url={https://doi-org.ezproxy.library.wisc.edu/10.1007/s00208-021-02271-x},
}

\bib{Cha03}{article}{
      author={Chai, C.},
       title={Families of ordinary \text{Abelian} varieties : canonical coordinates, $p$-adic monodromy, \text{Tate}-linear subvarieties and \text{Hecke} orbits},
        date={2003},
}

\bib{Cha97}{article}{
      author={Chavdarov, Nick},
       title={{The generic irreducibility of the numerator of the zeta function in a family of curves with large monodromy}},
        date={1997},
     journal={Duke Mathematical Journal},
      volume={87},
      number={1},
       pages={151 \ndash  180},
         url={https://doi.org/10.1215/S0012-7094-97-08707-X},
}

\bib{DJ95}{article}{
      author={de~Jong, A.~J.},
       title={Crystalline \text{Dieudonné} module theory via formal and rigid geometry},
        date={1995},
     journal={Publications Mathématiques de l'Institut des Hautes Études Scientifiques},
      volume={82},
      number={1},
       pages={5\ndash 96},
}

\bib{DJ98}{article}{
      author={de~Jong, A.~D.},
       title={\text{Homomorphisms of Barsotti--Tate groups and crystals in positive characteristic}},
        date={1998},
     journal={Inventiones mathematicae},
      volume={134},
       pages={301\ndash 333},
}

\bib{DJ99}{article}{
      author={{De Jong}, {A. J.}},
      author={Messing, W.},
       title={Crystalline \text{Dieudonné} theory over excellent schemes},
    language={English (US)},
        date={1999},
        ISSN={0037-9484},
     journal={Bulletin de la Societe Mathematique de France},
      volume={127},
      number={2},
       pages={333\ndash 348},
}

\bib{Fal99}{article}{
      author={Faltings, Gerd},
       title={Integral crystalline cohomology over very ramified valuation rings},
        date={1999},
     journal={Journal of the AMS, v.12, 117-144 (1999)},
      volume={12},
}

\bib{FC}{book}{
      author={Faltings, Gerd},
      author={Chai, Ching-Li},
       title={Degeneration of abelian varieties},
      series={Ergebnisse der Mathematik und ihrer Grenzgebiete (3) [Results in Mathematics and Related Areas (3)]},
   publisher={Springer-Verlag, Berlin},
        date={1990},
      volume={22},
        ISBN={3-540-52015-5},
         url={https://doi-org.ezproxy.library.wisc.edu/10.1007/978-3-662-02632-8},
        note={With an appendix by David Mumford},
}

\bib{EGAIV2}{article}{
      author={Grothendieck, A.},
       title={\'{E}l\'{e}ments de g\'{e}om\'{e}trie alg\'{e}brique. {IV}. \'{E}tude locale des sch\'{e}mas et des morphismes de sch\'{e}mas. {II}},
        date={1965},
        ISSN={0073-8301},
     journal={Inst. Hautes \'{E}tudes Sci. Publ. Math.},
      number={24},
       pages={231},
         url={http://www.numdam.org/item?id=PMIHES_1965__24__231_0},
}

\bib{Han04}{article}{
      author={Hanke, Jonathan},
       title={Local densities and explicit bounds for representability by a quadratic form},
        date={2004},
     journal={Duke Mathematical Journal},
      volume={124},
}

\bib{HBK17}{article}{
      author={Howard, Benjamin},
      author={Madapusi~Pera, Keerthi},
       title={Arithmetic of {B}orcherds products},
        date={2020},
        ISSN={0303-1179},
     journal={Ast\'{e}risque},
      number={421, Diviseurs arithm\'{e}tiques sur les vari\'{e}t\'{e}s orthogonales et unitaires de Shimura},
       pages={187\ndash 297},
         url={https://doi-org.ezproxy.library.wisc.edu/10.24033/ast},
}

\bib{HP17}{article}{
      author={Howard, Benjamin},
      author={Pappas, Georgios},
       title={\text{Rapoport--Zink} spaces for spinor groups},
        date={2017},
     journal={Compositio Mathematica},
      volume={153},
       pages={1050 \ndash  1118},
}

\bib{Kk89}{article}{
      author={Kato, Kazuya},
       title={Logarithmic structures of {F}ontaine-{I}llusie},
        date={1989},
       pages={191\ndash 224},
      review={\MR{1463703}},
}

\bib{logDiu}{article}{
      author={Kato, Kazuya},
       title={{L}ogarithmic {D}ieudonn\'e theory},
      eprint={https://arxiv.org/pdf/2306.13943},
}

\bib{KM09}{article}{
      author={Kisin, Mark},
       title={\text{Integral models for Shimura varieties of Abelian type}},
        date={2010},
        ISSN={08940347, 10886834},
     journal={Journal of the American Mathematical Society},
      volume={23},
      number={4},
       pages={967\ndash 1012},
         url={http://www.jstor.org/stable/20753577},
}

\bib{LogAV}{article}{
      author={Kajiwara, Takeshi},
      author={Kato, Kazuya},
      author={Nakayama, Chikara},
       title={Logarithmic abelian varieties. {II}. {A}lgebraic theory},
        date={2008},
     journal={Nagoya Math. J.},
      volume={189},
      number={1},
       pages={63\ndash 138},
}

\bib{LKW}{book}{
      author={Lan, Kai-Wen},
       title={Arithmetic compactifications of {PEL}-type {S}himura varieties},
      series={London Mathematical Society Monographs Series},
   publisher={Princeton University Press, Princeton, NJ},
        date={2013},
      volume={36},
        ISBN={978-0-691-15654-5},
         url={https://doi-org.ezproxy.library.wisc.edu/10.1515/9781400846016},
}

\bib{BMI98}{article}{
      author={Moonen, Ben},
       title={Linearity properties of {S}himura varieties. {I}},
        date={1998},
        ISSN={1056-3911},
     journal={J. Algebraic Geom.},
      volume={7},
      number={3},
       pages={539\ndash 567},
      review={\MR{1618140}},
}

\bib{Moo98}{incollection}{
      author={Moonen, Ben},
       title={Models of {S}himura varieties in mixed characteristics},
        date={1998},
   booktitle={Galois representations in arithmetic algebraic geometry ({D}urham, 1996)},
      series={London Math. Soc. Lecture Note Ser.},
      volume={254},
   publisher={Cambridge Univ. Press, Cambridge},
       pages={267\ndash 350},
         url={https://doi-org.ezproxy.library.wisc.edu/10.1017/CBO9780511662010.008},
}

\bib{MP08}{article}{
      author={Murty, V.~Kumar},
      author={Patankar, Vijay~M.},
       title={{Splitting of Abelian Varieties}},
        date={2008},
        ISSN={1073-7928},
     journal={International Mathematics Research Notices},
      volume={2008},
         url={https://doi.org/10.1093/imrn/rnn033},
}

\bib{MP11}{article}{
      author={Madapusi~Pera, Keerthi},
       title={Toroidal compactifications of integral models of {S}himura varieties of {H}odge type, {P}h{D} {T}hesis},
        date={2011},
}

\bib{MP19}{article}{
      author={Madapusi~Pera, Keerthi},
       title={Toroidal compactifications of integral models of {S}himura varieties of {H}odge type},
        date={2019},
        ISSN={0012-9593,1873-2151},
     journal={Ann. Sci. \'Ec. Norm. Sup\'er. (4)},
      volume={52},
      number={2},
       pages={393\ndash 514},
         url={https://doi.org/10.24033/asens.2391},
      review={\MR{3948111}},
}

\bib{MP15}{article}{
      author={Madapusi~Pera, Keerthi},
       title={Integral canonical models for spin \text{Shimura} varieties},
        date={2016},
        ISSN={0010-437X},
     journal={Compositio Mathematica},
      volume={152},
      number={4},
       pages={769\ndash 824},
}

\bib{MAT}{article}{
      author={Maulik, Davesh},
      author={Shankar, Ananth~N.},
      author={Tang, Yunqing},
       title={Reductions of \text{Abelian} surfaces over global function fields},
        date={2020},
     journal={Compositio Mathematica},
}

\bib{MST22}{article}{
      author={Maulik, Davesh},
      author={Shankar, Ananth},
      author={Tang, Yunqing},
       title={Picard ranks of \text{K3} surfaces over function fields and the \text{Hecke} orbit conjecture},
        date={2022},
     journal={Inventiones mathematicae},
      volume={228},
}

\bib{Ogus01}{incollection}{
      author={Ogus, Arthur},
       title={Singularities of the height strata in the moduli of {$K3$} surfaces},
        date={2001},
   booktitle={Moduli of abelian varieties ({T}exel {I}sland, 1999)},
      series={Progr. Math.},
      volume={195},
   publisher={Birkh\"{a}user, Basel},
       pages={325\ndash 343},
}

\bib{KR00}{article}{
      author={Rapoport, Michael},
      author={Kudla, Stephen},
       title={Cycles on \text{Siegel threefolds and derivatives of Eisenstein} series},
        date={2000},
     journal={Annales Scientifiques De L Ecole Normale Superieure - ANN SCI ECOLE NORM SUPER},
      volume={33},
       pages={695\ndash 756},
}

\bib{Shi00}{article}{
      author={Shiho, Atsushi},
       title={Crystalline fundamental groups. {I}. {I}socrystals on log crystalline site and log convergent site},
        date={2000},
        ISSN={1340-5705},
     journal={J. Math. Sci. Univ. Tokyo},
      volume={7},
      number={4},
       pages={509\ndash 656},
      review={\MR{1800845}},
}

\bib{ST20}{article}{
      author={Shankar, Ananth~N.},
      author={Tang, Yunqing},
       title={Exceptional splitting of reductions of \text{Abelian} surfaces},
        date={2020Feb},
        ISSN={0012-7094},
     journal={Duke Mathematical Journal},
      volume={169},
      number={3},
         url={http://dx.doi.org/10.1215/00127094-2019-0046},
}

\bib{stacks-project}{misc}{
      author={{Stacks Project Authors}, The},
       title={\textit{Stacks Project}},
         how={\url{https://stacks.math.columbia.edu}},
        date={2018},
}

\bib{ST22}{article}{
      author={Tayou, Salim},
       title={Picard rank jumps for \text{K3} surfaces with bad reduction},
        date={2022},
     journal={Algebra \& Number Theory},
         url={https://arxiv.org/abs/2203.09559},
}

\bib{Sheer}{article}{
      author={W\"urthen, Matti},
      author={Zhao, Heer},
       title={Log {$p$}-divisible groups associated with log 1-motives},
        date={2024},
        ISSN={0008-414X,1496-4279},
     journal={Canad. J. Math.},
      volume={76},
      number={3},
       pages={946\ndash 983},
         url={https://doi.org/10.4153/S0008414X23000287},
      review={\MR{4747297}},
}

\bib{Zin01}{article}{
      author={Zink, Thomas},
       title={{On the slope filtration}},
        date={2001},
     journal={Duke Mathematical Journal},
      volume={109},
      number={1},
       pages={79 \ndash  95},
         url={https://doi.org/10.1215/S0012-7094-01-10913-7},
}

\bib{Zyw14}{article}{
      author={Zywina, David},
       title={{The Splitting of Reductions of an Abelian Variety}},
        date={2013},
        ISSN={1073-7928},
     journal={International Mathematics Research Notices},
      volume={2014},
      number={18},
       pages={5042\ndash 5083},
         url={https://doi.org/10.1093/imrn/rnt113},
}

\end{biblist}
\end{bibdiv}
\bibliographystyle{numeric}
\end{document}